\documentclass[12pt]{amsart}

\usepackage{DKstyle}
\usepackage{bm}
\usepackage{float}

\newcommand{\Mod}[1]{\ (\mathrm{mod}\ #1)}
\newcommand{\vertiii}[1]{{\left\vert\kern-0.25ex\left\vert\kern-0.25ex\left\vert #1
    \right\vert\kern-0.25ex\right\vert\kern-0.25ex\right\vert}}
\theoremstyle{plain}
\newtheorem*{theorem*}{Theorem}
\makeatletter
\newcommand*{\rom}[1]{\expandafter\@slowromancap\romannumeral #1@}
\makeatother
\newcommand{\fg}{\frak{g}}

\newcommand{\bsm}{\left(\begin{smallmatrix}}
\newcommand{\esm}{\end{smallmatrix}\right)}

\usepackage{caption}
\usepackage[labelfont=rm]{subcaption}

\usepackage[pagebackref=true, colorlinks]{hyperref}

\hypersetup{pdffitwindow=true,linkcolor=blue,citecolor=blue,urlcolor=blue,menucolor=blue}

\usepackage{comment}


\date{\today}%
\dedicatory{}%
\commby{}%

\title{Translates of rational points along expanding closed horocycles on the modular surface}
\author{Claire Burrin}
\address{Department of Mathematics, Rutgers University, 110 Frelinghuysen Rd, Piscataway, NJ 08854}
\email{claire.burrin@rutgers.edu}
\address{Current address: Department Mathematik, ETH Z\"urich, 8092 Zurich, Switzerland}
\email{claire.burrin@math.ethz.ch}

\author{Uri Shapira}
\address{Department of Mathematics, Technion, Haifa, Israel}
\email{ushapira@tx.technion.ac.il}
\author{Shucheng Yu}
\address{Department of Mathematics, Technion, Haifa, Israel}
\email{yushucheng@campus.technion.ac.il}
\address{Current address: Department of Mathematics, Uppsala University, Box 480, SE-75106, Uppsala, SWEDEN}
\email{shucheng.yu@math.uu.se}
\thanks{U.S. and S.Y. acknowledge the support of ISF grant number  871/17. U.S. and S.Y. acknowledge that this project has received funding from the European Research Council (ERC) under the European Union's Horizon 2020 research and innovation program (grant agreement No.\ 754475).}

\begin{document}

\begin{abstract}
We study the limiting distribution of the rational points under a horizontal translation along a sequence of expanding closed horocycles on the modular surface. Using spectral methods we confirm equidistribution of these sample points for any translate when the sequence of horocycles expands within a certain polynomial range. We show that the equidistribution fails for generic translates and a slightly faster expanding rate. We also prove both equidistribution and non-equidistribution results by obtaining explicit limiting measures while allowing the sequence of horocycles to expand arbitrarily fast. Similar results are also obtained for translates of primitive rational points.
\end{abstract}
\maketitle

\tableofcontents

\section{Introduction}

Let $\{S_n\}_{n\in\N}$ be a sequence of ``nice'' subsets that become equidistributed in their ambient space. Given a sequence of discrete subsets $\{R_n\}_{n\in\N}$ with $R_n\subset S_n$, an interesting question is to study to what extent does the distribution behavior of $\{R_n\}_{n\in \N}$ mimic that of $\{S_n\}_{n\in\N}$. One naturally expects that when the size of $R_n$ is relatively large, it is more likely that $\{R_n\}_{n\in\N}$ inherits some distribution property from $\{S_n\}_{n\in\N}$; on the other hand if $R_n$ lies on $S_n$ sparsely, then it is more likely that points in $\{R_n\}_{n\in\N}$ become decorrelated and distribute like random points on the ambient space.

In the setting of unipotent dynamics, the most typical example of a sequence $\{S_n\}_{n\in\N}$ is a sequence of \textit{expanding closed horocycles} on a non-compact finite-area hyperbolic surface $\cM$. More precisely, we can realize $\cM$ as a quotient $\G\bk\bH$ where $\G$ is a co-finite Fuchsian subgroup and $\bH=\{z=x+iy\in\C:y>0\}$ is the Poincar\'e upper half-plane, equipped with the hyperbolic metric $ds=|dz|/y$, where $dz=dx+idy$ is the complex line element. Up to conjugating by an appropriate isometry, we may assume that $\cM=\G\bk\bH$ has a width one cusp at infinity, that is, that the isotropy group $\G_{\infty}< \G$ is generated by the translation sending $z\in\bH$ to $z+1$. A \textit{closed horocycle of height $y>0$} is a closed set of the form 
$$
\mathcal{H}_y:=\{\Gamma(x+iy): x\in\R/\Z\} \subset \cM,
$$ 
and its period, i.e., its hyperbolic length, is $y^{-1}$. As $\mathcal{H}_y$ gets longer, that is, as $y\to0^+$, it becomes equidistributed on $\cM$ with respect to the hyperbolic area $d\mu(z)=y^{-2}dxdy$. The first effective version of this result is due to Sarnak \cite{Sarnak1981} who, using spectral arguments, proved that for every $\Psi\in C_c^\infty(\Gamma\backslash\bH)$ and any $y>0$,
\begin{equation}\label{low-equid}
\int_0^1 \Psi(x+iy) dx =\frac{\int_{\cM}\Psi(z)d\mu(z)}{\mu(\cM)} +O\left(\mathcal{S}(\Psi) y^{\alpha}\right),
\end{equation}
where $\mathcal{S}$ is some Sobolev norm, and $0<\alpha<1$ is a constant depending on the first non-trivial residual hyperbolic Laplacian eigenvalue of $\Gamma$. %see \propref{prop:equiclohor} for a special case of \eqref{low-equid} when $\G$ is a congruence subgroup with an explicit Sobolev norm. %where the Sobolev norm is much more explicit. %Here and thereafter for any measure $\nu$ on $\cM$ we denote by $\nu(\Psi):=\int_{\cM}\Psi(z)d\nu(z)$. 
In the case of the modular surface $\SL_2(\Z)\backslash\bH$, $\alpha=\frac12$, while Zagier \cite{Zagier1981} observed that the Riemann hypothesis is equivalent to the equidistribution rate $O_{\e}\left(y^{3/4-\e}\right)$. 

In this setting, this problem was first investigated by Hejhal in \cite{Hejhal1996} with a heuristic and numerical study of the value distribution of the sample points
\begin{equation}\label{equ:samplehe}
\G\left(\tfrac{x+j}{n}+iy\right) : 0\leq j\leq n-1
\end{equation}
for some Hecke triangle groups $\Gamma=\mathbb{G}_q$ under the assumption that $ny$ is small. Set
%\begin{equation}\label{equ:hejhal}
$$S_{y,n,\Psi}(x):=\sum_{j=0}^{n-1}\Psi\left(\tfrac{x+j}{n}+iy\right),$$
%\end{equation}
where $\Psi$ is some mean-zero step function on a fixed fundamental domain for $\G\bk\bH$ (automorphically extended to $\bH$). The numerics show that the value distribution of $n^{-1/2}S_{n,y,\Psi}(x)$ with respect to $x\in [0,1)$
 approaches a Gaussian curve for the non-arithmetic Hecke triangle groups $\mathbb{G}_5$ and $\mathbb{G}_7$, while this phenomenon breaks down for $\mathbb{G}_3=\PSL_2(\Z)$.  
 Hejhal gave an explanation of this difference based on the existence of Hecke operators on $\mathbb{G}_3$. %see \rmkref{rmk:hejhal} for a more explicit form. 
The convergence to a Gaussian distribution for general non-arithmetic Fuchsian groups was later confirmed by Str\"ombergsson \cite[Corollary 6.5]{Strombergsson2004}, under the assumption that the sequence $\{y_n\}_{n\in\N}$ decays sufficiently rapidly.
%$$S_{y,n}(x):=\sum_{j=0}^{n-1}F(\tfrac{x+j}{n}+iy),$$
%for some Hecke triangle groups. 
%Here $F$ is a test function on $\cM$ and $n\to\infty$, $y\to 0$ in a way such that $ny\to 0$.

Other such problems have since been investigated. Marklof and Str\"ombergsson \cite{MarklofStrombergsson2003} proved the equidistribution of generic Kronecker sequences
\begin{equation}\label{equ:kroneseq}
\{\G(j\beta +iy_n)\in \cM:1\leq j\leq n\} \subset \cM
\end{equation}
 along a sequence of closed horocycles expanded at a certain rate $y_n$ on $T_1\cM$, the unit tangent bundle of $\cM$. The equidistribution of Hecke points proved by Clozel--Ullmo \cite{ClozelUllmo2004} (see also \cite{GoldsteinMayer2003}, \cite{ClozelOhUllmo2001}) implies the equidistribution of the primitive rational points 
$$\left\{\G\left(\tfrac{j}{n}+\tfrac{i}{n}\right) : 1\leq j\leq n-1,\ \gcd(j,n)=1\right\}$$ 
at prime steps on the modular surface, see \cite[Remark on p. 171]{GoldsteinMayer2003}. More recently, the equidistribution of the above sequence along the full sequence of positive integers was proved by Einsiedler--Luethi--Shah \cite{EinsiedlerLuethiShah2020} in a slightly more general setting, namely on the product of the unit tangent bundle of the modular surface and a torus. Various sparse equidistribution results have also been obtained for expanding horospheres in the space of lattices $\SL_n(\R)/\SL_n(\Z)$ for $n\geq3$ \cite{Marklof2010,Han2015,EinsiedlerMozesShahShapira2016,LeeMarklof2018,BazHuangLee2019} and in Hilbert modular surfaces \cite{Luethi2021}.
 %which stand as higher dimensional analogues of closed horocycles\footnote{ In this setting, hororspheres naturally embed in a proper sub-manifold of the ambient homogeneous space, and the equidistribution results mentioned here are with respect to this sub-manifold.}.  Marklof \cite{Marklof2010} proved the equidistribution of Farey sequences on expanding horospheres, see also \cite{Han2015} for an effective result. The equdistribution of primitive rational points on horospheres was proved by Einsiedler--Mozes--Shah--Shapira \cite{EinsiedlerMozesShahShapira2016}, and effective analogues were obtained in the recent work of Marklof--Lee \cite{LeeMarklof2018} and El-Baz--Huang--Lee \cite{BazHuangLee2019}. In another direction, Luethi \cite{Luethi2019} proved an effective equidistribution result for primitive rational points on expanding horospheres on Hilbert modular surfaces. 
For each of these equidistribution results, assumptions on the expanding rate of the sequence $\{S_n\}_{n\in\N}$ are crucial; the discrete subsets $\{R_n\}_{n\in\N}$ lying on $\{S_n\}_{n\in\N}$ can not be too sparse. %In the case of closed horocycles, the expanding rate corresponds to the decaying rate of heights of these horocycles.

%We also mention here another well-studied sparse equidistribution problem in the setting of unipotent dynamics (raised by Shah\cite{Shah1994}), that is, to study the distribution behavior of a sparse subset of a fixed dense unipotent orbit, see e.g. \cite{Venkatesh2010,Zheng2016,SarnakUbis2015,McAdam2019} for a non-exhaustive list. 
%In this paper, we consider the {\em sparse equidistribution problem} for the subset of $n$ evenly spaced points along the horocycle $\cH_{y_n}$ on the modular surface.
% {{These sample points, which we denote by $\cR_{n}(x, y_n)$ (cf. \eqref{equ:samplepoints} below), can be parameterized by the set of rational points on $\cH_{y_n}$ (with denominator $n$) under a horizontal translation $x\in\R/\Z$. }
{{
This paper emerged from an attempt to prove a result which turned out to be false. We consider the {\em sparse equidistribution problem} for the subset of rational points (with denominator $n$) under a horizontal translation $x\in\R/\Z$ on a horocycle $\cH_{y}$ on the modular surface; 
we denote this subset by $\cR_{n}(x,y_n)$ (cf. \eqref{equ:samplepoints}). We thought that since the closed horocycles $\cH_y$ equidistribute as $y\to 0^+$, if we fix a sequence $\{y_n\}_{n\in\N}$ approaching zero, then the normalized counting measures on $\cR_n(x,y_n)$ (and its primitive counterpart) should equidistribute for Lebesgue almost every $x$ as $n\to\infty$. See the recent paper of Bersudsky \cite[Theorem 1.5]{Bersudsky2020} for an analogue situation where such a result is true. Note the order of quantifiers; we first fix the sequence $\{y_n\}_{n\in\N}$ and only then choose the horizontal translation $x$. It is not hard to see that if one flips the quantifiers, for any fixed horizontal translation $x$, there are sequences $\{y_n\}_{n\in\N}$ (approaching zero rapidly) such that equidistribution fails. We were very surprised to learn though, that in stark contrast to our initial expectation, equidistribution fails. The main novel result of this paper (\thmref{thm:mainthm}) says that there are sequences $\{y_n\}_{n\in\N}$ approaching zero arbitrarily fast such that for almost every horizontal translation $x$ the normalized counting measures $\cR_n(x,y_n)$ and its primitive counterpart do not equidistribute. In fact, we show the collection of limit measures contains the uniform measure $\mu_{\cM}$, the zero measure and certain singular measures. Although these should be considered as the main contribution of this paper, we also complement our analysis with answering natural questions concerning sequences $\{y_n\}_{n\in\N}$ approaching zero in a polynomial rate.
}} 
The next subsections describe more precisely the setting and results obtained.

\subsection{Context of the present paper}\label{sec:context}
Let $\G=\SL_2(\Z)$ and let $\cM=\G\bk\bH$ be the modular surface. In this paper, generalizing the setting of \cite{EinsiedlerLuethiShah2020}, we study the equidistribution problem for the sets of rational and primitive rational points under an arbitrary horizontal translation $x\in\R/\Z$ along a given sequence of expanding closed horocycles on $\cM$. The set of rational points is the obvious choice of a sparse set with identical spacings, while primitive rational points constitute the simplest pseudorandom sequence (via the linear congruential generator). 
%by studying the \textit{limiting distribution} of \textit{translates} of (primitive) rational points on expanding closed horocycles on the modular surface. 
For any $n\in\N$, $x\in\R/\Z$ and $y>0$ we denote by
\begin{equation}\label{equ:samplepoints}
\cR_n(x,y):=\left\{\G(x+\tfrac{j}{n}+iy)\in\cH_y:0\leq j\leq n-1\right\} %\subset \cH_y
\end{equation}
and respectively
\begin{equation}\label{equ:prisamplepoints}
\cR_n^{\rm pr}(x,y):=\left\{\G(x+\tfrac{j}{n}+iy)\in\cH_y:j\in (\Z/n\Z)^{\times}\right\},%\subset \cH_y,
\end{equation}
the set of rational and respectively primitive rational points with denominator $n$ on the closed horocycle $\cH_y$ translated to the right by $x$. As usual, $(\Z/n\Z)^{\times}$ denotes here the multiplicative group of integers modulo $n$. %Note that $S_n(x,y)$ (respectively $S_n^{\rm pr}(x,y)$) is the set of (respectively primitive) rational points with denominator $n$ on the closed horocycle $\cH_y$ translated to the right by $x$, 
%and we call $x\in\R/\Z$ the \textit{the initial translate} of $S_n(x,y)$ and $S_n^{\rm pr}(x,y)$. 

Let $\{y_n\}_{n\in\N}$ be a sequence of positive numbers such that $y_n\to0$ as $n\to\infty$. We investigate the \textit{limiting distribution} of the sequences of sample points $\left\{\cR_n(x,y_n)\right\}_{n\in\N}$ and $\left\{\cR^{\rm pr}_n(x,y_n)\right\}_{n\in\N}$ under various assumptions on the expanding rate of the sequence of horocycles $\{\cH_{y_n}\}_{n\in \N}$, or equivalently, the decay rate of $\{y_n\}_{n\in\N}$. 

This problem is naturally easier when the sequence $\{y_n\}_{n\in\N}$ decays slowly since then at each step we have relatively more sample points on the underlying horocycle. For instance, if $ny_n\to\infty$ as $n\to\infty$, the hyperbolic distance between two adjacent points in $\cR_n(x,y_n)$ decays to zero as $n\to\infty$. Since the points in $\cR_n(x,y_n)$ distribute evenly on $\cH_{y_n}$, the distribution behavior of $\cR_n(x,y_n)$ then mimics that of $\cH_{y_n}$. In particular, for any $x\in\R/\Z$ the sequence $\left\{\cR_n(x,y_n)\right\}_{n\in\N}$ becomes equidistributed on $\cM$ with respect to the hyperbolic area $\mu$ as $n\to\infty$, following from the equidistribution of the sequence $\{\cH_{y_n}\}_{n\in\N}$. 

Regarding $\left\{\cR^{\rm pr}_n(x,y_n)\right\}_{n\in\N}$, its distribution behavior 
%On the other hand, the distribution behavior of the sample points $\left\{\cR^{\rm pr}_n(x,y_n)\right\}_{n\in\N}$ 
is well understood when $x=0$. Indeed, it was shown by Luethi \cite{Luethi2021} that if $y_n=c/n^{\alpha}$ for some $c>0$ and some $\alpha\in (0,1)$, then $\cR_n^{\rm pr}(0,y_n)$ becomes equidistributed on $\cM$ with respect to $\mu$ as $n\to\infty$. Moreover, under the simple symmetry relation that for $\gcd(j,n)=1$ and $y>0$
\begin{equation}\label{equ:symmetryzeor}
\G\left(\tfrac{j}{n}+iy\right)=\G \left(-\tfrac{\overline{j}}{n}+ \tfrac{i}{n^2y}\right),
\end{equation} 
%for any $\G (j/n+iy)\in \cR_n^{\rm pr}(0,y)$, 
one can extend this equidistribution result to the range $\alpha\in (1,2)${{; this improves the previous work of Demirci Akarsu \cite[Theorem 2]{Demirci2014} which confirms equdistribution of $\{\cR_n^{\rm pr}(0, c/n^{\alpha})\}_{n\in\N}$ for $\alpha\in (\frac32, 2)$.}} Here $\overline{j}\in (\Z/n\Z)^{\times}$ denotes the multiplicative inverse of $j\in (\Z/n\Z)^{\times}$. The equidistribution for the case $\alpha=1$ was later proved by Einsiedler--Luethi--Shah \cite{EinsiedlerLuethiShah2020}{{; Jana \cite[Theorem 1]{Jana2021} recently gave an alternative spectral proof to this equidistribution result. We also mention that both \cite[Theorem 2]{Demirci2014} and \cite[Theorem 1]{Jana2021} are valid in the same setting as \cite{EinsiedlerLuethiShah2020}, namely, on the product of the unit tangent bundle of the modular surface and a torus.}} When $\alpha=2$ the equidistribution fails as the aforementioned symmetry implies that $\cR^{\rm pr}_n(0,c/n^2)=\cR_n^{\rm pr}(0, 1/c)$ is always trapped in the closed horocycle $\cH_{1/c}$. For the same reason, when $\alpha>2$ (or more generally for any sequence satisfying $n^2y_n\to0$), one has with $\cR^{\rm pr}_n(0,c/n^{\alpha})=\cR_n^{\rm pr}(0, n^{\alpha-2}/c)\subset \cH_{n^{\alpha-2}/c}$ a full escape to the cusp of $\cM$ as $n\to\infty$. It is worth noting that while the symmetry \eqref{equ:symmetryzeor} still holds %(qualitative) 
for rational translates (cf. \lemref{lem:symmetryrational}), it breaks down for irrational translates.

\subsection{Statements of the results} 
We will state here the main results of this paper, and postpone the discussion of their proofs to the next subsection. Let $\mu_{\cM}:=\mu(\cM)^{-1}\mu$ be the normalized hyperbolic area on $\cM$. For any $n\in\N$, $x\in\R/\Z$ and $y>0$ let $\delta_{n,x,y}$ and $\delta_{n,x,y}^{\rm pr}$ denote the normalized probability counting measure supported on $\cR_n(x,y)$ and $\cR_n^{\rm pr}(x,y)$ respectively. That is, for any $\Psi\in C_c^{\infty}(\cM)$,
$$\delta_{n,x,y}(\Psi)=\frac{1}{n}\sum_{j=0}^{n-1}\Psi(x+\tfrac{j}{n}+iy),$$
and 
$$\delta_{n,x,y}^{\rm pr}(\Psi)=\frac{1}{\varphi(n)}\sum_{j\in (\Z/n\Z)^{\times}}\Psi(x+\tfrac{j}{n}+iy),$$
where $\varphi$ is Euler's totient function. Here and throughout, for any measure $\nu$ on $\cM$, we set $\nu(\Psi):=\int_{\cM}\Psi(z)d\nu(z)$.

Using spectral expansion and collecting estimates on the Fourier coefficients of Hecke--Maass forms and Eisenstein series, we obtain the following effective result, which yields equidistribution when the sequence is within a certain polynomial range.
\begin{Thm}\label{thm:equ}
Let $\cM$ be the modular surface. For any $\Psi\in C_c^{\infty}(\cM)$, for any $n\in\N$, $x\in \R/\Z$ and $y>0$ we have
$$\left|\delta_{n,x,y}(\Psi)-\mu_{\cM}(\Psi)\right|\ll_{\e}\cS_{2,2}(\Psi)\left(y^{1/2}+ n^{-1}y^{-(1/2+\theta+\e)}\right),$$
and
$$\left|\delta_{n,x,y}^{\rm pr}(\Psi)-\mu_{\cM}(\Psi)\right|\ll_{\e}\cS_{2,2}(\Psi)\left(y^{1/2}+n^{-1+\e}y^{-(1/2+\theta+\e)}\right),$$
where $\theta=7/64$ is the current best known bound towards the Ramanujan conjecture (which implies $\theta=0$) and $\cS_{2,2}$ is a ''$L^2$, order-$2$'' Sobolev norm on $C_c^{\infty}(\cM)$, see \secref{sec:soblev}. 
%In particular, let $\{y_n\}_{n\in\N}$ be a decreasing sequence satisfying that $\lim\limits_{n\to\infty}y_n=0$ and $y_n\gg n^{-\alpha}$ for some $\alpha\in (0,\frac{2}{1+2\theta})$. Then for any $0\leq x<1$ the sets 
%$$\cS_n(x,y_n):=\left\{\G (x+\frac{j}{n}+iy_n):1\leq j\leq n\right\}$$ 
%become equidistributed on $\G\bk \bH$ as $n\to \infty$.
\end{Thm}
If $\{y_n\}_{n\in\N}$ is a sequence of positive numbers satisfying $\lim\limits_{n\to\infty}y_n=0$ and $y_n\gg 1/n^{\alpha}$ for some fixed $\alpha\in \left(0, \frac{2}{1+2\theta}\right)=(0,\tfrac{64}{39})$, then \thmref{thm:equ} implies that for any translate $x\in \R/\Z$, both $\left\{\cR_n(x,y_n)\right\}_{n\in\N}$ and $\left\{\cR_n^{\rm pr}(x,y_n)\right\}_{n\in\N}$ become equidistributed on $\cM$ with respect to $\mu_{\cM}$ as $n\to \infty$. In particular, it gives an alternative -- spectral -- proof to the aforementioned results of Luethi \cite{Luethi2021} and Einsiedler--Luethi--Shah \cite{EinsiedlerLuethiShah2020}. The upper bound $\tfrac{2}{1+2\theta}$ is the natural barrier for our spectral methods. Nevertheless, when $x$ is a rational translate, a generalization of the symmetry \eqref{equ:symmetryzeor} allows to go beyond this barrier, and to prove unconditionally the remaining range $\alpha\in[\tfrac{2}{1+2\theta},2)$, as holds in the case of $\{\cR_n^{\rm pr}(0,y_n)\}_{n\in\N}$. %The key ingredient is the following symmetry lemma for rational translates, which generalizes the symmetry 

\begin{Thm}\label{thm:fullrange}
Let $x=p/q$ be a primitive rational number, i.e. $\gcd(p,q)=1$. Let $\{y_n\}_{n\in\N}$ be a sequence of positive numbers satisfying $y_n\asymp 1/n^{\alpha}$ for some fixed $\alpha\in [\tfrac{2}{1+2\theta},2)$. Then both $\left\{\delta_{n,x,y_n}\right\}_{n\in\N_{q}}$ and $\left\{\delta_{n,x,y_n}^{\rm pr}\right\}_{n\in\N_{q}^{\rm pr}}$ weakly converge to $\mu_{\cM}$ as $n$ goes to infinity, where 
$$\N_q:=\{n\in\N: \gcd(n^2, q)\mid n\}\quad \textrm{and}\quad \N_q^{\rm pr}:=\{n\in\N: \gcd(n,q)=1\}.$$
%$$
%\N_\alpha = \begin{cases} \N & \text{ if } 0<\alpha< \tfrac{2}{1+2\theta},\\
%\{ n\in\N: \gcd(n^2,q)\mid n\} & \text{ if } \tfrac{2}{1+2\theta}\leq \alpha<2,
%\end{cases}
%$$
%and
%$$
%\N^{\rm pr}_\alpha = \begin{cases} \N & \text{ if } 0<\alpha< \tfrac{2}{1+2\theta},\\
%\{ n\in\N: \gcd(n,q)=1\} & \text{ if } \tfrac{2}{1+2\theta}\leq \alpha<2.
%\end{cases}
%$$
%$\N_{\alpha}=\N_{\alpha}^{\rm pr}=\N$ if $0<\alpha< \tfrac{64}{39}$ and $\N_{\alpha}=\left\{n\in\N: \gcd(n^2, q)\mid n\right\}$ and $\N_{\alpha}^{\rm pr}=\{n\in\N:\gcd(n,q)=1\}$ if $\tfrac{64}{39}\leq \alpha<2$.
\end{Thm}
\begin{rmk}
If $q$ is squarefree, then the condition $\gcd(n^2, q)\mid n$ is void. Thus for such $q$, \thmref{thm:fullrange} (together with \thmref{thm:equ}) confirms the equidistribution of the sample points $\cR_n(p/q, y_n)$ (with $y_n\asymp 1/n^{\alpha}$) along the full set of positive integers for any $0<\alpha<2$. 
\end{rmk}

As a byproduct of our analysis, we also have the following non-equidistribution result for rational translates, giving infinitely many explicit limiting measures. Let us first fix some notation. For each $m\in\N$, let 
\begin{equation}\label{equ:primesub}
\bP_{m}:=\{n=m\ell\in\N : \textrm{$\ell$ is a prime number and $\ell\nmid m$}\}.
\end{equation} 
For each $Y>0$, we denote by $\mu_{Y}$ the uniform probability measure supported on the closed horocycle $\cH_Y$. For each $m\in\N$ and $Y>0$, we define the probability measure $\nu_{m,Y}$ on $\cM$ by
\begin{equation}\label{equ:limmea}
\nu_{m,Y}:=\frac{1}{m}\sum_{d\mid m}\varphi(\tfrac{m}{d})\mu_{d^2Y}.
\end{equation}

\begin{Thm}\label{thm:rationalnonequ}
Keep the notation as above. Let $x=p/q$ be a primitive rational number and let $\{y_n\}_{n\in\N}$ be a sequence of positive numbers. 
\begin{enumerate}
\item If $y_n=c/n^2$ for some constant $c>0$, then for any $m\in\N_q$ and for any $\Psi\in C_c^{\infty}(\cM)$
$$\lim\limits_{\substack{n\to\infty\\ \gcd(n,q)=1}}\delta_{n,x,y_n}^{\rm pr}(\Psi)=\mu_{\tfrac{1}{cq^2}}(\Psi)\quad \textrm{and}\quad \lim\limits_{\substack{n\to\infty\\ n\in\bP_{m}}}\delta_{n,x,y_n}(\Psi)=\nu_{m,\tfrac{\gcd(m,q)^2}{cq^2}}(\Psi).$$
\item If $\lim\limits_{n\to\infty}n^2y_n=0$, then both sequences $\{\cR_n(x,y_n)\}_{n\in\N}$ and $\{\cR^{\rm pr}_n(x,y_n)\}_{n\in\N}$ fully escape to the cusp of $\cM$.
\end{enumerate}
\end{Thm}
%Then \thmref{thm:equ} implies that for any translate $x\in \R/\Z$, both $\left\{\cR_n(x,y_n)\right\}_{n\in\N}$ and $\left\{\cR_n^{\rm pr}(x,y_n)\right\}_{n\in\N}$ become equidistributed on $\cM$ with respect to $\mu_{\cM}$ as $n\to \infty$. In particular, assuming the Ramanujan conjecture (which implies $\theta=0$) the equidistribution holds for any $\alpha\in (0,2)$. 
Our next result shows that, similar to the rational translate case, %similarly to the $\{\cR_n^{\rm pr}(0,y_n)\}_{n\in\N}$ case, 
equidistribution fails for generic translates as soon as $\{y_n\}_{n\in\N}$ decays logarithmically faster than $1/n^2$. %In fact, we will prove non-equidistribution by proving a much stronger result that the sample points fully escape to the cusp along some subsequences (with an explicit rate measured by a distance function). More precisely, we prove the following:

\begin{Thm}\label{thm:nonequ2intro}
Let $d_{\cM}(\cdot,\cdot)$ be the distance function on $\cM$ induced from the hyperbolic distance function on $\bH$. Fix $\G z_0\in\cM$. Let $\{y_n\}_{n\in\N}$ be a sequence of positive numbers satisfying $y_n\asymp 1/(n^2\log ^{\beta} n)$ for some fixed $0< \beta< 2$. Then for almost every $x\in\R/\Z$ 
\begin{equation}\label{equ:loglawdio}
\limsup_{n\to\infty}\frac{\inf_{\G z\in \cR_n(x,y_n)}d_{\cM}\left(\G z_0, \G z\right)}{\log\log n}\geq \min\{\beta, 2-\beta\}.
\end{equation}
%In particular, for such sequence $\{y_n\}_{n\in\N}$ we have for almost every $x\in [0,1)$ there exists a subsequence $\cN_x$ $($which may depend on $x$$)$ such that $\delta_{n,x,y_n}\xrightarrow{w^*} 0$ as $n\in \cN_x$ goes to infinity.
%\begin{equation}\label{equ:fullescape}
%\delta_{x,y_n}:=\frac{1}{n}\sum_{j=1}^n\delta_{\G u_{x+j/n}a_{y_n}}\ \xrightarrow[n\to\infty]{n\in \cN_x}\ 0.
%\end{equation}
%Let $y_n=\frac{1}{n^2\log^{1-\e} n}$ for some $\e\geq 0$. Then for almost every $x\in [0, 1)$ there exists an infinite subsequence $\{n_k\}\subset \N$ (which may depend on $x$) such that both $\c\cR^{\rm pr}_{n_k}(x,y_{n_k})$ and $\cS_{n_k}(x,y_{n_k})$ escape to infinity as $k\to\infty$.
\end{Thm}

This implies that for almost every $x\in\R/\Z$, there exists an unbounded subsequence of $\N$ such that along this subsequence
$$\inf_{\G z\in \cR_n(x,y_n)}d_{\cM}\left(\G z_0, \G z\right)\geq \left(\alpha-\e\right)\log\log n,$$
where $\alpha=\min\{\beta,2-\beta\}$.
%\begin{rmk}\label{rmk:escape}
%We note that since $\inf_{\G z\in \cR^{\rm pr}_n(x,y_n)}d_{\cM}\left(\G z_0, \G z\right)\geq \inf_{\G z\in S_n(x,y_n)}d_{\cM}\left(\G z_0, \G z\right)$, \thmref{thm:nonequ2intro} also holds for $\cR^{\rm pr}_n(x,y_n)$.
%\end{rmk}
That is, for almost every $x\in\R/\Z$, all the sample points $\cR_{n}(x,y_n)$ (and hence also $\cR_n^{\rm pr}(x,y_n)$) are moving towards the cusp of $\cM$ along this subsequence, and eventually escape to the cusp as $n$ in this subsequence goes to infinity.

Our proof of \thmref{thm:nonequ2intro} relies on connections to Diophantine approximation theory. This viewpoint comes with inherent limitations; in the specific setting $y_n\asymp 1/(n^2\log^\beta n)$, Khintchine's approximation theorem guarantees full escape to the cusp almost surely, but this argument does not extend to any sequence $\{y_n\}_{n\in\N}$ that decays polynomially faster than $1/n^2$, see \secref{sec:dis} for a more detailed discussion.
%As a result, the counterexamples we have in \thmref{thm:nonequ2intro} are very specific. In particular, the argument falls short of implying the almost sure escape to the cusp whenever the sequence $\{y_n\}_{n\in\N}$ decays polynomially faster than $1/n^{2}$. 
It is thus interesting to study the cases when $\{y_n\}_{n\in\N}$ is beyond the ranges in \thmref{thm:equ} and \thmref{thm:nonequ2intro}. %for example when $\{y_n\}_{n\in\N}$ decays super-polynomially.

Indeed, the rest of our results deal with sequences $\{y_n\}_{n\in\N}$ that can decay \textit{arbitrarily fast}, and give both positive and negative results. This is the main novelty of this paper; the handling of cases in which the sample points can be \textit{arbitrarily sparse} on the closed horocycles they lie on. 
{{We now state the main novel aspect of this paper:
\begin{Thm}\label{thm:mainthm}
For any sequence of positive numbers $\{c_n\}_{n\in\N}$, there exists a sequence $\{y_n\}_{n\in\N}$ satisfying $0<y_n<c_n$ for each $n\in\N$ and such that for almost every $x\in\R/\Z$ the set of limiting measures of $\{\delta_{n,x,y_n}\}_{n\in\N}$ and $\{\delta_{n,x,y_n}^{\rm pr}\}_{n\in\N}$ both contain the uniform measure $\mu_{\cM}$, the zero measure, and singular probability measures.
\end{Thm}

}}
{{\thmref{thm:mainthm} is a sum of three more precise theorems, which each handles a specific limiting measure, {{and which we discuss in the next subsection.}} %see also \rmkref{rmk:disthimma} for a discussion how \thmref{thm:mainthm} follows from these three theorems.
}}

\subsection{Discussion of the results}\label{sec:dis}
%Before ending the introduction, let us make some comments on our results and the proofs.\\
%
Our proofs of \thmref{thm:equ} and \thmref{thm:fullrange} rely on spectral estimates collected in the recent paper of Kelmer and Kontorovich \cite{KelmerKontorovich2020}, with a necessary refinement of \cite[(3.6)]{KelmerKontorovich2020} in the form of \propref{prop:fe2}, which comes at the cost of a higher degree Sobolev norm. 
{{This strategy is standard and is also found in \cite{ClozelUllmo2004,MarklofStrombergsson2003,SarnakUbis2015,Jana2021}, to name just a few recent papers on related problems.}}
%{{The recent article by Jana \cite{Jana2021} proves the equidistribution for $\{\cR_{n}^{\rm pr}(0, 1/n)\}_{n\N}$ using similar spectral arguments. In fact he proved this equidistribution result in the setting of \cite{EinsiedlerLuethiShah2020}, namely on the product of the unit tangent bundle of the modular surface and a torus. His argument also works for any translate $x\in \R/\Z$.}} 
The analysis in \cite{KelmerKontorovich2020} was carried out in a more general setting, namely for the congruence covers $\G_0(p)\bk\bH$ with $p$ a prime number. \thmref{thm:equ} can be extended to that more general setting, see \rmkref{rmk:extension}. With these spectral estimates in hand, we further prove an effective non-equidistribution result for rational translates from which part (1) of \thmref{thm:rationalnonequ} follows, %when $\{y_n\}_{n\in\N}$ is beyond the range in \thmref{thm:fullrange}, 
see \thmref{thm:effnonequ}. Part (2) of \thmref{thm:rationalnonequ} is an easy application of the symmetry \eqref{equ:symmetryzeor}.
{{
\begin{remark}%\label{rmk:sar-ubi}
As was pointed out to us by Asaf Katz, we could also have used the  estimates from \cite[Proposition 3.1]{SarnakUbis2015} in place of \cite[Proposition 3.4]{KelmerKontorovich2020}, which in our specific setting, give the same equidistribution range (with a higher degree Sobolev norm). We also mention that the estimates in \cite[Proposition 3.1]{SarnakUbis2015} are valid in the setting of $\G_0(q)\bk \SL_2(\R)$ with $q\in\N$, and thus imply an effective equidistribution result analogous to \thmref{thm:equ} in this generality.
\end{remark}%
}}
%\\

As mentioned earlier, a generalization of the symmetry \eqref{equ:symmetryzeor} is available for rational translates but breaks down for irrational translates. To handle irrational translates, we approximate them by rational ones to apply the symmetry relation, see \lemref{lem:counter1}. This is where Diophantine approximation kicks in. Similar ideas were also used in \cite[Section 7]{MarklofStrombergsson2003} to construct counterexamples in their setting. In fact, we prove \thmref{thm:nonequ2intro} by proving a more general result that captures the cusp excursion rates of the sample points $\cR_n(x,y_n)$ in terms of the Diophantine properties of the translate $x$, see \thmref{thm:nonequ2}. \thmref{thm:nonequ2intro} will then follow from \thmref{thm:nonequ2} by imposing a Diophantine condition which ensures cusp excursion, while also holds for almost every translate thanks to Khintchine's approximation theorem. This Diophantine condition accounts for the tight restrictions on $\{y_n\}_{n\in\N}$ in \thmref{thm:nonequ2intro}. On the other hand, assuming an even stronger Diophantine condition (which holds for a null set of translates), we can handle sequences decaying polynomially faster than $1/n^2$ with a much faster excursion rate towards the cusp, see \thmref{thm:negfas}. We also prove a non-equidistribution result (which, this time, holds for \emph{every} $x$) when $y_n=c/n^2$ and the constant $c$ is restricted to some range, see \thmref{thm:nonequ1}. The trade-off of this upgrade from \thmref{thm:nonequ2intro} to the everywhere non-equidistribution result is that we can no longer prove the full escape to the cusp along subsequences as in \thmref{thm:nonequ2intro}.\\

{{As mentioned before, \thmref{thm:mainthm} follows from three more precise theorems which each handles a specific limiting measure.}}
Our first result confirms equidistribution almost surely along a fixed subsequence of $\N$ for \textit{any} sequence $\{y_n\}_{n\in\N}$ decaying {{at least polynomially}}.
%faster than a certain explicit rate. %then the sample points $\left\{S_n(x,y_n)\right\}_{n\in\N}$ and $\left\{S_n^{\rm pr}(x,y_n)\right\}_{n\in\N}$ become equidistributed on $\cM$ along this subsequence almost surely. 
\begin{Thm}\label{thm:equipar}
{{Fix $\alpha>0$. Then there exists a fixed unbounded subsequence $\cN\subset \N$ such that for any sequence of positive numbers $\{y_n\}_{n\in\N}$ satisfying $y_n\ll n^{-\alpha}$ and for almost every $x\in \R/\Z$, both $\delta_{n,x,y_n}$ and $\delta^{\rm pr}_{n,x,y_n}$ weakly converge to $\mu_{\cM}$ as $n\in \cN$ goes to infinity. }}
%There exists a fixed unbounded subsequence $\cN\subset \N$ such that for any sequence of positive numbers $\{y_n\}_{n\in\N}$ satisfying $y_n\ll n^{-2+4\theta}$, and for almost every $x\in \R/\Z$, both $\delta_{n,x,y_n}$ and $\delta^{\rm pr}_{n,x,y_n}$ weakly converge to $\mu_{\cM}$ as $n\in \cN$ goes to infinity.
%Let $\cN\subset \N$ be an infinite subsequence such that $\sum_{n\in \cN}n^{-\alpha}$ converges for some constant $0<\alpha<1/2$. Then for any sequence  with $y_n\ll 1/n$, we have for almost every $x\in [0,1)$
%$$\delta_{x,y_n}:=\frac{1}{n}\sum_{j=1}^n\delta_{\G u_{x+j/n}a_{y_n}}\ \xrightarrow[n\to\infty]{n\in \cN}\ \mu_{\G}.$$
\end{Thm}

\begin{rmk}\label{rmk:subse}
It will be clear from our proof that one can take $\cN\subset \N$ to be any subsequence satisfying $\sum_{n\in\cN}n^{-c}<\infty$ for some positive {{$c<\min\{\frac{\alpha}{2}, 1-2\theta\}$}}, e.g. {{we may take}} $\cN=\left\{\left \lfloor{n^{\kappa}}\right \rfloor\right\}_{n\in\N}$ for {{any $\kappa>1/\min\{\frac{\alpha}{2}, 1-2\theta\}$}}. %$\beta:=\min\{\frac{\alpha}{2}, 1-2\theta\}$%\frac{32}{25}$. 
%In contrast, to include the case of $\{\delta^{\rm pr}_{n,x,y_n}\}_{n\in\N}$, we need the extra assumption that $\cN$ being a subsequence of the set of prime numbers.
%If we are only concerned about the limiting distribution of $\{\cR_n(x, y_n)\}_{n\in\N}$, then in \thmref{thm:equipar} we can take any subsequence $\cN\subset \N$ satisfying $\sum_{n\in\cN}n^{-\alpha}<\infty$ for some $\alpha\in (0,1-2\theta)$. For instance we can take $\cN=\left\{\left \lfloor{n^{\beta}}\right \rfloor\right\}_{n\in\N}$ for any $\beta>1/(1-2\theta)$. In contrast, to include the case of $\{\cR^{\rm pr}_n(x, y_n)\}_{n\in\N}$, the subsequence we can take is then much more sparse with the extra assumption that $\cN$ being a subsequence of the set of prime numbers.
\end{rmk}

\thmref{thm:equipar} follows from a second moment estimate for the discrepancies $|\delta_{n,x,y}-\mu_{\cM}|$ and $|\delta_{n,x,y}^{\rm pr}-\mu_{\cM}|$ along the closed horocycle $\cH_y$ (\thmref{thm:secmom}) together with a standard Borel-Cantelli type argument. This was also the strategy used in \cite{MarklofStrombergsson2003} when studying the Kronecker sequences in \eqref{equ:kroneseq}.
%$$\{\G(j\beta +iy_n)\in \cM:1\leq j\leq n\}.$$
Along these lines, they deduce from spectral estimates the equidistribution for almost every $\beta\in\R$ along a fixed subsequence $\{n^k\}_{n\in\N}$ when $y_n\asymp n^{-\alpha}$ with $k\in\N$ depending on $\alpha>0$. Then, using a continuity argument, this result is upgraded to the equidistribution along the full sequence of positive integers, see \cite[Section 4]{MarklofStrombergsson2003}. This continuity argument fails in our situation. Instead of applying directly spectral estimates to the second moment formulas, we express the latter in terms of certain Hecke operators (\propref{prop:secondmoment}), and rely on available (spectral) bounds for their operator norm, see \cite{GoldsteinMayer2003}. Contrarily to spectral estimates, the recourse to Hecke operators allows us to %handle sequences decaying arbitrarily fast and for
have a uniform subsequence $\cN$ which is valid for all $\{y_n\}_{n\in\N}$ decaying {{at least polynomially.}}
%faster than $n^{-2+4\theta}$. %See also a recent work of Bersudsky \cite[Theorem 1.5]{Bersudsky2020} on translates of rational points on dilations of analytic curves projected to the torus using similar moment arguments. 
\\

{{Next, we show that there exists a sequence $\{y_n\}_{n\in\N}$ decaying arbitrarily rapidly such that for almost every $x$, $\cR_n(x,y_n)$ (and thus also $\cR_n^{\rm pr}(x,y)$) escapes to the cusp with a certain rate along subsequences.}}

%In view of \thmref{thm:equipar} one may ask whether, for almost every $x\in\R/\Z$ and all sequences $\{y_n\}_{n\in\N}$ decaying sufficiently fast, the set of limiting measures of $\{\delta_{n,x,y_n}\}_{n\in\N}$ (and respectively $
%\{\delta_{n,x,y_n}^{\rm pr}\}_{n\in\N}$) consists only of $\mu_{\cM}$.
%sample points $S_n(x,y_n)$ and $S_n^{\rm pr}(x,y_n)$ always equidistribute on $\cM$ almost surely, or in other words, whether for almost every $x\in\R/\Z$, the set of limiting measures of $\delta_{n,x,y_n}$ and $
%\delta_{n,x,y_n}^{\rm pr}$ only contain $\mu_{\cM}$. 
%We answer this question negatively by showing that there will always exist sequences decaying faster than any prescribed sequence such that, almost surely, the set of limiting measures also contains the trivial measure. %Similar as in \thmref{thm:nonequ2intro} we prove this result by proving the following logarithm law for the simultaneous cusp excursions of our sample points. %for sequences that can decay arbitrarily fast.%$S_n(x,y_n)$ which can hold for sequences of arbitrarily fast decaying rates.
\begin{Thm}\label{thm:negative2}
Fix $\G z_0\in\cM$. For any sequence of positive numbers $\left\{c_n\right\}_{n\in\N}$, there exists a sequence $\{y_n\}_{n\in \N}$ satisfying $0<y_n< c_n$ for each $n\in\N$ and such that for almost every $x\in \R/\Z$ %there exists an infinite subsequence $\mathcal{N}_x\subset \mathbb{P}$ $($which may depend on $x$$)$ such that 
%$$\delta_{x,y_p}:=\frac{1}{p}\sum_{j=1}^p\delta_{\G u_{x+j/p}a_{y_p}} \xrightarrow[p\to \infty]{p\in \mathcal{N}_x}\ 0.$$
%as $p\in \mathcal{N}_x$ goes to infinity.
\begin{equation}\label{equ:loglaw}
\limsup_{n\to\infty}\frac{\inf_{\G z\in \cR_n(x,y_n)}d_{\cM}\left(\G z_0, \G z\right)}{\log\log n}\geq 1.
\end{equation}
%In particular, for such sequence $\{y_n\}_{n\in\N}$ we have for almost every $x\in [0,1)$ there exists a subsequence $\cN_x$ $($which may depend on $x$$)$ such that $\delta_{n,x,y_n}\xrightarrow{w^*} 0$ as $n\in \cN_x$ goes to infinity.
\end{Thm}
Finally, we show that escape to the cusp is not the only obstacle to equidistribution. %We now state our final result:
\begin{Thm}\label{thm:compactsupport}
Let $m\in\N$ and $Y>0$ satisfy $m^2Y>1$. Let $\bP_m\subset \N$ and $\nu_{m,Y}$ be as defined in \eqref{equ:primesub} and \eqref{equ:limmea} respectively. For any sequence of positive numbers $\{c_n\}_{n\in\bP_{m}}$, there exists a sequence $\{y_n\}_{n\in\bP_m}$ satisfying $0<y_n< c_n$ for all $n\in\bP_m$ such that for almost every $x\in\R/\Z$, the set of limiting measures of $\{\delta_{n,x,y_n}\}_{n\in\bP_m}$ contains $\nu_{m,Y}$.
\end{Thm} 
\begin{rmk}\label{rmk:disthimma}
We note that $\bP_1$ is the set of prime numbers and $\nu_{1,Y}=\mu_Y$. Since 
$$\delta^{\rm pr}_{p,x,y}(\Psi)=\tfrac{p}{p-1}\delta_{p,x,y}(\Psi)+O(p^{-1}\|\Psi\|_{\infty})$$ 
whenever $p$ is a prime number, when $m=1$ the conclusion of \thmref{thm:compactsupport} also holds for the sequence $\{\delta_{n,x,y_n}^{\rm pr}\}_{n\in\bP_1}$.
%We will prove \thmref{thm:compactsupport} by showing that for almost every $x\in\R/\Z$, all (but one) points in $\cR_p(x,y_p)$ evenly distributed (up to the one missing point) on some closed horocycles which converge to $\cH_Y$ along some subsequences of $\bP$. In particular, since at each prime step $\cR_p^{\rm pr}(x,y_p)$ only differs from $\cR_p(x,y_p)$ by one point, \thmref{thm:compactsupport} also holds for $\{\delta_{p,x,y_p}^{\rm pr}\}_{p\in\mathbb{P}}$. 
We also note that it will be clear from our proof that \thmref{thm:negative2} and \thmref{thm:compactsupport} can be combined. In fact, our argument shows that there always exists a sequence $\{y_n\}_{n\in\N}$ decaying faster than any prescribed sequence such that for almost every $x\in\R/\Z$ the set of limiting measures of $\left\{\delta_{n,x,y_n}\right\}_{n\in\N}$ contains the trivial measure and $\nu_{m,Y}$ for any finitely many pairs $(m,Y)\in\N\times \R_{>0}$ with $m^2Y>1$, see \rmkref{rmk:combine}. Moreover, in view of \thmref{thm:equipar} if {{$y_n\ll n^{-\alpha}$ for some $\alpha>0$}}, then it also contains the hyperbolic area $\mu_{\cM}$ almost surely.
\end{rmk}

For the rest of this introduction we describe the strategy of our proof to \thmref{thm:negative2} (\thmref{thm:compactsupport} follows from similar ideas). %The starting point of our analysis is the simple observation that for any $n\in\N$, for any $0\leq j\leq n-1$ and any $z\in\bH$, $z+j/n=u_{j/n}z$. %that is we can view our sample points $S_n(x,y)$ as the orbit $\langle u_{1/n}\rangle \cdot (x+iy)$ under the left regular action. 
To detect cusp excursions, we study for each $n\in\N$ the occurrence of the events 
\begin{equation}\label{equ:cuspex}
\G \left(x+\tfrac{j}{n}+iy_n\right)\in \cC\quad \textrm{for all $0\leq j\leq n-1$},
\end{equation}
where $\cC\subset \cM$ is some fixed cusp neighborhood of $\cM$. More precisely, we determine when the limsup set $I_{\infty}=\limsup_{n\to\infty}I_n$ is of full measure, where for each $n\in\N$,
$$I_n:=\{x\in \R/\Z:\cR_n(x,y_n)\subset \cC\}$$
consists of translates $x\in \R/\Z$ for which the events in \eqref{equ:cuspex} occur. This requires to study the left regular $u_{1/n}$-action on $\cC\subset \cM$ 
%(noting that $\G(x+j/n+iy)=\G u_{j/n} (x+iy)$) %However, this reformulation requires 
and thus calls for the underlying lattice to be normalized by $u_{1/n}$. Therefore, we construct an explicit tower of coverings $\{\Gamma_n\backslash\bH\}_{n\in\N}$ in which each $\Gamma_n$ is a congruence subgroup normalized by $u_{1/n}$. 
{{We note that the existence of such $\G_n<\G$ is the starting point of our proof and it relies on the assumption that $\G=\SL_2(\Z)$; this construction would fail for $\G$ replaced by a non-arithmetic lattice.}}
%However, this left regular action only makes sense when the underlying lattice is normalized by $u_{1/n}$. To overcome this difficulty, we lift this problem to certain congruence covers. More precisely, we will define a sequence of congruence subgroups $\{\G_n\}_{n\in\N}$ (see \eqref{equ:latdes}) such that $\G_n$ is normalized by $u_{1/n}$ (hence also by $u_{j/n}$ for any $0\leq j\leq n-1$), see \lemref{lem:explicitdes}. We will then study the left regular $u_{1/n}$-action on certain cusp neighborhoods on $\G_n\bk\bH$ and then descend these cusp excursions back to the modular surface. 

The key ingredient of the proof will be a sufficient condition 
%\footnote{It was communicated to us by Str\"ombergsson that using a number theoretic interpretation of this sufficient condition and some elementary estimates one can prove \thmref{thm:negative2} without going into these congruence covers, see \rmkref{rmk:andreas}.} 
which states that if a point $\G_n (x+iy_n)\in \G_n\bk \bH$ visits a certain cusp neighborhood $\cC_n$ on $\G_n\bk \bH$, then the events in \eqref{equ:cuspex} will be realized for $x\in \R/\Z$, that is, $x\in I_n$, see \lemref{lem:cuspanalysis2}. 
%Next, note that for our purpose we need to study, at the $n$-th step, the events in \eqref{equ:cuspex} for a generic point on the closed horocycle (of height $y_n$), rather than a generic point on the modular surface. That is, we need to know, for fixed $0<y_n<1$, the measure of the subset
%$$\cI_n:=\{x\in \R/\Z:\textrm{events in \eqref{equ:cuspex} happen for $z=x+iy_n$}\}\subset \R/\Z.$$
%More precisely, we want to understand when the limsup set $\cI_{\infty}:=\limsup_{n\to\infty}\cI_n$ is of full measure. 
Using this sufficient condition, we can then relate the measure of $I_n$ to the proportion of certain closed horocycles on $\G_n\bk\bH$ visiting the cusp neighborhood $\cC_{n}\subset \G_n\bk\bH$, which in turn, using the equidistribution of expanding closed horocycles on $\G_n\bk\bH$, can be estimated for $y_n$ sufficiently small. 
Since the sets $I_n$ also need to satisfy certain quasi-independence conditions for $I_\infty$ to have full measure (\lemref{lem:bcd}), we need to apply the equidistribution of certain subsegments of the expanding closed horocycles on $\G_n\bk\bH$.
%However, knowing the measure of these sets $\cI_n$ is still not enough since, to prove $\cI_{\infty}$ has full measure, one also needs to make sure that these sets satisfy certain quasi-independence conditions. To ensure such conditions we need to apply the equidistribution of certain subsegments of the closed horocycles on $\G_n\bk\bH$. 
More precisely, at the $n$-th step these subsegments will be taken to be the sets $I_m$ for all $m<n$. These subsegment are finite disjoint unions of subintervals whose number and size  depend sensitively on the height parameters $\{y_m\}_{m<n}$, see \rmkref{rmk:disopen}. If there would exist an effective equidistribution result which would be insensitive to the geometry of these subsegments, that is, for which the error term depends only on the measure of these subsegments, then we would have an effective control on the sequence $\{y_n\}_{n\in\N}$ in \thmref{thm:negative2} (and similarly also in \thmref{thm:compactsupport}). However, it is not clear to us whether one should expect such an effective equidistribution result.

{{Finally, we note that it was communicated to us by Str\"ombergsson that using a number theoretic interpretation of the aforementioned sufficient condition and some elementary estimates, one can alternatively prove \thmref{thm:negative2} without going into these congruence covers, see \rmkref{rmk:andreas}.
}}
%\begin{rmk}
%It was communicated to us by Str\"ombergsson that there is an elementary proof to \thmref{} without 
%\end{rmk}
%the lack of good effective equidistribution results for these subsegments is the main reason for which we do not have an effective control on the sequence $\{y_n\}_{n\in\N}$ in \thmref{thm:negative2}.
%is a disjoint union of finitely many open intervals with the number and size of these open intervals depending very senstively on the sequence $\{y_m\}_{m<n}$, see \rmkref{rmk:disopen}.We note that this is the main reason for which we do not have an effective control on the sequence $\{y_n\}_{n\in\N}$ in \thmref{thm:negative2}. %We note that  %and as far as we know, there is currently no very good effective equidistribution result handling these subsegments.
%We note that our proof is of algebraic nature relying on the simple fact that for each $n\in\N$ there is a large subgroup of $\SL_2(\Z)$ normalized by $u_{1/n}$. This fact fails when $\G$ is replaced by a non-arithmetic lattice.
%To deduce such a sufficient condition, we need to study the left regular $u_{1/n}$-action on the cusp neighborhoods of $\G_n\bk\bH$. For this we will study, in more generality, the left regular $u_{1/n}$-action on incomplete Eisenstein series on $\G_n\bk \bH$, namely, the generating functions of the continuous and residual spectrum of $L^2(\G_n\bk \bH)$, see \secref{sec:leftregular}.
\subsection*{Structure of the paper}
In \secref{sec:pre} we collect some preliminary results that will be needed in the rest of the paper. In \secref{sec:equran}, we prove a key spectral estimate (\propref{prop:fe2}) and proceed 
%after collecting necessary estimates on Fourier coefficients we proceed 
to prove \thmref{thm:equ} and \thmref{thm:fullrange}. In \secref{sec:cou1} we prove \thmref{thm:nonequ2} and \thmref{thm:nonequ1} by examining the connections between Diophantine approximations and cusp excursions on the modular surface. In \secref{sec:secmom} we prove \thmref{thm:equipar} by proving a second moment bound using Hecke operators. In \secref{sec:leftregular} we study the left regular action of a normalizing element on the set of cusp neighborhoods of a congruence cover of the modular surface. Building on the results, we prove  \thmref{thm:negative2} and \thmref{thm:compactsupport} in \secref{sec:neg2}.

\subsection*{Notation}
For two positive quantities $A$ and $B$, we will use the notation $A\ll B$  {or $A=O(B)$} to mean that there is a constant $c>0$ such that $A\leq cB$, and we will use subscripts to indicate the dependence of the constant on parameters. We will write $A\asymp B$ for $A\ll B\ll A$. For any $z\in \bH$ we denote by $e(z):=e^{2\pi iz}$. For any $n\in \N$, we denote by $\prod_{d\mid n}$  the product over all positive divisors of $n$, and by $\prod_{\substack{p\mid n\\ prime}}$ the product over all prime divisors of $n$. %For the remaining of the paper we will denote by $\G_1=\SL_2(\Z)$ and we will see that this notation is consistent with our definition of the congruence subgroups $\G_n$, see \eqref{}. 
For any $x\geq 0$ and $n\in \N$, $\sigma_x(n):=\sum_{d\mid n}d^x$ is the power-$x$ divisor function which satisfies the estimate $\sigma_x(n)\ll_{\e}n^{x+\e}$ for any small $\e>0$.

\subsection*{Acknowledgements}
The first named author would like to thank Alex Kontorovich for explanations and references on Sobolev norms. The second and third named authors would like to thank Michael Bersudsky and Rene R\"{u}hr for various discussions on this problem. The third named author would also like to thank Dubi Kelmer for answering some questions and pointing out a reference to him regarding the residual spectrum of congruence subgroups. We would also like to thank Str\"ombergsson for his comments on an earlier version of this paper, especially for suggesting an alternative elementary proof to \thmref{thm:negative2}. {{We would also like to thank Asaf Katz for pointing out to us some references on spectral estimates.}}

\section{Preliminaries}\label{sec:pre}
Let $G=\SL_2(\R)$. We consider the Iwasawa decomposition $G=NAK$ with
$$N=\left\{u_x:x\in\R\right\},\quad A=\left\{a_y:y>0\right\},\quad K=\left\{k_{\theta}:0\leq \theta<2\pi\right\},$$ 
where $u_x=\left(\begin{smallmatrix}
1 & x\\
0 & 1\end{smallmatrix}\right)$, $a_y=\left(\begin{smallmatrix}
y^{1/2} & 0\\
0 & y^{-1/2}\end{smallmatrix}\right)$ and $k_{\theta}=\left(\begin{smallmatrix}
\cos\theta & \sin\theta\\
-\sin\theta & \cos\theta\end{smallmatrix}\right)$ respectively. Under the coordinates $g=u_xa_yk_{\theta}$ on $G$, the Haar measure is given (up to scalars) by
$$d g=y^{-2}dxdyd\theta.$$
%The group $G$ acts on the upper half plane $\bH$ as isometries via the M\"obius transformation: $g z=\frac{az+b}{cz+d}$ for $g=\left(\begin{smallmatrix}
%a & b\\
%c & d\end{smallmatrix}\right)\in G$ and $z\in \bH$.
The group $G$ acts on the upper half plane $\bH=\{z=x+iy\in\C:y>0\}$ via the M\"{o}bius transformation: $g z=\frac{az+b}{cz+d}$ for any $g=\left(\begin{smallmatrix}
a & b\\
c & d\end{smallmatrix}\right)\in G$ and $z\in \bH$. This action preserves the hyperbolic area $d\mu(z)=y^{-2}dxdy$ and induces an identification between $G/K$ and $\bH$.% identifying $gK\in G/K$ with $gi\in\bH$. 

Let $\G< G$ be a \textit{lattice}, that is, $\G$ is a discrete subgroup of $G$ such that the corresponding hyperbolic surface $\G\bk \bH$ has finite area (with respect to $\mu$). We denote by $\mu_{\G}:=\mu(\G\bk \bH)^{-1}\mu$ the normalized hyperbolic area on $\G\bk \bH$ such that $\mu_{\G}(\G\bk \bH)=1$. We note that when $\G=\SL_2(\Z)$ then $\mu_{\G}=\mu_{\cM}$ with $\mu_{\cM}$ the normalized hyperbolic area on the modular surface $\cM$ given as in the introduction. We note that in this case it is well known $\mu(\cM)=\pi/3$, and hence
\begin{equation}\label{equ:norhyarea}
d\mu_{\cM}(z)=\frac{3}{\pi}\frac{dxdy}{y^2}.
\end{equation}

Using the above identification between $\bH$ and $G/K$ we can identify the hyperbolic surface $\G\bk\bH$ with the locally symmetric space $\G\bk G/K$. 
%We can identify the locally symmetric space $\G\bk G/K$ with the hyperbolic surface $\G\bk \bH$ by identifying $\G gK$ with $\G g i$, where for $g=\left(\begin{smallmatrix}
%a & b\\
%c & d\end{smallmatrix}\right)\in G$ and $z\in \bH$, $g z=\frac{az+b}{cz+d}$ is the M\"obius transformation. Hence we have a natural projection $\G\bk G\to \G\bk \bH$ ($\cong \G\bk G/K$). 
We can thus view subsets of $\G\bk \bH$ as right $K$-invariant subsets of $\G\bk G$. Similarly, we can view functions on $\G\bk \bH$ as right $K$-invariant functions on $\G\bk G$. We note that using the above description of the Haar measure, the probability Haar measure on $\G\bk G$ (when restricted to the sub-family of right $K$-invariant subsets) coincides with the normalized hyperbolic area $\mu_{\G}$ on $\G\bk \bH$.
%This way then the normalized hyberbolic area $d\mu_{\cM}(z)$ on the modular surface $\cM$ coincides with the normalized Haar measure here (restricted to the family of right $K$-invariant subsets of $\SL_2(\Z)\bk G$).

\subsection{Sobolev norms}\label{sec:soblev}
In this subsection we record some useful properties about Sobolev norms. Let $\fg=\mathfrak{sl}_2(\R)$ be the Lie algebra of $G$.  Fix a basis $\mathscr{B}=\{X_1, X_2, X_3\}$ for $\fg$, and given a smooth test function $\Psi\in C^{\infty}(\G\bk G)$ we define the ``$L^p$, order-$d$" Sobolev norm $\cS^{\G}_{p,d}(\Psi)$ as
$$\cS^{\G}_{p,d}(\Psi):=\sum_{\textrm{ord}(\mathscr{D})\leq d}\|\mathscr{D}\Psi\|_{L^p(\G\bk G)},$$
where $\mathscr{D}$ runs over all monomials in $\mathscr{B}$ of order at most $d$, and the $L^p$-norm is with respect to the normalized Haar measure on $\G\bk G$. 

For any $\Psi\in C^{\infty}(\G\bk G)$ (which we think of a smooth left $\G$-invariant function on $G$) and for any $h\in G$ we denote by $L_h\Psi(g):=\Psi(h^{-1}g)$ the left regular $h$-action on $\Psi$. It is easy to check that $L_h\Psi\in C^{\infty}(h\G h^{-1}\bk G)$, and since taking Lie derivatives commutes with the left regular action, we have
\begin{equation}\label{equ:sobolevconj}
\cS_{p,d}^{\G}(\Psi)=\cS_{p,d}^{h\G h^{-1}}(L_h\Psi).
\end{equation}
Next we note that using the product rule for Lie derivatives (see e.g. \cite[p. 90]{LangSL2R}), the triangle inequality and the Cauchy-Schwarz inequality, for any monomial $\mathscr{D}$ of order $k{{\leq d}}$ we have for any smooth functions $\Psi_1,\Psi_2\in C^{\infty}(\G\bk G)$
$$\|\mathscr{D}\Psi_1\Psi_2\|_{L^p(\G\bk G)}\ll_{k}\cS_{2p,{{k}}}^{\G}(\Psi_1)\cS_{2p,{{k}}}^{\G}(\Psi_2){{\leq \cS_{2p,d}^{\G}(\Psi_1)\cS_{2p,d}^{\G}(\Psi_2).}}$$
%where the bounding {{constant}} only depends on the order of $\mathscr{D}$. 
In particular this implies that
\begin{equation}\label{equ:sobolev}
\cS^{\G}_{p,d}(\Psi_1\Psi_2)\ll_{d}\cS_{2p,d}^{\G}(\Psi_1)\cS_{2p,d}^{\G}(\Psi_2).
\end{equation}
Finally, we note that if $\G'<\G$ is a finite-index subgroup of $\G$, then there is a natural embedding $C^{\infty}(\G\bk G)\hookrightarrow C^{\infty}(\G'\bk G)$ since each $\Psi\in C^{\infty}(\G\bk G)$ can be viewed as a smooth left $\G'$-invariant function on $G$. Since the Sobolev norms are defined with respect to the normalized Haar measure on the corresponding homogeneous space, we have for $\G'< \G$ of finite index and $\Psi\in C^{\infty}(\G\bk G)$
\begin{equation}\label{equ:soblev3}
\cS_{p,d}^{\G'}(\Psi)=\cS_{p,d}^{\G}(\Psi).
\end{equation}
%Finally, to ease the notation, when $\G=\SL_2(\Z)$ we abbreviate $\cS_{p,d}^{\G}$ by $\cS_{p,d}$.
\subsection{Spectral decomposition}
Let $\G< G$ be a \textit{non-uniform} lattice, that is, $\G$ is a lattice and $\G\bk \bH$ is not compact.
%$\G$ contains a \textit{principal congruence subgroup} 
%$$\G(n):=\left\{\gamma\in \SL_2(\Z):\gamma\equiv I_2 \Mod{n}\right\}$$
%of some level $n\in \N$. 
%Let $\bH:=\{z=x+iy\in \C :y>0\}$ denote the hyperbolic upper half plane and the group $G$ acts on $\bH$ as isometries via the M\"obius transformation. We can identify $\bH$ with the symmetric space $G/K$ via the map sending $gK\in G/K$ to $g\cdot i$. Similarly, the modular surface $\G\bk \bH$ can be identified with the locally symmetric space $\G\bk G/K$, and we identify $L^2(\G\bk \bH)$ as the subspace of right $K$-invariant functions in $L^2(\G\bk G)$.
Let $\Delta=-y^2(\frac{\partial}{\partial x^2}+\frac{\partial}{\partial y^2})$ be the hyperbolic Laplace operator. It is a second order differential operator acting on $C^{\infty}(\G\bk \bH)$ and extends uniquely to a self-adjoint and positive semi-definite operator on $L^2(\G\bk \bH)$. Since $\G$ is non-uniform, the spectrum of $\Delta$ is composed of a continuous part (spanned by Eisenstein series) and a discrete part (spanned by Maass forms) which further decomposes as the cuspidal spectrum and the residual spectrum. The residual spectrum always contains the constant functions (coming from the trivial pole of the Eisenstein series). If $\G$ is a \textit{congruence subgroup}, that is, $\G$ contains a \textit{principal congruence subgroup} 
$$\G(n):=\left\{\gamma\in \SL_2(\Z):\gamma\equiv I_2 \Mod{n}\right\}$$
for some $n\in \N$, then the residual spectrum consists only of the constant functions, see e.g. \cite[Theorem 11.3]{Iwaniec2002}.

Let $\{\phi_k\}$ be an orthonormal basis of the space of cusp forms that are eigenfunctions of the Laplace operator $\Delta$.  Explicitly, for each $\phi_k$ there exists $\lambda_k\geq0$ such that
%some $s_k=\frac12+ir_k$ with $r_k\in i(0, 1/2)\cup [0,\infty)$ such that
$$\Delta\phi_k=\lambda_k \phi_k = s_k(1-s_k)\phi_k=\left(\tfrac14+r_k^2\right)\phi_k.
$$  
%Selberg's eigenvalue conjecture states that for congruence subgroups there are no cusp forms with eigenvalue $s_k(1-s_k)<1/4$, or equivalently, there is no $r_k\in i(0,1/2)$. We note that the best currently known bound towards this conjecture is due to Sarnak and Kim \cite{KimSarnak03} which implies that for every cusp form $\phi_k$, its eigenvalue $s_k(1-s_k)\geq 1/4-\theta^2$ with $\theta=7/64$. Here we record the spectral theorem on the modular surface $\cM=\SL_2(\Z)\bk\bH$ on which the Selberg's eigenvalue conjecture is known to be true. 
Selberg's eigenvalue conjecture states that for congruence subgroups, $\lambda_k\geq 1/4$, or equivalently, there is no $r_k\in i(0,1/2)$. Selberg's conjecture is known to be true for the modular surface $\mathcal{M}$, and more generally, the best known bound towards this conjecture is currently $\lambda_k\geq \tfrac14-\theta^2$, with $\theta=7/64$, which follows from the bound of Kim and Sarnak towards the Ramanujan conjecture, see \cite[p. 176]{KimSarnak03}.

%\subsubsection{Eisenstein series}
%Let us review some backgronds on the continuous spectrum, namely the Eisenstein series. Since we assume $\G$ is a congruence subgroup, the set of cusps of $\G$ can be parameterized by the set of $\G$-equivalent classes of $\Q\cup \{\infty\}$, where $s_1, s_2\in \Q\cup\{\infty\}$ is called \textit{$\G$-equivanlent} if and only if there exists some $\gamma\in \G$ such that $\gamma s_1=s_2$ with the action of $\G$ on $\Q\cup \{\infty\}$ given by the M\"obius transformation. %We denote by $\cC_{\G}$ the set of cusps of $\G$ and denote by $\varepsilon_{\infty}(\G)$ its the number of cusps. 

%For every cusp $\fa$ of $\G$ let $\G_{\fa}<\G$ be the stabilizer of $\fa$ in $\G$ and let $\sigma_{\fa}\in G$ be the scaling matrix such that $\tau_{\fa}\infty=\fa$ and
%$$\tau_{\fa}^{-1}\G\tau_{\fa}\cap N=\langle u_1\rangle.$$
%\subsubsection{Hecke operators and spectral decomposition on the modular surface}
%\{\pm\left(\begin{smallmatrix}
%1 & n\\
%0 & 1\end{smallmatrix}\right):n\in\Z\}\leq \G$ be the stabilizer of $\infty$. 
Let now $\Gamma=\SL_2(\Z)$. In the notation introduced at the beginning of this section, the Eisenstein series for the modular group $\Gamma$ at the cusp $\infty$ is defined for $\Re(s)>1$ by
\begin{equation}\label{equ:eisen}
E(z,s)=\sum_{\gamma\in (\G\cap \pm N)\bk \G}\Im(\gamma z)^s
\end{equation}
with a meromorphic continuation to $s\in \C$. Moreover, for any $s\in \C$, $E(\cdot, s)$ is an eigenfunction of the Laplace operator with eigenvalue $s(1-s)$.

Let $\Psi\in L^2(\cM)$ and we have the following spectral decomposition (see \cite[Theorems 4.7 and 7.3]{Iwaniec2002})
\begin{equation}\label{equ:specdec}
\Psi(z)=\mu_{\cM}(\Psi)+\sum_{r_k\geq 0}\langle \Psi,\phi_k\rangle \phi_k(z)+\frac{1}{4\pi}\int_{-\infty}^{\infty}\langle \Psi, E(\cdot, \tfrac12+ir)\rangle E(z,\tfrac12+ir)dr,
\end{equation}
where the convergence holds in the $L^2$-norm topology, and is pointwise if $\Psi\in C_c^{\infty}(\cM)$. As a direct consequence we have for $\Psi\in  L^2(\cM)$,
\begin{equation}\label{equ:specde2}
\|\Psi\|_2^2=\left|\mu_{\cM}(\Psi)\right|^2+\sum_{r_k\geq 0}\left|\langle\Psi,\phi_k\rangle\right|^2+\frac{1}{4\pi}\int_{-\infty}^{\infty}\left|\langle\Psi,E(\cdot,\tfrac12+ir)\rangle\right|^2dr.
\end{equation}

\subsection{Hecke operators}\label{sec:hecopr}
The spectral theory of $\cM$ has extra structure due to the existence of \textit{Hecke operators}. The main goal of this subsection is to prove an operator norm bound for Hecke operators and the main reference is \cite[Section 8.5]{Iwaniec2002}. For any $n\in\N$ define the set
\begin{equation}\label{equ:hecdef1}
\cL_n:=\left\{n^{-1/2}g :g\in M_2(\Z),\; \det(g)=n\right\} \subset G,
\end{equation}
where $M_2(\Z)$ is the space of two by two integral matrices. The \textit{$n$-th Hecke operator $T_n$} is defined by that for any $\Psi\in L^2(\cM)$
$$T_{n}(\Psi)(z)=\frac{1}{n^{1/2}}\sum_{\g\in \G\bk \cL_n}\Psi(\g z).$$ 

The Hecke operator $T_n$ is a self-adjoint operator on $L^2(\cM)$ and since $T_n$ commutes with the Laplace operator $\Delta$ (since $\Delta$ is defined via right multiplication and $T_n$ is defined via left multiplication) the orthonormal basis of the space of cusp forms $\{\phi_k\}$ can be chosen consisting of joint eigenfunctions of all $T_n$, that is,
$$T_n \phi_k=\lambda_{\phi_k}(n)\phi_k.$$
On the other hand, for any $r\in\R$ the Eisenstein series $E(z, 1/2+ir)$ is an eigenfunction of $T_n$ with eigenvalue $\lambda_r(n):=\sum_{d\mid n}\left(\frac{n}{d^2}\right)^{ir}$, see \cite[Equation (8.33)]{Iwaniec2002}. It is clear that $\left|\lambda_r(n)\right|\leq \sigma_0(n)$ with $\sigma_0(n)$ the divisor function. For the eigenvalue of cusp forms it is conjectured (Ramanujan-Petersson) that for any above $\phi_k$ and for any $n\in\N$
$$\left|\lambda_{\phi_k}(n)\right|\leq \sigma_0(n).$$
The aforementioned bound of Sarnak and Kim \cite{KimSarnak03} implies that
$$\left|\lambda_{\phi_k}(n)\right|\leq \sigma_0(n)n^{7/64}.$$
Using these bounds on eigenvalues and the above spectral decomposition \eqref{equ:specdec} and \eqref{equ:specde2} we have the following bound on the operator norm of the Hecke operator, see also \cite[pp. 172-173]{GoldsteinMayer2003}.
\begin{Prop}\label{prop:heckeclass}
For any $\Psi\in L^2(\cM)$ and for any $n\in\N$ we have
$$\langle \Psi_0, T_n(\Psi_0)\rangle_{L^2(\cM)}\ll_{\e}n^{\theta+\e}\|\Psi\|_2^2,$$
where $\Psi_0:=\Psi-\mu_{\cM}(\Psi)$ and $\theta=7/64$ as before.
\end{Prop}

\subsubsection{Hecke operators attached to a group element}
Let $\G=\SL_2(\Z)$ and let $\cM=\G\bk\bH$ be the modular surface as above. There is another type of Hecke operators on $L^2(\cM)$ defined via a group element in $\SL_2(\Q)$. Namely, for each $h\in \SL_2(\Q)$ the \textit{Hecke operator attached to $h$}, denoted by $\widetilde{T}_h$, is defined by that for any $\Psi\in L^2(\cM)$
\begin{equation}\label{equ:hecdef2}
\widetilde{T}_h(\Psi)(z)=\frac{1}{\# (\G\bk \G h\G)}\sum_{g\in \G\bk \G h\G}\Psi(g z),
\end{equation}
where $\G h\G=\left\{\gamma_1 h\gamma_2:\gamma_1,\gamma_2\in\G\right\}$ is the double coset attached to $h$. We note that $\widetilde{T}_h$ is well-defined since $\Psi$ is left $\G$-invariant. 

For our purpose, we will need another expression for $\widetilde{T}_h$. For any $h\in \SL_2(\Q)$ we denote by $\G^h:=\G\cap h^{-1}\G h$. We note that the map from $\G$ to $\G\bk \G h\G$ sending $\g\in \G$ to $\G h\g$ induces an identification between $\G^h\bk \G$ and $\G\bk \G h\G$. This identification induces the following alternative expression for $\widetilde{T}_h$:
\begin{equation}\label{equ:Hecke}
\widetilde{T}_h(\Psi)(g)=\frac{1}{[\G: \G^h]}\sum_{\g\in \G^h\bk \G}\Psi(h\g g).
\end{equation}

It is clear from the definition that $\widetilde{T}_h$ is defined only up to representatives for the double coset $\G h\G$, that is, $\widetilde{T}_h=\widetilde{T}_{h'}$ whenever $\G h\G=\G h'\G$. For a fixed $h\in\SL_2(\Q)$, we call $n\in\N$ \textit{the degree of $h$} if $n$ is the smallest positive integer such that $nh\in M_2(\Z)$. Using elementary column and row operations one can see that for $h\in\SL_2(\Q)$ with degree $n$
\begin{equation}\label{equ:heckeorb}
\G h\G=\G \diag (1/n,n)\G=\left\{n^{-1}g : g\in M_2(\Z),\; det(g)=n^2,\; \gcd (g)=1\right\} {{\subset G}},
\end{equation}
where $\gcd(g)$ is the greatest common divisor of the entries of $g$. Thus we can parameterize the Hecke operators by their degrees, that is, we will denote by $\widetilde{T}_n:=\widetilde{T}_h$ for any $h\in \SL_2(\Q)$ with degree $n$. We also note that by direct computation when $h=\diag (1/n,n)$ we have $\G^h=\G_0(n^2)$, implying that for any $h\in \SL_2(\Q)$ with degree $n$ (see e.g. \cite[Section 1.2]{DiamondShurman2005})
\begin{equation}\label{equ:indexhec}
\nu_n:=\# (\G\bk \G h\G)=[\G : \G^h]=[\G : \G_0(n^2)]=n^2\prod_{\substack{p\mid n\\ \textrm{prime}}}\left(1+p^{-1}\right).
\end{equation}
Now using the description \eqref{equ:heckeorb} we have the double coset decomposition
$$\cL_{n^2}=\bigsqcup_{d\mid n}\G\begin{pmatrix}
d^{-1} & 0\\
0 & d\end{pmatrix}\G.$$
This decomposition together with the definitions \eqref{equ:hecdef1}, \eqref{equ:hecdef2} and \eqref{equ:indexhec} implies the relation
%recall the Hecke points $\cM_n\subset Z\G\bk \GL_2^+(\R)$ is defined the finite set of sub-lattices of $\Z^2$ of index $n$. 
%we can view a function on $\G\bk G$ as a left $Z\G$-invariant function on $\GL_2^+(\R)$. 
%
%it has a natural action on $X$ (and hence also on $L^2(\G\bk G)$) by identifying $X$ as the space of rank two lattices up to scales. More precisely, $g\in \GL_2^+(\R)$ acts on $\Psi\in L^2(\G\bk G)$ via $g\cdot \Psi (x)=\Psi(x\tilde{g})$ with $\tilde{g}=\left(\det(g)\right)^{-1/2}g\in G$. For any positive integer $n$ let 
%$$\cM_n:=\left\{g\in M_2(\Z):\det(g)=n\right\}$$
%be the finite set of integral two by two matrices with determinant equalling $n$. Recall that the classical Hecke operator $\widetilde{T}_n$ is defined such that 
$$nT_{n^2}=\sum_{d| n}\nu_d\widetilde{T}_{d}.$$
Thus by the M\"obius inversion formula we have
\begin{equation}\label{equ:heckerelation}
\widetilde{T}_{n}=\frac{n}{\nu_n}\sum_{d\mid n}\frac{\mu(d)}{d}T_{n^2/d^2}.
\end{equation}
Using this relation and \propref{prop:heckeclass} we can prove the following operator norm bounds for $\widetilde{T}_n$ which we will later use, see also \cite[Theorem 1.1]{ClozelOhUllmo2001} for such bounds in a much greater generality.
\begin{Prop}\label{prop:hecbound2}
Keep the notation as in \propref{prop:heckeclass}. For any $\Psi\in L^2(\cM)$ and for any $n\in\N$ we have
$$\langle \Psi_0, \widetilde{T}_n(\Psi_0)\rangle_{L^2(\cM)}\ll_{\e}n^{-1+2\theta+\e}\|\Psi\|_2^2.$$
\end{Prop}
\begin{proof}
By \propref{prop:heckeclass} and using the relation \eqref{equ:heckerelation}, the trivial estimates $|\mu(d)|\leq 1$ and $\nu_n\geq n^2$ and the triangle inequality we have
\begin{align*}
\langle \Psi_0, \widetilde{T}_n(\Psi_0)\rangle&\leq n^{-2}\sum_{d| n}(n/d)\langle \Psi_0, T_{n^2/d^2}(\Psi_0)\rangle\ll_{\e} n^{-2}\sum_{d\mid n}(n/d)^{1+2\theta+2\e}\|\Psi\|_2^2\\
&=n^{-1+2\theta+2\e}\sigma_{-1+2\theta+2\e}(n)\|\Psi\|_2^2\ll_{\e}n^{-1+2\theta+\e}\|\Psi\|_2^2.\qedhere
\end{align*}
\end{proof}

\subsection{Equidistribution of subsegments of expanding closed horocycles}\label{sec:longhoro}
%Let $\G< \SL_2(\Z)$ be a congruence subgroup and 
%that is $\G$ contains a \textit{principal congruence subgroup} 
%$$\G(n):=\left\{\gamma\in \SL_2(\Z):\gamma\equiv I_2 \Mod{n}\right\}$$
%of some level $n\in \N$. 
%In particular $\G$ has a cusp at $\infty$, and we 
%assume that the width of the cusp at $\infty$ is one, that is, $u_1\in \G$. %Since $u_1\in \G$, for any $y>0$ the projection of the horizontal horocycle $\{z\in \bH:\Im(z)=y\}$ on $\G\bk \bH$ is closed and is of hyperbolic length $1/y$. Similar as in the case on the modular surface, these closed horocycles 
%$$\G Na_y=\{\G u_ta_y:t\in [0,1)\}$$ 
%become equidistributed on $\G\bk \bH$ with respect to $\mu_{\G}$ as $y\to 0^+$, see \cite{Sarnak1981} for an effective equidistribution on a general non-uniform lattice using spectral arguments. %Moreover, we record an effective equidistribution result for congruence subgroups which is due to Str\"{o}mbergsson \cite[Theorem 1 and Remark 3.4]{Strombergsson2013}: there exists a constant $C_{\G}$ depending only on $\G$ such that for any $\Psi\in C^{\infty}(\G\bk G)$ with $\cS_{2,4}^{\G}(\Psi)<\infty$ and for any $0<y<1$
We record a special case of Sarnak's result \cite[Theorem 1]{Sarnak1981} on effective equidistribution of expanding closed horocycles, %This restriction to congruence subgroups gives us the advantage of not to consider the nontrivial residual spectrum, i.e. the contribution from the exceptional poles of the Eisenstein series (since in this case they do not exist, see \cite[Theorem 11.3]{Iwaniec2002}). As a result we can make the dependence on the testing function explicit. 
namely:
\begin{Prop}\label{prop:equiclohor}
Let $\G<\SL_2(\Z)$ be a congruence subgroup and assume that $\G$ has a cusp at $\infty$ with width one. Then for any $\Psi\in C^{\infty}(\G\bk \bH)\cap L^2(\G\bk \bH)$ satisfying $\|\Delta \Psi\|_2<\infty$ and for any $0<y<1$ we have 
\begin{equation}\label{equ:effequ}
\left|\int_0^1\Psi(x+ iy)dx-\mu_{\G}(\Psi)\right|\ll \|\Psi\|_2^{3/4}\|\Delta\Psi\|_2^{1/4}y^{1/2},
\end{equation}
where the implied constant is absolute, independent of $\G$, $\Psi$ and $y$, and the $L^2$-norm is with respect to the normalized hyperbolic area $\mu_{\G}$.
\end{Prop}
\begin{rmk}
We omit the proof here and refer the reader to \cite[(3.5)]{KelmerKontorovich2020}. We note that while \cite{KelmerKontorovich2020} only deals with the case when $\G=\G_0(p)$ with $p$ a prime number, 
%(taking the cusp $\fa=\infty$ and the scaling matrix $\tau_{\fa}=I_2$ in their setting), 
the proof there works for general congruence subgroups, given that they have trivial residual spectrum; see \cite[Theorem 11.3]{Iwaniec2002}.
%We note that their proof uses only some general theory of Eisenstein series, more specifically, the unitarity of the scattering matrix on the critical line, and can be easily generalized to a general congruence subgroup .
\end{rmk}

%\begin{equation}\label{equ:effequ}
%\left|\int_0^1\Psi(u_t a_y)dt-\mu_{\G}(\Psi)\right|\leq C_{\G}\cS_{2,4}^{\G}(\Psi)y^{1/2-\theta},
%\end{equation}
%where $\theta=7/64$ is the best known bound towards the Ramanujan conjecture obtained by Kim and Sarnak \cite{KimSarnak03} and the bounding constant depends only on $\G$. 
{{We will also need the following (non-effective) equidistribution result replacing the whole closed horocycle by a fixed subsegment: 
\begin{Prop}\label{prop:noneffequ}
Let $\G<\SL_2(\Z)$ be as in \propref{prop:equiclohor}. Let $I\subset (0,1)$ be an open interval, then for any $\Psi\in C_c(\G\bk \bH)$ we have
\begin{equation}\label{equ:equipie}
\lim\limits_{y\to 0^+}\frac{1}{|I|}\int_{I}\Psi(x+iy)dx=\mu_{\G}(\Psi).
\end{equation}
\end{Prop}
The proof of \propref{prop:noneffequ} uses}}
Margulis' thickening trick \cite{Margulis2004} and mixing {{property of the geodesic flow on the unit tangent bundle of $\G\bk \bH$};}
%one can prove a more general equidistribution result replacing the whole closed horocycle by a fixed subsegment, see \cite{EskinMargulisMozes1998}: 
%Let $I\subset (0,1)$ be an open interval, then for any $\Psi\in C_c(\G\bk \bH)$ we have
%\begin{equation}\label{equ:equipie}
%\lim\limits_{y\to 0^+}\frac{1}{|I|}\int_{I}\Psi(x+iy)dx=\mu_{\G}(\Psi).
%\end{equation}
this approach is also effective, see e.g. \cite[Proposition 2.3]{KelmerKontorovich2018}. A proof of \eqref{equ:equipie} using spectral methods was also sketched in \cite[Theorem 1$'$]{Hejhal1996}. 
{{
We also note that both equidistribution results in \propref{prop:equiclohor} and \propref{prop:noneffequ} can be lifted to the unit tangent bundle of $\G\bk \bH$ (with necessary modifications to the error term in \eqref{equ:effequ}); since we will be only working in the hyperbolic surface level, we state these two results in the current format for convenience of our discussion.
}}We further refer the reader to \cite{Hejhal2000a,Strombergsson2004} for some much stronger effective equidistribution results regarding long enough (varying) subsegments on expanding closed horocycles. 
{{
\begin{remark}\label{rmk:apparg}
\propref{prop:noneffequ} can be equivalently stated as following: For any fixed open interval $I\subset (0,1)$, the measures $\mu_{I,y}$ weakly converge to $\mu_{\G}$ as $y\to 0^+$, where for any $y\in (0,1)$ and $\Psi\in C_c(\G\bk \bH)$, $\mu_{I,y}(\Psi):=\frac{1}{|I|}\int_I\Psi(x+iy)dx$. Thus by the Portmanteau theorem, \eqref{equ:equipie} extends to $\Psi=\chi_B$ with $B\subset \G\bk\bH$ a Borel subset with boundary of measure zero. 
More generally, let $\rho: [0, 1)\to \R$ be a Riemann integrable function. Since $\rho$ can be weakly approximated from both above and below by step functions, we have
% and let $\Psi$ be a non-negative function on $\G\bk\bH$ such that there exists a sequence $\{\Psi_j^{\pm}\}_{j\in\N}\subset C_c(\G\bk \bH)$ satisfying $\Psi_j^-\leq \Psi\leq \Psi_j^+$ for every $j\in\N$ and $\lim\limits_{j\to\infty}\mu_{\G}\left(\Psi_j^{\pm}\right)=\mu_{\G}(\Psi)$. Then \propref{prop:noneffequ}, together with a standard approximation argument implies that
\begin{align*}
\lim\limits_{y\to 0^+}\int_{0}^1\rho(x)\chi_B(x+iy)dx=\mu_{\G}(B)\int_0^1\rho(x)dx
\end{align*}
with $B\subset \G\bk\bH$ a Borel set with boundary of measure zero.
\end{remark}
}}

\subsection{A quantitative Borel-Cantelli lemma}
Finally we record here a quantitative Borel-Cantelli lemma which ensures for the limsup set of certain sequence of events to have full measure given certain quasi-independence conditions.
\begin{Lem}{\cite[Chapter I, Lemma 10]{Sp79}}\label{lem:bcd}
Let $(X,\cB, \nu)$ be a probability space with $\cB$ a $\sigma$-algebra of subsets of $X$ and $\nu: X\to [0,1]$ a probability measure on $X$ with respect to $\cB$. Let $\{A_i\}_{i\in\N}$ be a sequence of measurable subsets in $\cB$. For any $n,m\in\N$ we denote by $R_{n,m}:=\nu(A_n\cap A_m)-\nu(A_n)\nu(A_m)$. Suppose that
\begin{equation}\label{equ:qbccond}
\textrm{$\exists\ C>0$ such that for all $k_2>k_1\geq 1$, $\sum_{n,m=k_1}^{k_2}R_{n,m}\leq C\sum_{n=k_1}^{k_2}\nu(A_n)$,}
\end{equation}
then $\sum_{n\in\N}\nu(A_n)=\infty$ implies that $\nu\left(\limsup_{n\to\infty}A_n\right)=1$.
\end{Lem}
\begin{rmk}\label{rmk:furcond}
Keep the notation as in \lemref{lem:bcd}. It was shown in \cite[Proposition 5.4]{KelmerYu2019a} that if
$$\textrm{$\exists C'>0$ and $\eta>1$ such that for any $n\neq m$, $R_{n,m}\leq C'\frac{\sqrt{\nu(A_n)\nu(A_m)}}{|n-m|^{\eta}}$},$$
then the sequence $\{A_i\}_{i\in\N}$ satisfies the condition \eqref{equ:qbccond}.
%In particular, using the trivial estimates $\nu(A_i)\leq 1$ and $|i-j|< \max\{i,j\}$, the condition 
%\begin{equation}\label{equ:quasicond2}
%\textrm{$\exists C'>0$ and $\eta>1$ such that for any $1\leq i< j$, $R_{i,j}\leq C'\frac{\nu(A_i)\nu(A_j)}{ j^{\eta}}$}
%\end{equation}
%$$\forall 1\leq i<j, R_{i,j}\leq C'\frac{\mu{A_i}\mu(A_j)}{ j^{\eta}}$$
%also implies \eqref{equ:qbccond}.
\end{rmk}
We will use the following slightly modified version of quantitative Borel-Cantelli lemma which has the flexibility to consider sequence of measurable sets $\{A_n\}_{n\in\mathbb{S}}$ indexed by a general unbounded subset $\mathbb{S}\subset \N$. 
\begin{Cor}\label{cor:quanbc}
Let $(X,\cB,\nu)$ be as in \lemref{lem:bcd}. Let $\mathbb{S}\subset\N$ be an unbounded subset and let $\{A_n\}_{n\in\mathbb{S}}$ be a sequence of measurable subsets in $\cB$. Suppose that
\begin{equation}\label{equ:quasicond2}
\textrm{$\exists C'>0$ and $\eta>1$ such that $\forall$ $n,m\in \mathbb{S}$ with $m< n$, $R_{n,m}\leq C'\frac{\nu(A_n)\nu(A_m)}{n^{\eta}}$,}
\end{equation}
then $\sum_{n\in\mathbb{S}}\nu(A_n)=\infty$ implies that $\nu\left(\limsup_{\substack{n\in \mathbb{S}\\ n\to\infty}}A_n\right)=1$.
\end{Cor}
\begin{proof}
%We first prove the case when $\mathbb{S}=\N$ is the whole set of natural numbers. In view of \lemref{lem:bcd} it suffices to show that condition \eqref{equ:quasicond2} implies the condition \eqref{equ:qbccond}. In fact, it was shown in \cite[Proposition 5.4]{KelmerYu2019a} that if
%$$\textrm{$\exists C'>0$ and $\eta>1$ such that for any $n\neq m\in \N$, $R_{n,m}\leq C'\frac{\sqrt{\nu(A_n)\nu(A_m)}}{|n-m|^{\eta}}$},$$
%then the sequence $\{A_n\}_{n\in\N}$ satisfies the condition \eqref{equ:qbccond}. Using the trivial estimates $\nu(A_n)\leq 1$ and $|n-m|< \max\{n,m\}$, we have for any $1\leq m< n$
%$$R_{n,m}<C'\frac{\nu(A_n)\nu(A_m)}{n^{\eta}}< C'\frac{\sqrt{\nu(A_n)\nu(A_m)}}{|n-m|^{\eta}}.$$
%This finishes the proof for the case when $\mathbb{S}=\N$. For a general unbounded subset $\mathbb{S}$, 
For any $i\in\N$ let $a_i\in \mathbb{S}$ be the $i$-th integer in $\mathbb{S}$ and let $B_i:=A_{a_i}$. For any $i,j\in\N$ let $R'_{i,j}:=\nu(B_i\cap B_j)-\nu(B_i)\nu(B_j)$ so that $R_{i,j}'=R_{a_i,a_j}$. Then by for any $i<j$ we have
%Using the relation $B_i=A_{a_i}$ and the trivial estimates $a_i\geq i$ we have for any $i<j$
$$R'_{i,j}=R_{a_i,a_j}\leq C'\frac{\nu(A_{a_i})\nu(A_{a_j})}{a_j^{\eta}}= C'\frac{\nu(B_i)\nu(B_j)}{a_j^{\eta}}<C'\frac{\sqrt{\nu(B_i)\nu(B_j)}}{|i-j|^{\eta}},$$
where for the first inequality we used the assumption \eqref{equ:quasicond2} and for the second inequality we used the estimates $a_j\geq j> j-i$ and $\sqrt{\nu(B_i)\nu(B_j)}\leq 1$. Thus in view of \rmkref{rmk:furcond} and \lemref{lem:bcd} we have $\sum_{i\in\N}\nu(B_i)=\infty$ implies that $\nu\left(\limsup_{i\to\infty}B_i\right)=1$ which is equivalent to the conclusion of this corollary in view of the relation $B_i=A_{a_i}$.
%That is, the sequence $\{B_i\}_{i\in\N}$ satisfies the condition \eqref{} (for the whole set of natural numbers). Hence by the 
\end{proof}
\section{Equidistribution range}\label{sec:equran}
%Let $\G=\SL_2(\Z)$ and let $\cM=\G\bk\bH$ be the modular surface. Let $\mu_{\cM}$ be the normalized hyperbolic area on $\cM$. 
Let $\cM=\SL_2(\Z)\bk\bH$. Since we fix $\G=\SL_2(\Z)$ throughout this section, we abbreviate the Sobolev norm $\cS_{p,d}^{\G}$ by $\cS_{p,d}$. In this section we prove \thmref{thm:equ} and \thmref{thm:fullrange}. The main ingredient of our proof is an explicit bound of Fourier coefficients which follows from a slight modification of the estimates obtained in \cite{KelmerKontorovich2020}. 

\subsection{Bounds on Fourier coefficients}\label{sec:fouriercoe}
Let $\Psi\in C_c^{\infty}(\cM)$. Since $\Psi$ is left $\G$-invariant, it is invariant under the transformation determined by $u_1:z\mapsto z+1$, and it thus has a Fourier expansion for $\Psi$ in the variable $x=\Re(z)$: %with respect to the variable $x$:
\begin{equation}\label{equ:foex}
\Psi(x+iy)=\sum_{m\in\Z}a_{\Psi}(m,y)e(mx),
\end{equation}
where
$$a_{\Psi}(m,y)=\int_0^1\Psi(x+iy)e(-mx)dx.$$
%Let $G=\PSL_2(\R)$ and it can be realized as the orientation preserving isometry group of the upper half plane $\bH:=\{z=x+iy\in \C :y>0\}$ via the M\"obius transformation. Let $\G=\PSL_2(\Z)$ and let $X=\G\bk \bH$ be the modular surface and it can be identified with the locally symmetric space $\G\bk G/\operatorname{PSO}_2(\R)$. 
%
%
%
%Let $d\mu(z)=\frac{3}{\pi}\frac{dxdy}{y^2}$ be the normalized hyperbolic area on $X$. Let $L^2(X)$ be the space of square-integrable automorphic functions (with respect to $\mu$), and it can be identified with the subspace of right $\operatorname{PSO}_2(\R)$-invariant functions in $L^2(\G\bk G)$. 
Similarly we denote by $a_{\phi_k}(m,y)$ and $a(s;m,y)$ the $m$th Fourier coefficients of the Hecke-Maass form $\phi_k$ and the Eisenstein series $E(\cdot, s)$ respectively. Estimates on these Fourier coefficients yield, via the spectral expansion \eqref{equ:specdec}, estimates on the Fourier coefficients of $\Psi$. Namely, 
$$a_{\Psi}(m,y)=\sum_{r_k\geq 0}\langle \Psi,\phi_k\rangle a_{\phi_k}(m,y)+\frac{1}{4\pi}\int_{-\infty}^{\infty}\langle \Psi, E(\cdot, \tfrac12+ir)\rangle a(\tfrac12+ir;m,y)dr.$$

%For any $m\in \Z$ and $y>0$ let 
%$$a_{\phi_k}(m,y)=\int_0^1\phi_k(x+iy)e(-mx)dx\quad \textrm{and}\quad a(s;m,y)=\int_0^1E(x+iy,s)e(-mx)dx$$ 
%be the $m$th Fourier coefficients of $\phi_k$ and $E(z,s)$ respectively. 
%\section{Estimates on Fourier coefficients} 
We record the following bounds for $a_{\phi_k}(m,y)$ and $a(s;m,y)$: 
\begin{Lem}[{{\cite[Lemmata 3.7 and 3.13]{KelmerKontorovich2020}}}]\label{lem:fe1}
For any $m\neq 0$ and for any $\e>0$ we have
\begin{equation}\label{equ:b1}
|a_{\phi_k}(m,y)|\ll_{\e}|m|^{\theta}y^{1/2-\e}(r_k+1)^{-1/3+\e}\min\{1,e^{\pi r_k/2-2\pi|m|y}\},
\end{equation}
and
\begin{equation}\label{equ:b2}
|a\left(\tfrac12+ir;m,y\right)|\ll_{\e}y^{1/2-\e}(1+|r|)^{-1/3+\e}\min\{1,e^{\pi |r|/2-2\pi|m| y}\},
\end{equation}
where $\theta=7/64$ is the best known bound towards the Ramanujan conjecture as before.
\end{Lem}
\begin{rmk}
Contrarily to \cite{KelmerKontorovich2020} that uses the trivial bound $\min\{1,e^{\pi r/2-2\pi|m|y}\}\leq 1$, we keep this term.
\end{rmk}

%\begin{proof}
%We note that \eqref{equ:b1} and \eqref{equ:b2} follow from the analysis in \cite[Lemma 3.7]{KelmerKontorovich2020} and \cite[Lemma 3.13]{KelmerKontorovich2020} respectively. The only difference is that in \cite{KelmerKontorovich2020} they used the trivial bound $\min\{1,e^{\pi r/2-2\pi|m|y}\}\leq 1$ while here we keep this term.
%\end{proof}

\begin{Prop}[{{\cite[Proposition 3.4]{KelmerKontorovich2020}}}]
For any $\Psi\in C_c^\infty(\cM)$, we have that 
\begin{equation}\label{equ:zerocoef}
a_{\Psi}(0,y)\ =\ \mu_{\mathcal{M}}(\Psi) + O\left(\|\Psi\|_2^{3/4}\|\Delta\Psi\|_2^{1/4}y^{1/2}\right).
\end{equation}
Moreover, for any $m\neq0$, and any $\epsilon>0$ and any $\alpha_0>5/3$, we have
\begin{align}\label{equ:smallm}
a_{\Psi}(m,y)\ \ll_{\alpha_0,\epsilon,p} 
\mathcal{S}_{\alpha_0}(\Psi)y^{1/2-\epsilon}|m|^\theta,
\end{align}
where $\mathcal{S}_{\alpha_0}$ is a Sobolev norm of degree $\alpha_0$. %Moreover, for each of the exceptional forms $a_{\phi_k,\frak{a}}(m,y)=O_\epsilon\left(|m|^{\theta-|r_k|+\epsilon}y^{1/2-|r_k|}e^{-2\pi |m|y}\right).$
\end{Prop}
\begin{rmk}\label{rmk:sobolev}
The Sobolev norm $\cS_{\alpha_0}$ is explicit from the proof of \cite[Proposition 3.4]{KelmerKontorovich2020}: Writing $\alpha_0=5/3+\e$ with $\e>0$, then $\cS_{\alpha_0}(\Psi)=\cS_{2,0}(\Psi)^{2/3-\e/2}\cS_{2,2}(\Psi)^{1/3+\e/2}$ for any $\Psi\in C_c^{\infty}(\cM)$. In particular, using the estimate $\cS_{2,0}(\Psi)\leq \cS_{2,2}(\Psi)$ we have $\cS_{\alpha_0}(\Psi)\leq \cS_{2,2}(\Psi)$. 
\end{rmk}
The following refinement of this last estimate allows to estimate the Fourier coefficients when $|m|>y^{-1}$ is large. This refinement is crucial for our later results, and the price we pay is a Sobolev norm of higher degree.

\begin{Prop}\label{prop:fe2}
Let $\Psi\in C_c^{\infty}(\cM)$. Whenever $|m|y>1$ and for any $\epsilon>0$, we have
$$|a_{\Psi}(m,y)|\ll_{\e}\cS_{2,2}(\Psi) |m|^{-4/3+\theta+\e}y^{-5/6}.$$
%$$|a_{\Psi}(m,y)|\ll_{\e}\|\Psi\|_2|m|^{\theta+\frac23+\frac{3}{2}\e}y^{\frac76+\frac{\e}{2}}e^{-\pi |m|y}+\|\Delta\Psi\|_2 |m|^{\theta-\frac43+\frac{3}{2}\e}y^{-\frac56+\frac{\e}{2}}.$$
\end{Prop}
\begin{proof}
For the contribution from the cusp forms we
%\begin{align*}
%\left|\sum_{r_k\geq 0}\langle \Psi,\phi_k\rangle a_{\phi_k}(m,y)\right|&\leq\left(\sum_{0\leq r_k\leq 2|m|y}+\sum_{r_k>2|m|y}\right)\left|\langle \Psi,\phi_k\rangle a_{\phi_k}(m,y)\right|.
%\end{align*}
apply the bound \eqref{equ:b1} to the Fourier coefficients and the bound 
\begin{equation}\label{equ:prebound}
\min\{1,e^{\pi r/2-2\pi|m|y}\}\leq \left\lbrace\begin{array}{ll} 
e^{-\pi|m|y} & 0\leq r\leq 2|m|y\\
1 & r> 2|m|y,
\end{array}\right.
\end{equation}
and the relation $\langle \Delta\Psi, \phi_k\rangle=\langle \Psi, \Delta\phi_k\rangle=(1/4+r_k^2)\langle \Psi, \phi_k\rangle$ to get that
\begin{align}\label{equ:cusppb}
\left|\sum_{r_k\geq 0}\langle \Psi,\phi_k\rangle a_{\phi_k}(m,y)\right|&\ll_{\e} \sum_{0\leq r_k\leq 2|m|y}\left|\langle \Psi, \phi_k\rangle\right| |m|^{\theta}y^{1/2-\e}(r_k+1)^{-1/3+\e}e^{-\pi |m|y}\\
&+\sum_{r_k>2|m|y}\left|\langle \Delta\Psi,\phi_k\rangle\right||m|^{\theta}y^{1/2-\e}r_k^{-7/3+\e}.\nonumber
\end{align}
Now using Cauchy-Schwarz followed by summation by parts (together with Weyl's law stating that $\#\{r_k:r_k\leq M\}\ll M^2$ (see e.g. \cite[Corollary 11.2]{Iwaniec2002}) we can bound
\begin{align*}
\sum_{0\leq r_k\leq 2|m|y}\left|\langle \Psi, \phi_k\rangle\right|(r_k+1)^{-1/3+\e}&\leq \|\Psi\|_2\left(\sum_{0\leq r_k\leq 2|m|y}\frac{1}{(r_k+1)^{2/3-2\e}}\right)^{1/2}\\
&\ll_{\e} \|\Psi\|_2\left(|m|y\right)^{2/3+\e}. 
\end{align*}
Similarly, for the second sum we can bound
$$\sum_{r_k>2|m|y}\left|\langle \Delta\Psi,\phi_k\rangle\right|r_k^{-7/3+\e}\leq \|\Delta\Psi\|_2\left(\sum_{r_k>2|m|y}r_k^{-14/3+2\e}\right)^{1/2}\ll_{\e}\|\Delta\Psi\|_2\left(|m|y\right)^{-4/3+\e}.$$
To summarize, the left-hand side of \eqref{equ:cusppb} is bounded by
\begin{equation}\label{equ:cuspb}
\ll_\epsilon\ \|\Psi\|_2 |m|^{2/3+\theta+\epsilon} y^{7/6} e^{-\pi|m|y} + \|\Delta\Psi\|_2 |m|^{-4/3+\theta+\epsilon} y^{-5/6}.
%\left|\sum_{r_k\geq 0}\langle \Psi,\phi_k\rangle a_{\phi_k}(m,y)\right|\ll
%\ll_{\e}|m|^{\theta}y^{1/2-\e}\left(\|\Psi\|_2\left(|m|y\right)^{2/3+\frac32\e}e^{-\pi|m|y}+\|\Delta\Psi\|_2\left(|m|y\right)^{-4/3+\frac32\e}\right).
\end{equation}
For the contribution from the continuous spectrum using the estimates \eqref{equ:b2}, \eqref{equ:prebound}, the relation 
$\langle \Delta\Psi,E(\cdot, \tfrac12+ir)\rangle=(\tfrac14+r^2)\langle\Psi, E(\cdot,\tfrac12+ir)\rangle$ and Cauchy-Schwarz
we can similarly bound $\left|\int_{-\infty}^{\infty}\langle \Psi, E(\cdot, \tfrac12+ir)\rangle a(\tfrac12+ir;m,y)dr\right|$ by
\begin{align*}
&\ll_{\e}e^{-\pi |m|y}y^{1/2-\e}\int_{|r|\leq 2|m|y}\left|\langle \Psi, E\left(\cdot, \tfrac12+ir\right)\rangle\right|(|r|+1)^{-1/3+\e}dr\\
&+y^{1/2-\e}\int_{|r|>2|m|y}\left|\langle\Delta\Psi, E\left(\cdot, \tfrac12+ir\right)\rangle\right||r|^{-7/3+\e}dr\\
&\ll_{\e}y^{1/2-\e}\left(\|\Psi\|_2\left(|m|y\right)^{1/6+\e}e^{-\pi|m|y}+\|\Delta\Psi\|_2\left(|m|y\right)^{-11/6+\e}\right),
\end{align*}
which is subsumed by the right-hand side of \eqref{equ:cuspb} (since $|m|y>1$). Finally, we conclude the proof by applying the bounds $\max\{\|\Psi\|_2, \|\Delta\Psi\|_2\}\leq \cS_{2,2}(\Psi)$ and $e^{-\pi |m|y}\ll (|m|y)^{-2}$ (again since $|m|y>1$) to the right hand side of \eqref{equ:cuspb}.
\end{proof}

%\begin{rmk}\label{rmk:smallm}
%\propref{prop:fe2} is a slight modification to \cite[(3.6)]{KelmerKontorovich2020} (with a different choice of parameter splitting the sum in \eqref{equ:cusppb}). We will use \propref{prop:fe2} to estimate the Fourier coefficients when $|m|> y^{-1}$ is large. Unconditionally, it was proved in \cite{KelmerKontorovich2020} that for any $|m|\neq 0$ and for any $0<y<1$ (noting that there are no exceptional cusp forms for $\G=\SL_2(\Z)$ since the Selberg's eigenvalue conjecture holds for $\G$)
%\begin{equation}\label{equ:smallm}
%\left|a_{\Psi}(m,y)\right|\ll_{\e} \cS_{2,0}(\Psi)^{2/3-\e/2}\cS_{2,2}(\Psi)^{1/3+\e/2}|m|^{\theta}y^{1/2-\e}. %\leq \cS_{2,2}(\Psi)|m|^{\theta}y^{1/2-\e}.
%\end{equation}
%When $|m|>y^{-1}$, \propref{prop:fe2} reads as $a_{\Psi}(m,y)=O_{\e, \Psi}\left((|m|y)^{-4/3+3\e/2} |m|^{\theta}y^{1/2-\e}\right)$, having an extra $(|m|y)^{-4/3+3\e/2}$ saving compared to \eqref{equ:smallm}. We note that this extra saving is crucial for our proof, and the price we pay here is that the Sobolev norm in \propref{prop:fe2} is of higher degree than that of \eqref{equ:smallm}.
%\end{rmk}

The following corollary of \propref{prop:fe2} is the key estimate that we will use to prove \thmref{thm:equ}. 
\begin{Cor}\label{cor:keyestimate}
Let $q$ be a positive integer. For any $\Psi\in C_c^{\infty}(\cM)$, $y>0$, and any $\epsilon>0$, we have
\begin{displaymath}
\sum_{m\neq 0}\left|a_{\Psi}(qm,y)\right|\ll_{\e}\cS_{2,2}(\Psi)q^{-1}y^{-(1/2+\theta+\e)}.
\end{displaymath}
\end{Cor}
\begin{proof}
If $qy\leq 1$ we can separate the above sum into two parts to get
\begin{align*}
\sum_{m\neq 0}\left|a_{\Psi}(qm,y)\right|&=\sum_{1\leq |m|\leq (qy)^{-1}}\left|a_{\Psi}(qm,y)\right|+\sum_{|m|>(qy)^{-1}}\left|a_{\Psi}(qm,y)\right|.
\end{align*}
Applying \eqref{equ:smallm} (and the estimate $\cS_{\alpha_0}(\Psi)\leq \cS_{2,2}(\Psi)$ by \rmkref{rmk:sobolev}) to the first sum and \propref{prop:fe2} to the second, we have
\begin{align*}
\sum_{m\neq 0}\left|a_{\Psi}(qm,y)\right|&\ll_{\e}\cS_{2,2}(\Psi)\left(\sum_{1\leq |m|\leq (ny)^{-1}}|qm|^{\theta}y^{1/2-\e}+\sum_{|m|>(qy)^{-1}}|qm|^{-4/3+\theta+\e}y^{-5/6}\right)\\
&\asymp \cS_{2,2}(\Psi)\left(q^{\theta}y^{1/2-\e}(qy)^{-(1+\theta)}+q^{-4/3+\theta+\e}y^{-5/6}(qy)^{1/3-\theta-\e}\right)\\
&=\cS_{2,2}(\Psi)q^{-1}y^{-(1/2+\theta+\e)},
\end{align*}
where for the second estimate we used that $4/3-\theta-\e>1$. If $qy> 1$ then we have $|qm|y>1$ for all $m\neq 0$. We can apply \propref{prop:fe2} to $a_{\Psi}(qm,y)$ for all integers $m\neq 0$ to get
\begin{align*}
\sum_{m\neq 0}\left|a_{\Psi}(qm,y)\right|&\ll_{\e} \cS_{2,2}(\Psi)\sum_{|m|\neq 0}|qm|^{-4/3+\theta+\e}y^{-5/6}\\
&\ll \cS_{2,2}(\Psi) q^{-4/3+\theta+\e} y^{-5/6}\ll \cS_{2,2}(\Psi)q^{-1}y^{-(1/2+\theta+\e)},
\end{align*}
where for the last estimate we used that $\theta<1/3-\e$.% $\sum_{k\neq 0}|k|^{\theta-4/3+3\e/2}\ll 1< (qy)^{-(\theta-1/3+3\e/2)}$ (since $\theta-4/3+3\e/2<-1$ and $qy>1$). This finishes the proof.
\end{proof}
\begin{rmk}\label{rmk:extension}
The estimates in \cite{KelmerKontorovich2020} hold more generally for any $\Gamma$ conjugate to some $\Gamma_0(p)$. In this generality, there might be (finitely many) exceptional cusp forms with $r_k\in i(0,\theta]$. For such forms, it was shown in \cite[Lemma 3.7]{KelmerKontorovich2020} that for any $m\neq 0$
$$\left|a_{\phi_k}(m,y)\right|\ll_{\e,p}\|\Psi\|_2|m|^{\theta}y^{1/2-\e}(|m|y)^{-|r_k|+\e}e^{-2\pi |m|y}.$$
Using the estimates $(|m|y)^{-|r_k|+\e}e^{-2\pi |m|y}< (|m|y)^{-\theta}$ when $|m|y\leq 1$ and $(|m|y)^{-|r_k|+\e}e^{-2\pi |m|y}\ll (|m|y)^{-2}$ when $|m|y>1$ one can easily recover \corref{cor:keyestimate} for $\phi_k$, and hence for a general $\Psi\in C_c^{\infty}(\G_0(p)\bk \bH)$.  Then one can easily deduce analogous estimates as in \thmref{thm:equ} for $\Psi$, see the arguments in the next subsection.
%Let $\Psi\in C_c^{\infty}(\G_0(p)\bk \bH)$ with $p$ a prime number. We note that if one can prove the estimate in \corref{cor:keyestimate} for $\Psi$, then one can easily deduce analogous estimates as in \eqref{equ:effbound} and \eqref{equ:prmeffbound} for $\Psi$, see the arguments in the next subsection. If $\Psi$ is orthogonal to all the (finitely many) exceptional cusp forms $\phi_k\in C_c^{\infty}(\G_0(p)\bk \bH)$, then \corref{cor:keyestimate} holds for $\Psi$ using the exact same arguments as above. It thus remains to deal with the exceptional cusp forms $\phi_k$. For any such $\phi_k$, it was shown in \cite[Lemma 3.7]{KelmerKontorovich2020} that for any $m\neq 0$
%$$\left|a_{\phi_k}(m,y)\right|\ll_{\e,p}\|\Psi\|_2|m|^{\theta}y^{1/2-\e}(|m|y)^{-|r_k|+\e}e^{-2\pi |m|y}$$
%with $r_k\in i(0,\theta]$ the parameter such that $\Delta\phi_k=(1/4+r_k^2)\phi_k$. We note that using the estimates $(|m|y)^{-|r_k|+\e}e^{-2\pi |m|y}< (|m|y)^{-\theta}$ when $|m|y\leq 1$ and $(|m|y)^{-|r_k|+\e}e^{-2\pi |m|y}\ll (|m|y)^{-2}$ when $|m|y>1$ one can easily recover \corref{cor:keyestimate} for $\phi_k$, and hence for a general $\Psi\in C_c^{\infty}(\G_0(p)\bk \bH)$. 
\end{rmk}
\subsection{Proof of \thmref{thm:equ}}
In this subsection we prove \thmref{thm:equ}. In view of \eqref{equ:zerocoef} it suffices to prove the following proposition.
\begin{Prop}
Let $\cM$ be the modular surface. For any $\Psi\in C_c^\infty(\cM)$, for any $x\in\R/\Z$ and $y>0$, we have
\begin{equation}\label{equ:effbound}
\delta_{n,x,y}(\Psi)=a_{\Psi}(0,y)+O_{\e}\left(\cS_{2,2}(\Psi)n^{-1}y^{-(1/2+\theta+\e)}\right)
\end{equation}
and
\begin{equation}\label{equ:prmeffbound}
\delta^{\rm pr}_{n,x,y}(\Psi)=a_{\Psi}(0,y)+O_{\e}\left(\cS_{2,2}(\Psi)n^{-1+\e}y^{-(1/2+\theta+\e)}\right).
\end{equation}
\end{Prop}
\begin{proof}
Let $J\subset \R/\Z\cong [0,1)$ be a finite subset and for any $m\in\Z$ denote by $W_J(m):=\frac{1}{|J|}\sum_{t\in J}e(mt)$. We note that $\frac{1}{|J|}\sum_{t\in J}\Psi(t+iy)$ equals $\delta_{n,x,y}(\Psi)$ when $J=\{x+j/n :0\leq j\leq n-1\}$ and equals $\delta_{n,x,y}^{\rm pr}(\Psi)$ when $J=\left\{x+j/n:0\leq j\leq n-1, \gcd(j,n)=1\right\}$. Applying the Fourier expansion \eqref{equ:foex} to $\Psi$ we get that
\begin{align*}
\frac{1}{|J|}\sum_{t\in J}\Psi(t+iy)&=\frac{1}{|J|}\sum_{t\in J}\sum_{m\in\Z}a_{\Psi}(m,y)e(mt)=\sum_{m\in\Z}a_{\Psi}(m,y)\frac{1}{|J|}\sum_{t\in J}e(mt)\\
&=a_{\Psi}(0,y)+\sum_{m\neq 0}a_{\Psi}(m,y)W_J(m).
\end{align*}
Now for \eqref{equ:effbound} we take $J=\{x+j/n:0\leq j\leq n-1\}$ and note that for such $J$, $|W_{J}(m)|$ equals $1$ if $n\mid m$ and equals $0$ otherwise. Hence 
$$\left|\sum_{m\neq 0}a_{\Psi}(m,y)W_J(m)\right|\ \leq\ \sum_{\substack{m\neq 0\\ n | m}}|a_{\Psi}(m,y)|\ \ll_{\e}\ n^{-1}y^{-(1/2+\theta+\e)},$$
where for the last estimate we applied \corref{cor:keyestimate}. 

For \eqref{equ:prmeffbound} we take $J=\left\{x+j/n:0\leq j\leq n-1, \gcd(j,n)=1\right\}$ and note the identity
$$\sum_{j\in (\Z/n\Z)^{\times}}e\left(\tfrac{mj}{n}\right)=\frac{\mu(n_m)\varphi(n)}{\varphi(n_m)}$$ 
for the Ramanujan's sum, where $n_m:=n/\gcd(n,m)$ and $\mu: \N\to \{0, \pm 1\}$ is the M\"{o}bius function; see e.g. \cite[Theorem 272]{HardyWright2008}. Then
$$|W_J(m)|\ =\ \left|\frac{1}{\varphi(n)}\sum_{j\in (\Z/n\Z)^{\times}}e\left(\tfrac{mj}{n}\right)\right|\ =\ \frac{|\mu(n_m)|}{\varphi\left(n_m\right)}\ \leq\ \frac{1}{\varphi(n_m)}.$$
Hence we have
\begin{align*}
\left|\sum_{m\neq 0}a_{\Psi}(m,y)W_J(m)\right|\ &\leq\ \sum_{m\neq 0}\frac{\left|a_{\Psi}(m,y)\right|}{\varphi(n_m)}=\sum_{d| n}\frac{1}{\varphi(d)}\sum_{\substack{m\neq 0\\ \gcd(m,n)=n/d}}\left|a_{\Psi}(m,y)\right|\\
&\leq\ \sum_{d| n}\frac{1}{\varphi(d)} \sum_{\substack{m\neq 0\\ (n/d) | m}}\left|a_{\Psi}(m,y)\right|\ll_{\e}\sum_{d\mid n}\frac{1}{\varphi(d)}\left(\frac{n}{d}\right)^{-1}y^{-(1/2+\theta+\e)}\\
&\ll_{\e}\ n^{-1}\sigma_{\e/2}(n)y^{-(1/2+\theta+\e)}\ll_{\e}n^{-1+\e}y^{-(1/2+\theta+\e)},
\end{align*}
where for the second inequality we used the fact that $\gcd(m,n)=n/d$ implies that $(n/d)\mid m$, for the third inequality we applied \corref{cor:keyestimate} and for the second last inequality we applied the estimate $\varphi(d)\gg_{\e} d^{1-\e/2}$. 
\end{proof}

\subsection{Full range equidistribution for rational translates}\label{sec:fullrange}
%Let $y_n\asymp 1/n^{\alpha}$ for some $\alpha>0$. As mentioned in the introduction, \thmref{thm:equ} confirms equidistribution of the two sequences $\{\cR_n(x,y_n)\}_{n\in\N}$ and $\{\cR_n^{\rm pr}(x,y_n)\}_{n\in\N}$ for any $0<\alpha< \tfrac{2}{1+2\theta}=\tfrac{64}{39}$ (for any $0<\alpha<2$ assuming the Ramanujan conjecture). This upper bound $\tfrac{64}{39}$ is a natural barrier for our argument using spectral estimates. In this subsection we show that one can go beyond this barrier to prove the full range equidistribution for rational translates. 
In this subsection we prove \thmref{thm:fullrange}. We fix $x=p/q$ a primitive rational number and let 
$$\N_q=\left\{n\in\N: \gcd(n^2, q)\mid n\right\}$$
be as in \thmref{thm:fullrange}.
As mentioned in the introduction, the key ingredient is a symmetry lemma for rational translates which generalizes the symmetry \eqref{equ:symmetryzeor}. Before stating the lemma, let us briefly explain why we need to restrict to the subsequence $\N_q$. Let $n\in\N$ and let $y>0$. We need to study the distribution of the points $\G(x+\tfrac{j}{n}+iy)=\G(\tfrac{p}{q}+\tfrac{j}{n}+iy)$ for $0\leq j\leq n-1$.
Let $\tfrac{p_j}{q_j}$ be the reduced form of $\tfrac{p}{q}+\tfrac{j}{n}$ and in view of the symmetry \eqref{equ:symmetryzeor} we have
$$\G\left(x+\tfrac{j}{n}+iy\right)=\G\left(\tfrac{p_j}{q_j}+iy\right)=\G\left(-\tfrac{\overline{p_j}}{q_j}+\tfrac{i}{q_j^2y}\right),$$
where $\overline{p_j}$ is the multiplicative inverse of $p_j$ modulo $q_j$. To further analyze the distribution of these points, we thus need to solve the congruence equation $xp_j\equiv 1\Mod{q_j}$ in $x$. Write $k=\gcd(n, q)$ and $q'=q/k$ and $n'=n/k$. Then
$$\tfrac{p}{q}+\tfrac{j}{n}=\tfrac{p}{kq'}+\tfrac{j}{kn'}=\tfrac{pn'+jq'}{kq'n'},$$
implying that 
$$q_j=\tfrac{kq'n'}{\gcd(pn'+jq', kq'n')}=\tfrac{kn'q'}{\gcd(pn'+jq', kn')}=q'\tfrac{n}{\gcd(pn'+jq',n)}$$
can be written canonically as a product of two integers. Here for the second equality we used that $\gcd(pn'+jq',q')=\gcd(pn',q')=1$. In view of the Chinese remainder theorem, the above congruence equation modulo $q_j$ is relatively easy to solve when the two factors $q'$ and $n/\gcd(pn'+jq', n)$ are coprime (see the proof of \lemref{lem:symmetryrational} for more details). This condition can be guaranteed for any $j$ if $\gcd(q', n)=\gcd(q/\gcd(q,n),n)=1$ which is equivalent to the condition $n\in\N_q$. Finally, we also note that by writing $n$ and $q$ in prime decomposition forms, it is not hard to check that $n\in\N_q$ is equivalent to $q=kl$ with $l=\gcd(n,q)\mid n$ and $\gcd(k,n)=1$.
We now state the symmetry lemma.
\begin{Lem}\label{lem:symmetryrational}
Let $\tfrac{m}{kl}$ be a primitive rational number and let $n\in\N$ such that $l\mid n$ and $\gcd(k,n)=1$. Then for any $0\leq j\leq n-1$ and for any $y>0$ we have
\begin{equation}\label{equ:symmetrywecond}
\G(\tfrac{m}{kl}+\tfrac{j}{n}+iy)=\G\left(-\tfrac{dl\overline{mn}a}{k}-\tfrac{\left(\left(m\tfrac{n}{l}+jk\right)/d\right)^*b}{n/d}+i\tfrac{d^2}{k^2n^2y}\right),
\end{equation} 
where $d=d_{j}:=\gcd(m\tfrac{n}{l}+jk, n)$ and $a=a_d$, $b=b_d\in\Z$ are some fixed integers such that $a\tfrac{n}{d}+bk=1$. Here, for any integer $x$, $\overline{x}$ denotes the multiplicative inverse of $x$ modulo $k$, $x^*$ denotes the multiplicative inverse of $x$ modulo $n/d$. If we further assume $\gcd(j,n)=l=1$, then $d_j=\gcd(mn+jk, n)=1$ and
\begin{equation}\label{equ:prsymm}
\G(\tfrac{m}{k}+\tfrac{j}{n}+iy)=\G\left(-\tfrac{\overline{mn}a}{k}-\tfrac{(jk)^*b}{n}+\tfrac{i}{k^2n^2y}\right).
\end{equation}
\end{Lem}
\begin{proof}
Since $l\mid n$, by direct computation we have $\frac{m}{kl}+\frac{j}{n}=\frac{mn/l+jk}{kn}$. Note that since $\gcd(k, mn)=1$ we have $\gcd(m\tfrac{n}{l}+jk, k)=\gcd(m\tfrac{n}{l}, k)=1$. This implies that $\gcd(m\tfrac{n}{l}+jk, kn)=\gcd(m\tfrac{n}{l}+jk, n)=d$. Hence let $\tfrac{p}{q}$ be the reduced form of $\tfrac{m}{kl}+\tfrac{j}{n}$, then we have $(p,q)=((m\tfrac{n}{l}+jk)/d, kn/d)$. Now since $\gcd(p,q)=1$, there exist some integers $v, w\in\Z$ such that $\gamma=\left(\begin{smallmatrix}
w & v\\
-q & p\end{smallmatrix}\right)\in \G$. By direct computation we have
$$\gamma\left(\tfrac{m}{kl}+\tfrac{j}{n}+iy\right)=\gamma\left(\tfrac{p}{q}+iy\right)=-\tfrac{w}{q}+\tfrac{i}{q^2y}.$$
implying that
\begin{equation}\label{equ:prerelation}
\G\left(\tfrac{m}{kl}+\tfrac{j}{n}+iy\right)=\Gamma\left(-\tfrac{w}{q}+\tfrac{i}{q^2y}\right)=\Gamma\left(-\tfrac{w}{kn/d}+i\tfrac{d^2}{k^2n^2y}\right),
\end{equation}
where for the second equality we used the relation $q=kn/d$. Moreover, since $\gamma\in \G$ we have $wp+vq=1$, implying that (again using the relation $(p,q)=((m\tfrac{n}{l}+jk)/d, kn/d)$)
$$w\left((m\tfrac{n}{l}+jk)/d\right)\equiv 1\Mod{k\tfrac{n}{d}}.$$
We claim that 
\begin{equation}\label{equ:solution}
w\equiv dl\overline{mn}\tfrac{n}{d}a+\left(\left(m\tfrac{n}{l}+jk\right)/d\right)^*kb \Mod{k\tfrac{n}{d}}.
\end{equation}
%$$w\equiv dl\overline{pn}\tfrac{n}{d}a+\left(\left(p\tfrac{n}{l}+jq\right)/d\right)^*qb \Mod{q\tfrac{n}{d}}.$$
In view of the Chinese Remainder Theorem, since $\gcd(k, n/d)=1$, it suffices to check
$$\left(dl\overline{mn}\tfrac{n}{d}a+\left(\left(m\tfrac{n}{l}+jk\right)/d\right)^*kb\right)\left((m\tfrac{n}{l}+jk)/d\right)\equiv 1\Mod{k}$$
and
$$\left(dl\overline{mn}\tfrac{n}{d}a+\left(\left(m\tfrac{n}{l}+jk\right)/d\right)^*kb\right)\left((m\tfrac{n}{l}+jk)/d\right)\equiv 1\Mod{\tfrac{n}{d}}.$$
For the first equation we have
\begin{align*}
\left(dl\overline{mn}\tfrac{n}{d}a+\left(\left(m\tfrac{n}{l}+jk\right)/d\right)^*kb\right)\left((m\tfrac{n}{l}+jk)/d\right)\equiv dl\overline{mn}\tfrac{n}{d}a mn\overline{ld}\equiv a\tfrac{n}{d}=1-bk\equiv 1\Mod{k},
\end{align*}
where for the first equality we used the fact that $\gcd(dl, k)=1$ (since $d\mid n$, $l\mid n$ and $\gcd(k,n)=1$). The second equation follows similarly. Now plugging relation \eqref{equ:solution} into \eqref{equ:prerelation} we get \eqref{equ:symmetrywecond}.

For the second half we note that $d_j=\gcd(mn+jk, n)=\gcd(jk,n)=1$. The first equality is true since $l=1$, and the second equality is true since by assumption $\gcd(k,n)=\gcd(j,n)=1$. Thus in view of \eqref{equ:symmetrywecond}, to prove \eqref{equ:prsymm} it suffices to note that $(mn+jk)^*\equiv (jk)^*\Mod{n}$, or equivalently, $mn+jk\equiv jk\Mod{n}$.
%Now since $\gcd(pn+qj, qn)=1$, there exist $c,d\in\Z$ such that $\gamma=\left(\begin{smallmatrix}
%c & d\\
%-qn & pn+qj\end{smallmatrix}\right)\in \G$. By direct computation we have
%$$\gamma(\tfrac{p}{q}+\tfrac{j}{n}+iy)=\gamma(\tfrac{pn+qj}{qn}+iy)=-\tfrac{c}{qn}+\tfrac{i}{q^2n^2y},$$
%implying that $\G(\tfrac{p}{q}+\tfrac{j}{n}+iy)=\G\left(-\tfrac{c}{qn}+\tfrac{i}{q^2n^2y}\right)$.
%Moreover, note that $\tfrac{\overline{pn}a}{q}+\tfrac{(qj)^*b}{n}=\tfrac{\overline{pn}an+(qj)^*bq}{qn}$, 
%%$$\G\left(-\tfrac{\overline{pn}a}{q}-\tfrac{(qj)^*b}{n}+\tfrac{i}{q^2n^2y}\right)=\G\left(-\tfrac{\overline{pn}an+(qj)^*bq}{qn}+\tfrac{i}{q^2n^2y}\right).$$ 
%it thus suffices to show that $c\equiv \overline{pn}an+(qj)^*bq\Mod{qn}$. Since $\gamma\in \G$ we have $c(pn+qj)=1-qnd\equiv 1\Mod{qn}$. We thus need to show
%$$\left(\overline{pn}an+(qj)^*bq\right)(pn+qj)\equiv 1\Mod{qn}.$$
%In view of the Chinese remainder theorem, since $\gcd(q,n)=1$ it suffices to show 
%$$\left(\overline{pn}an+(qj)^*bq\right)(pn+qj)\equiv 1\Mod{q}\ \textrm{and}\ \left(\overline{pn}an+(qj)^*bq\right)(pn+qj)\equiv 1\Mod{n}.$$
%For the first equality we have 
%$$\left(\overline{pn}an+(qj)^*bq\right)(pn+qj)\equiv \overline{pn}an pn\equiv an=1-bq\equiv 1\Mod{q}.$$ 
%Similarly, for the second equality we have
%$$\left(\overline{pn}an+(qj)^*bq\right)(pn+qj)\equiv (qj)^*bq qj\equiv bq=1-an\equiv 1\Mod{n}.$$
%This finishes the proof.
\end{proof}
\begin{rmk}
When $k=1$ we can take $(a,b)=(0,1)$, then \eqref{equ:prsymm} recovers the symmetry \eqref{equ:symmetryzeor}. We also note that for the point $\G(x+j/n+iy)$ with $x$ irrational, the above symmetry clearly breaks.
\end{rmk}

%In the rest of this section we consider the effect of rational translates $x=p/q$ taken in reduced form, i.e., $q\in\N$, $\gcd(p,q)=1$. To apply the symmetry relation in \lemref{lem:symmetryrational}, we will start from the trivial factorization $q= \gcd(q,n)\cdot\tfrac{q}{\gcd(q,n)}$. Clearly $\gcd(q,n)\mid n$, and we easily check that $\gcd(\tfrac{q}{\gcd(q,n)},n)=1$ is guaranteed for all $n\in\N_q$, given by
%$$
%\N_q = \{ n\in\N: (n^2,q) \mid n\}.
%$$
%Indeed, 
%$$
%\gcd(\tfrac{q}{\gcd(q,n)},n) = \gcd(\tfrac{q}{\gcd(q,n)},\tfrac{n}{\gcd(q,n)}\cdot\gcd(q,n)) = \gcd(\tfrac{q}{\gcd(q,n)},\gcd(q,n))= 1
%$$
%if and only if for each prime $\ell\mid n$, $\ell^s \parallel n$ implies that $\ell^{s+1} \nmid q$, or equivalently, $\gcd(n^2,q)\mid n$. Here, $\ell^s\parallel n$ means that $\ell^s$ is the highest power of $\ell$ dividing $n$.
%As a direct corollary of this symmetry we have the following relation:
\begin{Prop}\label{prop:symmetryset}
Let $p/q$ be a primitive rational number and let $n\in\N_q$. Then for any $y>0$ we have
\begin{equation}\label{equ:symreg}
\cR_n\left(\tfrac{p}{q},y\right)=\bigcup_{d| n}\cR_{n/d}^{\rm pr}\left(x_{d}, \tfrac{d^2}{k^2n^2y}\right),
\end{equation}
where $x_{d}\in\R/\Z$ is some number depending on $d$ $($and also on $p, q, n$$)$ and $k:=q/\gcd(n, q)$. If we further assume $\gcd(n,q)=1$, then
\begin{equation}\label{equ:sympri}
\cR_n^{\rm pr}\left(\tfrac{p}{q}, y\right)=\cR_n^{\rm pr}\left(-\tfrac{\overline{pn}a}{q}, \tfrac{1}{q^2n^2y}\right),
\end{equation}
where $\overline{x}$ denotes the multiplicative inverse of $x$ modulo $q$ and $a\in \Z$ is as in \lemref{lem:symmetryrational}.
\end{Prop}
\begin{proof}
Relation \eqref{equ:sympri} follows immediately from \eqref{equ:prsymm} by taking $(m,k)=(p,q)$ and noting that 
$$\{(-[(qj)^*b]\in (\Z/n\Z)^{\times}: j\in (\Z/n\Z)^{\times}\}=(\Z/n\Z)^{\times},$$
which follows from the fact that $\gcd(bq,n)=1$ (since $\gcd(bq, n)=\gcd(1-an, n)=1$). Here $(qj)^*$ denotes the multiplicative inverse of $qj$ modulo $n$ and $b\in\Z$ is as in \lemref{lem:symmetryrational}.

For \eqref{equ:symreg}, we set $m=p$, $l=\gcd(n, q)$ (so that $k=q/l$). As mentioned above, the condition $\gcd(n^2, q)\mid n$ implies that $\gcd(k,n)=1$. Thus the pair $(\tfrac{m}{kl}, n)$ satisfies the assumptions in \lemref{lem:symmetryrational} and we can apply \eqref{equ:symmetrywecond} for the points 
$$\G\left(\tfrac{p}{q}+\tfrac{j}{n}+iy\right)=\G\left(\tfrac{m}{kl}+\tfrac{j}{n}+iy\right), 0\leq j\leq n-1.$$
Now for any $d\mid n$ define
$$D_d:=\left\{0\leq j\leq n-1: d_{j}=\gcd(m\tfrac{n}{l}+jk, n)=d\right\}$$
so that
\begin{equation}\label{equ:predecom}
\cR_n\left(\tfrac{p}{q},y\right)=\bigcup_{d\mid n}\left\{\G\left(\tfrac{p}{q}+\tfrac{j}{n}+iy\right)\in\cM: j\in D_d\right\}.
\end{equation}
Moreover, we note that since $\gcd(k,n)=1$, we have $\left\{[m\tfrac{n}{l}+jk]\in \Z/n\Z: 0\leq j\leq n-1\right\}=\Z/n\Z$ and hence
\begin{equation}\label{equ:decomcongru}
\left\{[m\tfrac{n}{l}+jk]\in \Z/n\Z: j\in D_d\right\}=\left\{[j]\in \Z/n\Z: \gcd(j, n)=d\right\}.
\end{equation}
On the other hand, by \eqref{equ:symmetrywecond} we have
$$\left\{\G\left(\tfrac{p}{q}+\tfrac{j}{n}+iy\right)\in\cM: j\in D_d\right\}=\left\{\G\left(-\tfrac{dl\overline{mn}a_d}{k}-\tfrac{\left(\left(m\tfrac{n}{l}+jk\right)/d\right)^*b_d}{n/d}+i\tfrac{d^2}{k^2n^2y}\right)\in\cM: j\in D_d \right\},$$
where for any integer $x$, $\overline{x}$ denotes the multiplicative inverse of $x$ modulo $k$, $x^*$ denotes the multiplicative inverse of $x$ modulo $n/d$, and $a_d, b_d\in\Z$ are some fixed integers such that $a_d\tfrac{n}{d}+b_dk=1$. Now for each $d\mid n$ we let $x_d\in[0,1)$, $x_d\equiv-\tfrac{dl\overline{mn}a_d}{k} \Mod{1}$ so that it remains to show
$$\left\{-[\left((m\tfrac{n}{l}+jk)/d\right)^*b_d]\in (\Z/(n/d)\Z)^{\times}: j\in D_d\right\}=(\Z/(n/d)\Z)^{\times}.$$
We can thus conclude the proof by noting that the above relation follows immediately from \eqref{equ:decomcongru} together with the fact $\gcd(b_d,\tfrac{n}{d})=1$ (since $\gcd(b_d,\tfrac{n}{d})=\gcd(b_dk, \tfrac{n}{d})=\gcd(1-a_d\tfrac{n}{d}, \tfrac{n}{d})=1)$.
%Relation \eqref{equ:symreg} follows from \eqref{equ:sympri} and the decomposition
%\begin{displaymath}
%\cR_n(p/q, y)=\bigcup_{d\mid n}\cR_{n/d}^{\rm pr}(p/q, y).\qedhere
%\end{displaymath}
\end{proof}
Using these two relations and the estimate \eqref{equ:prmeffbound} one gets the following effective estimates.
\begin{Prop}\label{prop:symmetry}
Let $x=p/q$ be a primitive rational number and let $n\in\N_q$. Then for any $\Psi\in C_c^{\infty}(\cM)$ and $y>0$ we have
$$\delta_{n,x,y}(\Psi)=\frac{1}{n}\sum_{d\mid n}\varphi\left(\tfrac{n}{d}\right)a_{\Psi}\left(0, \tfrac{d^2}{k^2n^2y}\right)+O_{\e,q}\left(\cS_{2,2}(\Psi)n^{2\theta+4\e}y^{1/2+\theta+\e}\right),$$
where {{$k:=q/\gcd(n, q)$}}. If we further assume that $\gcd(n,q)=1$, then
$$\delta_{n,x,y}^{\rm pr}(\Psi)=a_{\Psi}\left(0, \tfrac{1}{q^2n^2y}\right)+O_{\e,q}\left(\cS_{2,2}(\Psi)n^{2\theta+3\e}y^{1/2+\theta+\e}\right).$$
\end{Prop}
\begin{proof}
For any positive divisor $d\mid n$, let $y_d=d^2/(k^2n^2y)$ with  {{$k:=q/\gcd(n, q)$}} as above and let $x_d\in\R/\Z$ be as in \eqref{equ:symreg}. Then by \eqref{equ:symreg} for $x=p/q$ we have
\begin{align*}
\delta_{n,x,y}(\Psi)&=\frac{1}{n}\sum_{d| n}\varphi\left(\tfrac{n}{d}\right)\delta^{\rm pr}_{n/d, x_d, y_d}(\Psi)\\
&=\frac{1}{n}\sum_{d\mid n}\varphi\left(\tfrac{n}{d}\right)\left(a_{\Psi}\left(0, y_d\right)+O_{\e}\left(\cS_{2,2}(\Psi)\left(\tfrac{n}{d}\right)^{-1+\e}y_d^{-(1/2+\theta+\e)}\right)\right)\\
&=\frac{1}{n}\sum_{d\mid n}\varphi\left(\tfrac{n}{d}\right)a_{\Psi}\left(0, y_d\right)+O_{\e}\left(\cS_{2,2}(\Psi)n^{-1}\sum_{d\mid n}\left(\tfrac{n}{d}\right)^{\e}y_d^{-(1/2+\theta+\e)}\right),
%&=\frac{1}{n}\sum_{d\mid n}\varphi\left(\frac{n}{d}\right)a_{\Psi}\left(0, \frac{d^2}{q^2n^2y}\right)+O_{\e,q}\left(n^
\end{align*}
where for the second estimate we applied \eqref{equ:prmeffbound} and for the third estimate we used the trivial estimate $\varphi(n/d)<n/d$. Now plugging $y_d=d^2/(k^2n^2y)$ into the above equation we get
\begin{align*}
\delta_{n,x,y}(\Psi)&=\frac{1}{n}\sum_{d\mid n}\varphi\left(\tfrac{n}{d}\right)a_{\Psi}\left(0, \tfrac{d^2}{k^2n^2y}\right)+O_{\e,q}\left(\cS_{2,2}(\Psi)n^{-1}\sigma_{1+2\theta+3\e}(n)y^{1/2+\theta+\e}\right)\\
&=\frac{1}{n}\sum_{d\mid n}\varphi\left(\tfrac{n}{d}\right)a_{\Psi}\left(0, \tfrac{d^2}{k^2n^2y}\right)+O_{\e,q}\left(\cS_{2,2}(\Psi)n^{2\theta+4\e}y^{1/2+\theta+\e}\right),
\end{align*}
where the dependence on $k$ in the first estimate is absorbed into the dependence on $q$ (since  ${{k:=q/\gcd(n, q)}}\leq q$).
The second estimate follows from similar (but easier) analysis with the relation \eqref{equ:sympri} in place of \eqref{equ:symreg}. 
\end{proof}

We are now in the position to prove \thmref{thm:fullrange}. We will prove the following proposition from which \thmref{thm:fullrange} follows, see also \rmkref{rmk:exponentde}.
%In view of \thmref{thm:equ}, \thmref{thm:fullrange} follows immediately from the following effective estimate. 
\begin{Thm}\label{thm:symmtry}
Let $x=p/q$ be a primitive rational number and let $n\in\N_q$. Let $y_n=c/n^{\alpha}$ for some $1<\alpha<2$ and $c>0$. Then for any $\Psi\in C_c^{\infty}(\cM)$ we have
$$\left|\delta_{n,x,y_n}(\Psi)-\mu_{\cM}(\Psi)\right|\ll_{\e,q, c,\Psi}n^{\alpha/2-1+\e}+n^{2\theta+4\e-\alpha(1/2+\theta+\e)}.$$
If we further assume $\gcd(n,q)=1$, then we have
$$\left|\delta^{\rm pr}_{n,x,y_n}(\Psi)-\mu_{\cM}(\Psi)\right|\ll_{\e, q, c}\cS_{2,2}(\Psi)\left(n^{\alpha/2-1}+n^{2\theta+3\e-\alpha(1/2+\theta+\e)}\right).$$
\end{Thm}
\begin{rmk}\label{rmk:exponentde}
The dependence on $\Psi$ in the first estimate can also be made explicit. In fact, we can remove this dependence by adding a factor of $\cS_{2,2}(\Psi)+\|\Psi\|_{\infty}$ to the right hand side of this estimate. We also note that since we may take $\theta=7/64$, the right hand side of these two estimates decays to zero as $n\to\infty$ for any $1<\alpha<2$.
\end{rmk}
\begin{proof}[Proof of \thmref{thm:symmtry}]
In view of \propref{prop:symmetry} and the assumption $y_n=c/n^{\alpha}$, it suffices to show that
$$\frac{1}{n}\sum_{d\mid n}\varphi\left(\tfrac{n}{d}\right)a_{\Psi}\left(0, \tfrac{d^2}{k^2n^2y_n}\right)=\mu_{\cM}(\Psi)+O_{\e, c,\Psi}\left(n^{\alpha/2-1+\e}\right)$$
with  {{$k:=q/\gcd(n, q)$}}, and that (under the extra assumption $\gcd(n,q)=1$)
$$a_{\Psi}\left(0, \tfrac{1}{q^2n^2y_n}\right)=\mu_{\cM}(\Psi)+O_{c}\left(\cS_{2,2}(\Psi)n^{\alpha/2-1}\right).$$

The second estimate follows immediately from \eqref{equ:zerocoef} and the trivial estimate $|q|\geq 1$. For the first estimate we separate the sum into two parts to get
\begin{align*}
\frac{1}{n}\sum_{d\mid n}\varphi\left(\tfrac{n}{d}\right)a_{\Psi}\left(0, \tfrac{d^2}{k^2n^2y_n}\right)&=\frac{1}{n}\left(\sum_{\substack{d\mid n\\ d< n^{1-\alpha/2}}}+\sum_{\substack{d| n\\ d\geq n^{1-\alpha/2}}}\right)\varphi\left(\tfrac{n}{d}\right)a_{\Psi}\left(0, \tfrac{d^2}{k^2n^2y_n}\right).
\end{align*}
Applying \eqref{equ:zerocoef} (and the trivial estimate $|k|\geq 1$) for the first sum and applying the estimate
$$\left|a_{\Psi}\left(0, \tfrac{d^2}{k^2n^2y_n}\right)\right|=\left|\int_0^1\Psi\left(t+i\tfrac{d^2}{k^2n^2y_n}\right)dt\right|\leq \|\Psi\|_{\infty}$$
for the second sum we get $\frac{1}{n}\sum_{d\mid n}\varphi\left(\tfrac{n}{d}\right)a_{\Psi}\left(0, \tfrac{d^2}{k^2n^2y_n}\right)$ equals
\begin{align*}
& \frac{1}{n}\left(\sum_{\substack{d\mid n\\ d< n^{1-\alpha/2}}}\varphi\left(\tfrac{n}{d}\right)\left(\mu_{\cM}(\Psi)+O_{c,\Psi}\left(\left(\tfrac{n}{d}\right)^{-1}n^{\alpha/2}\right)\right)+O_{\Psi}\left(\sum_{\substack{d| n\\ d\geq n^{1-\alpha/2}}}\varphi\left(\tfrac{n}{d}\right)\right)\right)\\
&=\mu_{\cM}(\Psi)+\frac{1}{n}O_{c, \Psi}\left(n^{\alpha/2}\sum_{\substack{d\mid n\\ d< n^{1-\alpha/2}}}1+\sum_{\substack{d| n\\ d\geq n^{1-\alpha/2}}}\tfrac{n}{d}\right)\\
&=\mu_{\cM}(\Psi)+O_{c, \Psi}\left(n^{\alpha/2-1}\sigma_0(n)\right)=\mu_{\cM}(\Psi)+O_{\e, c, \Psi}\left(n^{\alpha/2-1+\e}\right),
\end{align*}
finishing the proof, where for the first estimate we used the identity that $\sum_{d\mid n}\varphi(n/d)=n$ and the estimate that $\varphi\left(n/d\right)<n/d$, and for the second estimate we used the estimates $\sum_{\substack{d\mid n\\ d< n^{1-\alpha/2}}}1\leq \sigma_{0}(n)$ and
\begin{displaymath}
\sum_{\substack{d| n\\ d\geq n^{1-\alpha/2}}}\frac{n}{d}=\sum_{\substack{d\mid n\\ d\leq n^{\alpha/2}}}d\leq n^{\alpha/2}\sum_{\substack{d\mid n\\ d\leq n^{\alpha/2}}}1\leq n^{\alpha/2}\sigma_0(n).\qedhere
\end{displaymath}
\end{proof}
%Combining \thmref{thm:equ} and \thmref{thm:symmtry} we have the following full range equidistribution result for rational translates.
%\begin{Cor}\label{cor:fullrange}
%Let $x=p/q$ be a rational number with $\gcd(p,q)=1$ and let $\{y_n\}_{n\in\N}$ be a sequence satisfying $y_n\asymp 1/n^{\alpha}$ for some $\alpha\in (0,2)$. Then both $\delta_{n,x,y_n}$ and $\delta_{n,x,y_n}^{\rm pr}$ weakly converge to $\mu_{\cM}$ as $n\in\N_{\alpha}$ goes to infinity, where $\N_{\alpha}=\N$ if $0<\alpha< \tfrac{64}{39}$ and $\N_{\alpha}=\{n\in\N:\gcd(n,q)=1\}$ if $\tfrac{64}{39}\leq \alpha<2$.
%\end{Cor}

\subsection{Quantitative non-equidistribution for rational translates}
As a direct consequence of the analysis in the previous subsection we also have the following quantitative non-equidistribution result for rational translates when $\{y_n\}_{n\in\N}$ is beyond the above range, generalizing the situation for $\{\cR_n^{\rm pr}(0,y_n)\}_{n\in\N}$.  As before, for any $Y>0$ we denote by $\mu_Y$ the probability uniform distribution measure supported on $\cH_Y$.
\begin{Thm}\label{thm:effnonequ}
Let $x=p/q$ be a primitive rational number and let $y_n=c/n^2$ for some constant $c>0$. Let $\Psi\in C_c^{\infty}(\Psi)$. Then for any $n\in\N_q$ we have
%\begin{enumerate}
%\item If $y_n=c/n^2$ for some constant $c>0$, then for any $n\in\N$ with $\gcd(n^2,q)\mid n$ we have
$$\delta_{n,x,y_n}(\Psi)=\frac{1}{n}\sum_{d\mid n}\varphi\left(\tfrac{n}{d}\right)\mu_{\tfrac{d^2}{ck_n^2}}(\Psi)+O_{\e,q, c}\left(\cS_{2,2}(\Psi)n^{-1+2\e}\right)$$
with $k_n=q/\gcd(n^2, q)$. If we further assume that $\gcd(n,q)=1$, then
$$\delta_{n,x,y_n}^{\rm pr}(\Psi)=\mu_{\tfrac{1}{cq^2}}(\Psi)+O_{\e,q,c}\left(\cS_{2,2}(\Psi)n^{-1+\e}\right).$$
\end{Thm}
\begin{proof}
These two effective estimates follow immediately from \propref{prop:symmetry} by plugging in $y_n=c/n^2$ and noting that $a_{\Psi}(0, Y)=\int_0^1\Psi(x+iY)dx=\mu_Y(\Psi)$. %For the first estimate in part (2) we apply the first estimate in \propref{prop:symmetry} for the sequence $\{y_n\}_{n\in\N}$ with $\gcd(n^2,q)\mid n$ and noting that for such $n\gg 1$ sufficiently large and for any positive divisor $d\mid n$, $\mu_{d^2/(k^2n^2y_n)}(\Psi)=0$ (since 
%%$$d^2/(k^2n^2y_n)=d^2\gcd(n^2,q)^2/(q^2n^2y_n)\geq 1/(q^2n^2y_n)\to \infty$$ 
%%as $n\to\infty$ by the assumption $\lim\limits_{n\to\infty}n^2y_n=0$ and $\Psi$ is compactly supported). The second estimate in part (2) follows similarly by applying the second estimate in \propref{prop:symmetry}.
\end{proof}
We can now give the
\begin{proof}[Proof of \thmref{thm:rationalnonequ}]
For part (1), in view of \thmref{thm:effnonequ} only the second equation needs a proof. Since we are taking $n\in\bP_m$ going to infinity, it is sufficient to consider $n=m\ell\in \bP_m$ with the prime number $\ell>q$ (so that $\ell\nmid q$). For such $n$, we have $\gcd(n^2, q)=\gcd(m^2\ell^2, q)=\gcd(m^2, q)$. Since by assumption $\gcd(m^2, q)\mid m$ and $m\mid n$, we can apply the first effective estimate in \thmref{thm:effnonequ} for such $n=m\ell\in \bP_m$. Moreover, for any such $n$ we have 
$$k_n=\frac{q}{\gcd(n^2, q)}=\frac{q}{\gcd(m^2, q)}=\frac{q}{\gcd(m,q)}$$
is a fixed number only depending on $m$ and $q$. Here for the last equality we used the assumption that $\gcd(m^2, q)\mid m$. Now let $n=m\ell\in \bP_{m}$ with $\ell\gg q$ sufficiently large such that $\mu_{Y}(\Psi)=0$ whenever $Y>\ell^2/(ck_n)^2$ (this can be guaranteed since $k_n$ is a fixed number and $\Psi$ is compactly supported). In particular, for any $d\mid n$, $\mu_{d^2/(ck_n^2)}(\Psi)=0$ whenever $\ell\mid d$. This, together with the first estimate in \thmref{thm:effnonequ} implies that for all such sufficiently large $n=m\ell\in \bP_{m}$
\begin{align*}
\delta_{n,x,y_n}(\Psi)&=\frac{1}{m\ell}\sum_{d\mid m}\varphi\left(\tfrac{m\ell}{d}\right)\mu_{\tfrac{d^2}{ck_n^2}}(\Psi)+O_{\e,q,c,\Psi,m}\left(\ell^{-1+2\e}\right)\\
&=\frac{\ell-1}{\ell}\nu_{m,\tfrac{1}{ck_n^2}}(\Psi)+O_{\e,q,c,\Psi,m}\left(\ell^{-1+2\e}\right),
\end{align*}
where for the second estimate we used that $\gcd(m,\ell)=1$ and $\ell$ is a prime number. We can now finish the proof by taking $n=m\ell\to\infty$ along the subsequence $\bP_{m}$ (equivalently, taking $\ell\to\infty$) and plugging in the relation $k_n=q/\gcd(m,q)$. 

For part (2), since $\cR_n^{\rm pr}(x,y_n)\subset \cR_n(x,y_n)$, we only need to prove the full escape to the cusp for the sequence $\{\cR_n(x,y_n)\}_{n\in\N}$. Identify (up to a null set) $\cM$ with the standard fundamental domain $\cF_{\G}:=\left\{z\in \bH: \Re(z)<\frac12, |z|>1\right\}$. For any $n\in\N$ and $0\leq j\leq n-1$ let $\tfrac{p_j}{q_j}$ be the reduced form of $x+\tfrac{j}{n}=\tfrac{p}{q}+\tfrac{j}{n}=\tfrac{pn+qj}{qn}$ so that by \eqref{equ:symmetryzeor}
$$\G\left(x+\tfrac{j}{n}+iy_n\right)=\G\left(-\tfrac{\overline{p_j}}{q_j}+\tfrac{i}{q_j^2y_n}\right).$$ 
Thus using the trivial inequality $|q_j|\leq |q|n$ for all $0\leq j\leq n-1$ and the assumption $\lim\limits_{n\to\infty}n^2y_n=0$, we have
\begin{displaymath}
\cR(x,y_n)\subset \left\{z\in \cF_{\G}: \Im(z)\geq \tfrac{1}{q^2n^2y_n}\right\}\xrightarrow{n\to\infty} \textrm{cusp of $\cM$}.\qedhere
\end{displaymath}
\end{proof}

\section{Negative results: in connection with Diophantine approximations}\label{sec:cou1}

Let $\G=\SL_2(\Z)$ and $\cM=\G\bk \bH$ be the modular surface. Let $\mu_{\cM}$ be the normalized hyperbolic area on $\cM$ as before. In this section we prove a general result which captures the cusp excursion rate for the sample points $\cR_n(x,y_n)$ in terms of the Diophantine property of the translate $x\in \R/\Z \cong [0,1)$, see \thmref{thm:nonequ2}. \thmref{thm:nonequ2intro} will then be an easy consequence of this result. %We also prove a non-equidistribution result which holds for any translate, see \thmref{thm:nonequ1}. 

\subsection{Notation and a preliminary result on cusp excursions}
In this subsection we prove a preliminary lemma relating cusp excursions on the modular surface to Diophantine approximations. Let us first fix some notation. For any $Y>0$, we denote by $\cC_Y\subset \cM$ the image of the region
$$\{z\in\bH:\Im(z)>Y\}$$
under the natural projection from $\bH$ to $\cM=\G\bk\bH$. 
%\begin{equation}\label{equ:cuspnbhd}
%\cC_Y:=\{\G z\in\cM:\Im(z)>Y\}\subset \cM
%\end{equation}
%be the image of the open region $\left\{z\in\bH:\Im(z)>Y\right\}$ onto the modular surface $\cM$ under the natural projection from $\bH$ to $\cM$. 
As $Y$ goes to infinity, the sets $\cC_Y$ diverge to the cusp of $\cM$, and we call $\cC_Y$ a \textit{cusp neighborhood of $\cM$}. Similarly, for any $Y'>Y>0$, we denote by $\cC_{Y,Y'}$ the projection onto $\cM$ of the open set
%\begin{equation}\label{equ:cuspnbhddiff}
%\cC_{Y,Y'}:=
$$\left\{z\in\bH:Y<\Im(z)< Y'\right\}.$$%\subset\cM
%\end{equation} 
%the projection of the open strip $\left\{z\in \bH:Y<\Im(z)<Y'\right\}$ onto $\cM$. 

For any primitive rational number $m/n$, and for any $r>0$ we denote by 
$$H_{m/n,r}:=\left\{z=x+iy\in\bH:(x-m/n)^2+(y-r)^2= r^2\right\}$$ 
the horocycle tangent to $\partial \bH$ at $m/n$ with Euclidean radius $r$. We denote by
$$H^{\circ}_{m/n,r}:=\left\{z=x+iy\in\bH:(x-m/n)^2+(y-r)^2< r^2\right\}$$
the open horodisc enclosed by $H_{m/n,r}$. We have the following geometric description of \lemref{lem:symmetryrational}: Let $\gamma=\left(\begin{smallmatrix}
m & * \\
n & *\end{smallmatrix}\right)$ be an element in $\G$. Then $\gamma$ sends the horizontal horocycle $\{z\in\bH:\Im(z)=Y\}$ 
to the horocycle $H_{m/n,r}$ with $r=1/(2Yn^2)$, while the open region $\left\{z\in\bH:\Im(z)>Y\right\}$ is mapped to the horodisc $H^{\circ}_{m/n,r}$. %Moreover, we also note that when $Y\geq 1$ all these horodiscs are disjoint from each other. 
On the other hand, for any primitive rational number $m/n$, there is $\gamma\in \G$ of the form $\gamma=\bsm m&*\\ n&*\esm$. Thus for any $Y>0$ and for any $z\in \bH$, $\G z\in \cC_Y$ if and only if $z\in H^{\circ}_{m/n,r}$ for some primitive rational number $m/n$ with $r=1/(2Yn^2)$.

Finally, we record a distance formula that we will later use. Let $d_{\cM}(\cdot,\cdot)$ be the distance function on $\cM$ induced from the hyperbolic distance function $d_{\bH}$ on $\bH$, i.e.,
$$d_{\cM}(\G z_1, \G z_2)=\inf_{\gamma\in\G}d_{\bH}(\gamma z_1, z_2).$$
%Let $d_{\cM}(z_0, z_1):=\inf_{\gamma\in\G_1}d_{\bH}(\gamma z_0, z_1)$ be the distance function on $\cM$ as before. Then f
\begin{Lem}\label{lem:distance}
Let $\G z_0\in\cM$ be a fixed base point. Then there exists a constant $c>0$ $($which may depend on $\G z_0$$)$ such that for any $Y>1$ and for any $\G z\in \cC_Y$
\begin{equation}\label{equ:distance}
d_{\cM}(\G z_0,\G z)\geq \log Y - c.
\end{equation}
\end{Lem}

The estimate \eqref{lem:distance} holds for a general non-compact finite-volume hyperbolic manifold using  reduction theory after Garland and Raghunathan \cite[Theorem 0.6]{GarlandRaghunathan1970} combined with a distance estimate by Borel \cite[Theorem C]{Borel1972}. We give here a self-contained elementary proof for the special case of the modular surface.  
\begin{proof}[Proof of \lemref{lem:distance}]
In view of the triangle inequality, we may assume $\G z_0=\G i$. Note that $d_{\bH}(i, z)\geq \log Y$ for any $z\in \bH$ with $\Im(z)\in (0, 1/Y)\cup (Y,\infty)$. Thus it suffices to show that if $\G z\in \cC_Y$, then $\Im(\gamma z)\in (0,1/Y)\cup (Y,\infty)$ for any $\gamma\in \G$. By the definition of $\cC_Y$, we may assume $z=x+iy\in\bH$ with $y>Y$. Now let $\gamma=\left(\begin{smallmatrix}
* & *\\
a & b\end{smallmatrix}\right)\in \G$. If $a=0$, then $\Im(\gamma z)=\Im(z)>Y$. If $a\neq 0$, then
\begin{displaymath}
\Im(\gamma z)\ =\ \frac{\Im(z)}{|az+b|^2}\ =\ \frac{y}{(ax+b)^2+a^2y^2}\ \leq\ \frac{1}{y}\ <\ \frac{1}{Y}.\qedhere
\end{displaymath}
\end{proof}
%where $d_{\cM}$ is the distance function on $\cM$ induced from the hyperbolic distance function $d_{\bH}$ on $\bH$.
%and we denote by $W_{\rm pr}(\Psi)\subset [0,1)$ the set of primitive $\psi$-approximable numbers. 
%The following lemma reveals the connection between primitive $\psi$-approximable numbers and cusp excursions for certain collection of points on the modular surface.
%The classical Khinchine's theorem together with the monotonicity of the function $\psi$ implies the following zero-one law for the Lebesgue measure of $W_{\rm pr}(\psi)$. 
%The counter examples we will construct in the section relies on the following simple lemma which says that for any initial translate $x\in [0,1)$, if $x$ can be approximated well by a primitive rational number $m/n$, then the sample points we consider at the $n$-th step are all high into the cusp of the modular surface.
The following simple lemma is the key observation relating cusp excursions with Diophantine approximation.
\begin{Lem}\label{lem:counter1}
Let $x\in[0,1)$ be a real number. Suppose there exist a primitive rational number $m/n$ and $n>0$, and a real number $Y>0$ satisfying 
$$\left|x-\frac{m}{n}\right|<\frac{1}{2Yn^2}.$$
Then for any $0\leq j\leq n-1$ we have
\begin{equation}\label{equ:counter2}
\G\left(x+\tfrac{j}{n}+\tfrac{i}{2Yn^2}\right)\in \cC_{Y_j,2Y_j},\quad \text{ where }\ Y_{j}=\gcd(n,m+j)^2Y.
\end{equation}
%where $Y_{j}=\gcd(n,m+j)^2Y$. 
In particular, we have
\begin{equation}\label{equ:counter1}
%\cR_n\left(x,\tfrac{1}{2Yn^2}\right)\subset \cC_Y.
\left\{\G\left(x+\tfrac{j}{n}+\tfrac{i}{2Yn^2}\right):0\leq j\leq n-1\right\}\subset \cC_{Y}.
\end{equation}
%$$\left\{ \G  (x+\frac{j}{q}+i\frac{1}{2Yq^2}):1\leq j\leq q\right\}\subset \cC_Y.$$
\end{Lem}
\begin{proof}
The in particular part follows immediately from the inclusion
$\cC_{Y_j,2Y_j}\subset \cC_{Y},$
which in turn follows from the trivial bound $Y_j\geq Y$. Hence it suffices to prove the first half of the lemma.
%Since $x\in [0,1)$ is primitive $\psi$-approximable, there exists infinitely many pairs of primitive integers $(m,n)\in \Z\times \N$ such that \eqref{equ:dio} is satisfied. Moreover, for each $n\in \N$ since $0<\psi(n)< 1/2$, there exists at most one $m\in \N$ such that \eqref{equ:dio} is satisfied. Hence it suffices to show that if a primitive pair $(m,n)\in \Z\times \N$ satisfies \eqref{equ:dio}, then for any $1\leq j\leq n$ we have $\G\left(x+j/n+i\psi(n)/n\right)\in \widetilde{C}_{2Y_{n,j},Y_{n,j}}$ with $Y_{n,j}=\gcd(n,m+j)^2/(2n\psi(n))$ as in the lemma. For this we first simplify notation to set $r_n=\psi(n)/n$. %and we note that $r_n$ and $Y_n$ satisfy the relation $2r_nY_nn^2=1$. 
For simplicity of notation, we set $r=1/(2Yn^2)$. Then by assumption $|x-\tfrac{m}{n}|<r$. Fix $0\leq j\leq n-1$, and let $\tfrac{p}{q}$ be the reduced form of $\tfrac{m+j}{n}$ (so that $q=\tfrac{n}{\gcd(n,m+j)}$). Then
% we have for each $0\leq j\leq n-1$
%$$\left|x+j/n- (m+j)/n\right|< r.$$ 
%It is then clear from the geometry that this inequality is equivalent to the facts that 
$x+\tfrac{j}{n}+ir \in H^{\circ}_{p/q,r}$ and $x+\tfrac{j}{n}+ir'\in H_{p/q,r}$ for some $r<r'<2r$. Take $\gamma\in\G$ sending $H^{\circ}_{p/q,r}$ to the region $\left\{z\in\bH:\Im(z)>1/(2rq^2)=Y_j\right\}$. Then we have $\Im\left(\gamma(x+\tfrac{j}{n}+ir)\right)>Y_j$ and $\Im\left(\gamma(x+\tfrac{j}{n}+ir')\right)=Y_j$. Since $r<r'<2r$ we can bound the hyperbolic distance
$$d_{\bH}\left(\gamma(x+\tfrac{j}{n}+ir), \gamma(x+\tfrac{j}{n}+ir')\right)=\log\left(\tfrac{r'}{r}\right)<\log 2,$$
implying that 
$$\gamma(x+\tfrac{j}{n}+ir)\in\left\{z\in\bH:Y_j<\Im(z)<2Y_j\right\},$$
which implies \eqref{equ:counter2}. 
%Hence we have $\G\left(x+j/q+ir_n\right)\in \widetilde{C}_{Y'}$ with $Y'=1/(2rn^2)=(Yq^2)/n^2$. Since $p/q$ is the primitive form of $(m+j)/n$, we have $q | n$ implying that $n\geq q$, or equivalently, $Y'\geq Y$. Hence we have $\G\left(x+j/n+ir\right)\in \widetilde{C}_{Y'}\subset \cC_Y$, finishing the proof.
\end{proof}

\subsection{Full escape to the cusp along subsequences for almost every translate}
%For any initial translate of $x\in[0,1)$ \lemref{lem:counter1} says that if $x$ can be approximated well by a primitive rational number, then there are certain sample points lying in a certain cusp neighborhood of the modular surface. In this subsection we will use classical results from Diophantine approximation theory, namely the Khinchin's theorem, to ensure that for some nice choice of approximating function this phenomenon happens infinitely often for almost every initial translate $x\in[0,1)$. 
In this subsection we prove \thmref{thm:nonequ2}. Before stating this theorem, we first recall a definition from Diophantine approximation.  %We will prove \thmref{thm:nonequ2} using \lemref{lem:counter1} together with the classical Khinchine's theorem on Diophantine approximations. 
Let $\psi : \N\to (0,1/2)$ be a non-increasing function. We say that $x\in \R$ is \textit{primitive $\psi$-approximable} if there exist infinitely many $n\in\N$ such that the inequality
\begin{equation}\label{equ:dio}
\left|x-\frac{m}{n}\right|< \frac{\psi(n)}{n}
\end{equation}
is satisfied by some $m\in\Z$ coprime to $n$.  Since we assume $\psi(\N)\subset (0,1/2)$, the existence of such an $m$ implies its uniqueness. %Khinchin's theorem together with the monotonicity of the function $\psi$ implies that the Lebesgue measure of the set of primitive $\psi$-approximable numbers is full (resp. null) if and only if the series $\sum_{n\in\N}\psi(n)$ diverges (resp. converges). 
We prove the following:
\begin{Thm}\label{thm:nonequ2}
Let $\psi : \N\to (0,1/2)$ be a non-increasing function such that $\lim\limits_{n\to\infty}n\psi(n)=0$. Let $\{y_n\}_{n\in\N}$ be a sequence of positive numbers satisfying  
\begin{equation}\label{equ:condony}
r_n:=\frac12\min\{\psi(n)^{-2}y_n, n^{-2}y_n^{-1}\}\xrightarrow{n\to\infty}\infty.%:=\lim\limits_{n\to\infty}\psi(n)^{-2}y_n=\lim\limits_{n\to\infty}n^{-2}y^{-1}_n=\infty.
\end{equation}
If $x\in [0,1)$ is primitive $\psi$-approximable, then $\cR_n(x,y_n)\subset \cC_{r_n}$ infinitely often.  
\end{Thm}
\begin{rmk}\label{rmk:nonequid}
Since $\cR_n^{\rm pr}(x,y)\subset \cR_n(x,y)$ for any $n\in\N$, $x\in\R$ and $y>0$, \thmref{thm:nonequ2} also holds for translates of the primitive rational points. %implies that for any $x\in [0,1)$, as $n$ goes to infinity neither the sample points $\cR^{\rm pr}_n(x,y_n)$ nor $S_n(x,y_n)$ equidistribute on $\cM$ with respect to $\mu_{\cM}$.
%$\delta_{n, x,y_n} \nrightarrow \mu_{\cM}$
%as $n\to \infty$.
\end{rmk}
\begin{proof}[Proof of \thmref{thm:nonequ2}]
%In view of Khinchine's theorem and the assumption that $\sum_n\psi(n)=\infty$, it suffices to show the conclusion of \thmref{thm:nonequ2} holds for all primitive $\psi$-approximable numbers. Now assume that 
Let $x\in [0,1)$ be primitive $\psi$-approximable. Then for $Y_n=1/\left(2n\psi(n)\right)$, we have by \eqref{equ:counter1} that
\begin{equation}\label{equ:cond1}
\left\{\G\left(x+\tfrac{j}{n}+i\tfrac{\psi(n)}{n}\right)\in\cM:0\leq j\leq n-1\right\}\subset \cC_{Y_n}
\end{equation}
for infinitely many $n$'s. For every $n\in\N$, set $d_n:= Y_n/r_n = \max\left\{\psi(n)/(ny_n), ny_n/\psi(n)\right\}$. Then
\begin{equation}\label{equ:cond2}
d_{\bH}(t+i\psi(n)/n, t+iy_n)=\log(d_n)
\end{equation}
for any $t\in\R$. As in the proof of \lemref{lem:counter1}, by \eqref{equ:cond1} and \eqref{equ:cond2} we have $\cR_n(x,y_n)\subset \cC_{Y_n/d_n}$ for any $n$ in \eqref{equ:cond1}.
\end{proof}
%\begin{Thm}\label{thm:nonequ2}
%Keep the notation as above. Let $\psi : \N\to (0,1)$ be a non-increasing function such that $\sum_n\psi(n)=\infty$ and $\lim\limits_{n\to\infty}n\psi(n)=0$. Let $\{y_n\}_{n\in\N}$ be a sequence of positive numbers satisfying that 
%\begin{equation}\label{equ:condony}
%r_n:=\min\{\psi(n)^{-2}y_n, n^{-2}y_n^{-1}\}\xrightarrow{n\to\infty}\infty.%:=\lim\limits_{n\to\infty}\psi(n)^{-2}y_n=\lim\limits_{n\to\infty}n^{-2}y^{-1}_n=\infty.
%\end{equation}
%Then for almost every $x\in \R/\Z$ we have $\cR_n(x,y_n)\subset \cC_{r_n}$ infinitely often.  
%%$$\limsup_{n\to\infty}\inf_{z\in S_n(x,y_n)}d_{\cM}(z_0, z)=\infty.$$
%%n particular, for such sequence $\{y_n\}_{n\in\N}$ we have for almost every $x\in [0,1)$ there exists a subsequence $\cN_x$ $($which may depend on $x$$)$ such that $\delta_{n,x,y_n}\xrightarrow{w^*} 0$ as $n\in \cN_x$ goes to infinity.
%\end{Thm}

We now give a short 
\begin{proof}[Proof of \thmref{thm:nonequ2intro}]
Let $\alpha=\min\{\beta, 2-\beta\}$. For each $n\geq2$, let $\psi(n)=1/(n\log n)$ and let $\{y_n\}_{n\in\N}$ be a sequence of positive numbers satisfying $y_n\asymp 1/(n^2\log^{\beta} n)$. Then $r_n$ as in (\ref{equ:condony}) is given by $r_n=\tfrac12\min\{\psi(n)^{-2}y_n, n^{-2}y_n^{-1}\}\asymp \log^{\alpha} n$. By \thmref{thm:nonequ2}, for any $x\in [0,1)$ primitive $\psi$-approximable, we have that $\cR_n(x,y_n)\subset \cC_{r_n}$ infinitely often. Hence by \eqref{equ:distance}, for each such $x\in\R/\Z$, we have 
$$
\inf_{\G z\in \cR_n(x,y_n)}d_{\cM}(\G z_0, \G z)\geq \log(r_n)+O(1)=\alpha\log\log n+O(1)
$$ 
infinitely often, implying the inequality \eqref{equ:loglawdio}. Finally, since $\sum_{n\in\N}\psi(n)=\infty$ and $\psi$ is decreasing, the set of primitive $\psi$-approximable numbers in $[0,1)$ is of full measure by Khintchine's approximation theorem.
\end{proof}
%\begin{rmk}\label{rmk:counter1}
%It is easy to check that $\psi(n)=1/(n\log n)$ satisfies the conditions in \thmref{thm:nonequ2}. In this case the condition \eqref{equ:condony} on the sequence $\{y_n\}_{n\in\N}$ reads as 
%$$\lim\limits_{n\to\infty}\frac{1}{y_n n^2\log^2 n}=\lim\limits_{n\to\infty}n^2y_n=0.$$
%In particular, if $y_n\asymp 1/(n^2\log^{\alpha} n)$ for some fixed $0<\alpha< 2$, then the sequence $\{y_n\}_{n\in\N}$ satisfies the above condition. Hence the conclusion of \thmref{thm:nonequ2} holds for such sequence.
%$$\log^{\beta}n\leq e^{(1/2-\e)\log^{\alpha}n}$$ 
%eventually for all large enough $n\in\N$, we can see that for any fixed $\beta\in\R$, if $y_n\asymp n^{-2}\log^{\beta}n$, then the sequence $\{y_n\}_{n\in\N}$ satisfies the condition \eqref{equ:condony}. Hence the conclusion of \thmref{thm:nonequ2} holds for such $\{y_n\}_{n\in\N}$.
%\end{rmk}
For every irrational $x\in \R$, the \textit{Diophantine exponent} $\kappa_x>0$ is the supremum of $\kappa'>0$ for which $x$ is primitive $n^{-\kappa'}$-approximable. Dirichlet's approximation theorem implies that $\kappa_x\geq 1$ for any irrational $x$ and by Khintchine's theorem, $\kappa_x=1$ for almost every $x\in\R$. When $\kappa_x>1$, we have the following result that yields much faster cusp excursion rates for our sample points while handling sequences $\{y_n\}_{n\in\N}$ decaying polynomially faster than $1/n^2$.

\begin{Thm}\label{thm:negfas}
Let $\G z_0\in\cM$ be a fixed base point. Let $x\in [0,1)$ with Diophantine exponent $\kappa_x>1$ and let $\{y_n\}_{n\in\N}$ be a sequence of positive numbers satisfying $y_n\asymp n^{-\beta}$ for some fixed $2< \beta <2\kappa_x$. Then 
$$\limsup_{n\to\infty}\frac{\inf_{\G z\in \cR_n(x,y_n)}d_{\cM}\left(\G z_0, \G z\right)}{\log n}\geq \min\{2\kappa_x-\beta, \beta-2\}.$$
\end{Thm}

\begin{proof}
Take $\kappa\in (1,\kappa_x)$ and set $\alpha=\min\{2\kappa-\beta,\beta-2\}$. Let $\psi(n)=1/n^{\kappa}$. Then $x$ is primitive $\psi$-approximable since $\kappa<\kappa_x$. By \thmref{thm:nonequ2}, we have $\cR_n(x,y_n)\subset \cC_{r_n}$ infinitely often with $r_n=\tfrac12\min\{\psi(n)^{-2}y_n, n^{-2}y_n^{-1}\}\asymp n^{\alpha}$. This implies that 
$$\limsup_{n\to\infty}\frac{\inf_{\G z\in \cR_n(x,y_n)}d_{\cM}\left(\G z_0, \G z\right)}{\log n}\geq \alpha=\min\{2\kappa-\beta, \beta-2\}.$$
Taking $\kappa\to\kappa_x$ finishes the proof.
\end{proof}

\subsection{A non-equidistribution result for all translates} 
In this subsection we prove the following result which, together with part (1) of \thmref{thm:rationalnonequ} implies non-equidistribution for all translates: 
\begin{Thm}\label{thm:nonequ1}
Let $1/\sqrt{5}\leq c< 3/2$ and let $y_n=c/n^2$. 
%Then for any $x\in [0,1)$ there exists infinitely many $n\in\N$ such that $\c\cR_n(x,y_n)\subset \widetilde{C}_1$. In particular, for such sequence $\{y_n\}_{n\in\N}$, for any $x\in [0,1)$
%$$\delta_{x,y_n}:=\frac{1}{n}\sum_{j=1}^n\delta_{\G u_{x+j/n}a_{y_n}} \nrightarrow \mu_{\G}$$
%as $n\to \infty$.
Then there exists a closed measurable subset $\cE_c\subset \cM$, depending only on $c$, with $\mu_{\cM}(\cE_c)< 1$, and such that for each irrational $x\in[0,1)$, $\cR_n(x,y_n)\subset \cE_c$ infinitely often. %$$\left\{\G(x+j/n+iy_n)\in\cM:1\leq j\leq n\right\}\subset E_c.$$
%In particular, neither $\c\cR^{\rm pr}_{n_k}(x,y_{n_k})$ nor $\cS_{n_k}(x,y_{n_k})$ become equidistributed as $k\to\infty$. 
\end{Thm}

The set $\cE_c$ in \thmref{thm:nonequ1} is explicit: For any $c>0$, $\cE_c\subset \cM$ is defined to be the image of the closed set
 $$\left\{z\in\bH:\Im(z)\in [1/(2c),1/c]\cup [2/c,4/c]\cup [9/(2c),\infty)\right\}$$
under the natural projection from $\bH$ to $\cM$. It is clear from the definition that $\cE_c\subset \cM$ is closed. \thmref{thm:nonequ1} is a direct consequence of the following two lemmas.
\begin{Lem}\label{lem:prenonequ}
For any $c>0$ let $y_n=c/n^2$ and let $\psi_c(n)=c/n$. Then if $x\in [0,1)$ is primitive $\psi_c$-approximable, we have $\cR_n(x, y_n)\subset \cE_c$ infinitely often.
\end{Lem}
\begin{proof}
Let $x\in[0,1)$ be primitive $\psi_c$-approximable, that is, there exist infinitely many $n\in \N$ satisfying $\left|x-m/n\right|<c/n^2=y_n$ with some uniquely determined $m\in\Z$ satisfying $\gcd(m,n)=1$. For each such $n$, and for any $0\leq j\leq n-1$, let $k=\gcd(n,m+j)^2$. Then by \eqref{equ:counter2},   $\G(x+j/n+iy_n)\in \cC_{k^2/(2c),k^2/c}$.
%$$\G(x+j/n+iy_n)\in %\widetilde{C}_{2Y_{n,j},Y_{n,j}}=
%\widetilde{C}_{k^2/(2c),k^2/c},$$
Moreover, since $(k^2/(2c), k^2/c)\subset [1/(2c), 1/c]\cup [2/c, 4/c]\cup [9/(2c),\infty)$ for any $k\in\N$, we have $\cC_{k^2/(2c),k^2/c}\subset \cE_c$ for any $k\in\N$, implying that $\cR_n(x,y_n)\subset \cE_c$ for these infinitely many $n\in\N$.
\end{proof}
%Before giving the proof of \thmref{thm:nonequ1}, we first prove a simple lemma on Diophantine approximations using continued fractions. We refer the reader to \cite[Chapter 3]{EinsiedlerWard2011} for more details on continued fractions.
%\begin{Lem}\cite[Theorem 193 ]{HardyWright2008}\label{lem:dioapp}
%Let $\psi_c(q)=c/q$ for some $c>0$. If $c>1/\sqrt{5}$, then any irrational number in $[0,1)$ is primitive $\psi_c$-approximable.
%\end{Lem}

\begin{Lem}\label{lem:prenonequ2}
For any $0<c< 3/2$, we have $\mu_{\cM}(\cE_c)\leq 1-\frac{3}{\pi}\left(\frac{1}{\max\{2c,4/c\}}-\frac{2c}{9}\right)<1$.
\end{Lem}
\begin{proof}
Let $\cU\subset \cM$ be the projection of the open set
$$\left\{z\in\bH: \max\left\{2c,4/c\right\}<\Im(z)<9/(2c)\right\}.$$
Since $0<c<3/2$ we have $\max\{2c, 4/c\}<9/(2c)$ implying that $\cU$ is nonempty. We will show that $\cE_c$ is disjoint from $\cU$. Let $I_1= [1/(2c),1/c]$, $I_2= [2/c,4/c]$ and $I_3=[9/(2c),\infty)$, and for $1\leq j\leq 3$, define 
$\cE_c^j$ to be the projection onto $\cM$ of $\{z\in\bH:\Im(z)\in I_j\}$ 
%$$\cE_c^j:=\left\{\G z\in \cM:\Im(z)\in I_j\right\}$$
such that $\cE_c=\bigcup_{j=1}^3\cE_c^j$. 
It thus suffices to show that $\cE_c^j\cap \cU=\emptyset$ for each $1\leq j\leq 3$. For this, we identify (up to a null set) $\cM$ with the standard fundamental domain $\cF_{\G}:=\left\{z\in\bH:\Re(z)<\frac12, |z|>1\right\}$. Since $0<c<3/2$, we have $\max\left\{2c,4/c\right\}>2/c>2/(3/2)>1$. Thus we have
$$\cU=\left\{z\in\cF_{\G}: \max\left\{2c,4/c\right\}<\Im(z)<9/(2c)\right\},\qquad \cE_c^j=\left\{z\in\cF_{\G}:\Im(z)\in I_j\right\}$$
for $j= 2, 3$. Moreover, since the interval $( \max\left\{2c,4/c\right\}, 9/2c)$ intersects $I_{2}$ and $I_3$ trivially, we have $\cE_c^j\cap \cU=\emptyset$ for $j=2,3$. It thus remains to show that $\cE_c^1\cap \cU=\emptyset$. For this we note that $z\in\cF_{\G}$ satisfies the property that 
$$\Im(z)=\max_{\gamma\in \G}\Im(\gamma z).$$ 
Hence to show $\cE_c^1\cap \cU=\emptyset$, it suffices to show that 
$\max_{\gamma\in \G}\Im(\gamma z)\leq \max\left\{2c,4/c\right\}$ 
for any $z=s+it\in \bH$ with $\Im(z)=t\in I_1=[1/(2c), 1/c]$. For this, using the same discussion as in the proof of \lemref{lem:distance} we have
%For any $\gamma=\left(\begin{smallmatrix}
%* & * \\
%a & b\end{smallmatrix}\right)\in \G$, we have 
%$$\Im(\gamma z)=\frac{\Im(z)}{\left|az+b\right|^2}=\frac{1}{(as+b)^2t^{-1}+a^2t}.$$
%If $a=0$, then $b=\pm 1$ and $\Im(\gamma z)=t$. If $a\neq 0$, the trivial estimate $(as+b)^2t^{-1}+a^2t\geq t$ implies that $\Im(\gamma z)\leq t^{-1}$. To conclude we have 
for any $z=s+it\in\bH$ with $t\in [1/(2c), 1/c]$
\begin{displaymath}
\max_{\gamma\in\G}\Im(\gamma z)\leq \max\left\{t,t^{-1}\right\}\leq \max\left\{1/c,2c\right\}\leq \max\left\{2c,4/c\right\}.
\end{displaymath}
Finally, using the above description of $\cU$ and \eqref{equ:norhyarea} we have by direct computation
$$\mu_{\cM}(\cU)=\frac{3}{\pi}\left(\frac{1}{\max\{2c,4/c\}}-\frac{2c}{9}\right)$$ 
implying that $\mu_{\cM}(\cE_c)\leq 1-\frac{3}{\pi}\left(\frac{1}{\max\{2c,4/c\}}-\frac{2c}{9}\right)<1$ (again since $0<c<3/2$).
\end{proof}

\begin{proof}[Proof of \thmref{thm:nonequ1}]
Let $\psi_c(n)=c/n$. Since $c\geq 1/\sqrt{5}$, any irrational number is primitive $\psi_c$-approximable by the Hurwitz's approximation theorem; see, e.g., \cite[Theorem 193]{HardyWright2008}. Hence by \lemref{lem:prenonequ}, for each irrational $x\in[0,1)$, we have $\cR_n(x,y_n)\subset \cE_c$ infinitely often. Moreover, since $c<3/2$ by \lemref{lem:prenonequ2} we have $\mu_{\cM}(\cE_c)<1$, finishing the proof.
\end{proof}
\begin{rmk}\label{rmk:explain}
%Recall that $x\in\R$ is called \textit{badly approximable} if there exists some $c>0$ such that $|x-m/n|>c/n^2$ for any rational number $m/n$. 
The condition on the sequence $\{y_n\}_{n\in\N}$ in \thmref{thm:nonequ1} is quite restrictive and the proof of \thmref{thm:nonequ1} is much more involved than that of \thmref{thm:nonequ2}. We note that this is because we need to take care of the \textit{badly approximable numbers}, that is, the set of irrational numbers that are not primitive $\psi_c$-approximable for some $c>0$. If $x\in[0,1)$ is not badly approximable, then a similar argument as in the proof of \thmref{thm:nonequ2} using only the crude estimate \eqref{equ:counter1} would already be sufficient to prove non-equidistribution of the sample points $\cR_n(x,y_n)$ for any sequence $\{y_n\}_{n\in\N}$ satisfying $y_n\asymp 1/n^2$.
 \end{rmk}
%\thmref{thm:nonequ1} shows that if $y_n=c/n^2$ with $(\sqrt{5}-1)/2<c<3/2$, then for any initial translate $x\in[0,1)$ there is no equidistribution for the sample points 
%$$\left\{\G(x+j/n+iy_n)\in\G\bk \bH:1\leq j\leq n\right\}$$
%as $n$ goes to infinity. The next theorem says that if we only consider the distribution behavior of the sample points for almost every initial translate, then we can have much stronger results (full escape of mass along subsequences) which holds for much more general sequences $\{y_n\}_{n\in\N}$, see also \rmkref{rmk:counter1}.

%\begin{Thm}\label{thm:nonequ3}
%For any unbounded positive sequence $\{a_n\}_{n\in\N}$ and for a.e $x\in (0,1)$, there exists a sequence $\{y_n\}_{n\in \N}$ $($which may depend on $x$$)$ such that $1/y_{n}>a_n$ and $\cS_{n}(x,y_n)$ $($and hence also $\cS_{n}^{\rm pr}(x,y_n)$$)$ escapes to infinity as $n\to \infty$.
%\end{Thm}
%\begin{proof}
%
%\end{proof}
%The proofs of \thmref{thm:nonequ1} and \thmref{thm:nonequ2} rely on relations between hyperbolic geometry and Diophantine approximations, and we refer the reader to \cite{Sullivan1982,AthreyaMargulis09} for more details on such relations.

\section{Second moments of the discrepancy}\label{sec:secmom}
Let $\G=\SL_2(\Z)$ and let $\cM=\G\bk \bH$ be the modular surface as before. In this section we prove \thmref{thm:equipar}. Our proof relies on a second moment computation of the discrepancies $|\delta_{n,x,y}-\mu_{\cM}|$ and $|\delta_{n,x,y}^{\rm pr}-\mu_{\cM}|$ along the closed horocycle $\cH_y$. Throughout this section, we abbreviate the second moments $\int_0^1\left|\delta_{n,x,y}(\Psi)-\mu_{\cM}(\Psi)\right|^2dx$ and $\int_0^1\left|\delta^{\rm pr}_{n,x,y}(\Psi)-\mu_{\cM}(\Psi)\right|^2dx$ by $D_{n,y}(\Psi)$ and $D_{n,y}^{\rm pr}(\Psi)$ respectively. Since we assume $\G=\SL_2(\Z)$ we will also use the notation $\mu_{\G}$ for $\mu_{\cM}$.
%To avoid some difficulties, for the primitive discrepancy, we only consider the special case when $n$ is a prime number. While this specialization makes the subsequence $\cN$ for $\{\delta_{n,x,y_n}^{\rm pr}\}_{n\in\N}$ much sparser than that for $\{\delta_{n,x,y_n}\}_{n\in\N}$, it it still sufficient to prove the existence of such a subsequence.
%We now state the main result of this section.
%Namely, we are interested in computing $$\int_0^1\left|\frac{1}{n}\sum_{j=1}^n\Psi(\G u_{+j/n}a_y-\mu(\Psi)\right|^2dt,$$

\subsection{Relation to Hecke operators}
In this subsection we prove two preliminary estimates relating these second moments to the Hecke operators defined in \secref{sec:hecopr}.  %As mentioned before, we view $\Psi\in C_c^{\infty}(\cM)$ as a right $K$-invariant function on $\G\bk G$ and then the hyperbolic area $\mu_{\cM}$ coincides with the probability Haar measure $\mu_{\G}$. Moreover, since for any $x\in \R, y>0$ $u_x a_y\cdot i=x+iy$ we have (when viewing $\Psi$ as a function on $\G\bk G$)
%\begin{equation}\label{equ:di1}
%D_{n,y}(\Psi)=\int_0^1\left|\frac{1}{n}\sum_{j=0}^{n-1}\Psi(u_{x+j/n}a_y)-\mu_{\G}(\Psi)\right|^2dx,
%\end{equation}
%and
%\begin{equation}\label{equ:di2}
%D_{n,y}^{\rm pr}(\Psi)=\int_0^1\left|\frac{1}{\varphi(n)}\sum_{j\in (\Z/n\Z)^{\times}}\Psi(u_{x+j/n}a_y)-\mu_{\G}(\Psi)\right|^2dx.
%\end{equation}
%We now state the results relating second moments with Hecke operators.
\begin{Prop}\label{prop:secondmoment}
For any $n\in \N$, $y>0$ and $\Psi\in C_c^{\infty}(\cM)$, we have
\begin{equation}\label{equ:prehecra}
D_{n,y}(\Psi)=\frac{1}{n}\sum_{j=0}^{n-1}\left\langle \Psi_0, \widetilde{T}_{u_{j/n}}(\Psi_0)\right\rangle+O\left(\cS(\Psi)y^{1/2}\right),
\end{equation}
and
\begin{equation}\label{equ:prehecpri}
D_{n,y}^{\rm pr}(\Psi)\leq \frac{1}{\varphi(n)}\sum_{j=0}^{n-1}\left|\left\langle \Psi_0, \widetilde{T}_{u_{j/n}}(\Psi_0)\right\rangle\right|+O\left(\cS(\Psi)y^{1/2}\right).
\end{equation}
where $\Psi_0=\Psi-\mu_{\G}(\Psi)$, $\widetilde{T}_{u_{j/n}}$ is the Hecke operator associated to $u_{j/n}\in \SL_2(\Q)$ defined as in \eqref{equ:hecdef2}, the Sobolev norm $\cS(\Psi)$ is defined by
\begin{equation}\label{equ:Sobolev}
\cS(\Psi):=\cS_{4,2}^{\G}(\Psi)^2+\cS_{2,2}^{\G}(\Psi)\cS_{1,0}^{\G}(\Psi),
\end{equation}
and the implied constants are absolute. 
\end{Prop}

\begin{proof}%[Proof of \propref{prop:secondmoment}]
 Without loss of generality we may assume that $\Psi$ is real-valued. Expanding the square in the left hand side of \eqref{equ:prehecra}, doing a change of variables, and using the left $u_1$-invariance of $\Psi$, we have that $D_{n,y}(\Psi)$ equals
\begin{align*}
&\frac{1}{n^2}\sum_{j_1,j_2=0}^{n-1}\int_0^1\Psi(x+\tfrac{j_1}{n}+iy)\Psi(x+\tfrac{j_2}{n}+iy)dx-2\mu_{\G}(\Psi)\frac{1}{n}\sum_{j=0}^{n-1}\int_0^1\Psi(x+\tfrac{j}{n}+iy)dx+\mu_{\G}(\Psi)^2\\
&=\frac{1}{n}\sum_{j=0}^{n-1}\int_0^1\Psi(x+iy)\Psi(x+\tfrac{j}{n}+iy)dx-2\mu_{\G}(\Psi)\int_0^1\Psi(x+iy)dx+\mu_{\G}(\Psi)^2.
\end{align*}
Applying \eqref{equ:effequ} to the term $\int_0^1\Psi(x+iy)dx$ and using the trivial estimate 
\begin{equation}\label{equ:trivialestimate}
\|\Psi\|_2^{3/4}\|\Delta\Psi\|^{1/4}\left|\mu_{\G}(\Psi)\right|\leq \cS^{\G}_{2,2}(\Psi)\cS_{1,0}^{\G}(\Psi)\leq \cS(\Psi),
\end{equation} 
we get%the effective equidistribution \eqref{equ:effequ} of the long horocycle $\G Na_y\subset \G\bk G$ to $\Psi$ we get
\begin{equation}\label{equ:prediscom}
D_{n,y}(\Psi)=\frac{1}{n}\sum_{j=0}^{n-1}\int_0^1\Psi(x+iy)\Psi(x+\tfrac{j}{n}+iy)dx-\mu_{\G}(\Psi)^2+O(\cS(\Psi)y^{1/2}).
\end{equation}
%It thus remains to show that for each $0\leq j\leq n-1$
%$$\int_0^1\Psi(u_{x}a_y)\Psi(u_{x+j/n}a_y)dx=\left\langle \Psi, T_{u_{j/n}}(\Psi)\right\rangle+O\left(\cS(\Psi)y^{1/2}\right).$$
For each $0\leq j\leq n-1$, let $\G^j:=\G^{u_{j/n}}=\G\cap u_{j/n}^{-1}\G u_{j/n}$ and define  
%and $n_j:=n/(n,j)$. We note that since $\G u_{j/n}\G=\G\diag(1/n_j,n_j)\G$, %\G\left(\begin{smallmatrix}
%%(n,j)/n & 0\\
%%0 & n/(n,j)\end{smallmatrix}\right)\G$, 
%$\G^j$ is $\G$-conjugate to $\G_0(n_j^2)$. %Now we view $\Psi\in C_c^{\infty}(\G\bk G)$ as a left $\G$-invariant function on $G$. In particular 
%For each $1\leq j\leq n$ 
 $F_j(\Psi):=\Psi L_{u_{j/n}^{-1}}\Psi\in C^{\infty}(\bH)$. Since $\Psi$ is left $\G$-invariant, 
%$L_{u_{j/n}}\Psi$ is left $u_{j/n}^{-1}\G u_{j/n}$-invariant. In particular, 
and $L_{u_{j/n}^{-1}}\Psi$ is left $u_{j/n}^{-1}\G u_{j/n}$-invariant, we have $F_j(\Psi)\in C^{\infty}(\G^j\bk \bH)$. Moreover,
$$F_j(\Psi)(x+iy)=\Psi(x+iy)\Psi(x+\tfrac{j}{n}+iy).$$ 
For each $0\leq j\leq n-1$, it is easy to check that $u_1\in \G^j$ and $\G^j$ contains the principal congruence subgroup $\G(n^2)$, hence $\G^j$ satisfies the assumptions in \propref{prop:equiclohor}. Then by \eqref{equ:effequ}, 
$$
\int_0^1 F_j(\Psi)(x+iy)dx = \int_{\Gamma^j\backslash\bH} F_j(\Psi)(z) d\mu_{\Gamma^j}(z) + O\left(\|F_j(\Psi)\|_2^{3/4}\|\Delta F_j(\Psi)\|_2^{1/4} y^{1/2}\right).
$$
%$$
%\int_0^1\Psi(x+iy)\Psi(x+\tfrac{j}{n}+iy)dx=\left\langle \Psi, L_{u_{j/n}^{-1}}\Psi\right\rangle_{L^2(\G^j\bk \bH)}+O\left(\|F_j(\Psi)\|_2^{3/4}\|\Delta F_j(\Psi)\|_2^{1/4} y^{1/2}\right).
%$$
Next we note that by \eqref{equ:sobolev},
$$\|F_j(\Psi)\|_2^{3/4}\|\Delta F_j(\Psi)\|_2^{1/4}\leq \cS^{\G^j}_{2,2}\left(F_j(\Psi)\right)=\cS^{\G^j}_{2,2}\left(\Psi L_{u_{j/n}^{-1}}\Psi\right)\leq \cS_{4,2}^{\G^j}\left(\Psi\right)\cS_{4,2}^{\G^j}\left(L_{u_{j/n}^{-1}}\Psi\right).$$
Using the fact that $\Psi$ is left $\G$-invariant and $\G^j$ is a finite-index subgroup of $\G$, by \eqref{equ:soblev3}, $\cS_{4,2}^{\G^j}(\Psi)=\cS_{4,2}^{\G}(\Psi)$. Similarly, we have 
$$\cS_{4,2}^{\G^j}\left(L_{u_{j/n}^{-1}}\Psi\right)=\cS_{4,2}^{u_{j/n}^{-1}\G u_{j/n}}\left(L_{u_{j/n}^{-1}}\Psi\right)=\cS_{4,2}^{\G}\left(\Psi\right),$$ 
where for the second equality we used \eqref{equ:sobolevconj}. Hence we have
\begin{equation}\label{equ:soblevcomp}
\|F_j(\Psi)\|_2^{3/4}\|\Delta F_j(\Psi)\|_2^{1/4}\leq\cS^{\G^j}_{2,2}\left(F_j(\Psi)\right)\leq \cS_{4,2}^{\G}(\Psi)^2\leq \cS(\Psi)<\infty.
\end{equation}
%Thus,
%\begin{align*}
%D_{n,y}(\Psi) = \frac{1}{n}\sum_{j=0}^{n-1} \int_{\Gamma^j\backslash\bH} \Psi_0(z)\Psi_0(u_{j/n}z)d\mu_{\Gamma^j}(z) + O\left(\mathcal{S}(\Psi)y^{1/2}\right).
%\end{align*}
Thus applying \eqref{equ:effequ} to $F_j\in C^{\infty}(\G^j\bk \bH)$ and using \eqref{equ:soblevcomp} we get
\begin{equation}\label{equ:suff}
\int_0^1\Psi(x+iy)\Psi(x+\tfrac{j}{n}+iy)dx=\left\langle \Psi, L_{u_{j/n}^{-1}}\Psi\right\rangle_{L^2(\G^j\bk \bH)}+O\left(\cS(\Psi)y^{1/2}\right).
\end{equation}
%\begin{align*}
%\int_0^1\left|D_n(\Psi,y)\right|^2dt&=\frac{1}{n}\sum_{j=1}^n\left(\mu_{\G^j}(F_j)+O(\cS_{p,d}^{\G^j}(\Psi L_{u_{j/n}}\Psi)y^{1/2})\right)-\mu(\Psi)^2+O(\cS_{p,d}(\Psi)y^{1/2})\\
%&=\frac{1}{n}\sum_{j=1}^n\left(\mu_{\G^j}(F_j)+O(\cS_{2p,d}^{\G^j}(\Psi)\cS_{2p,d}^{\G^j}(L_{u_{j/n}}\Psi) y^{1/2})\right)-\mu(\Psi)^2+O(\cS_{p,d}(\Psi)y^{1/2})\\
%&=\frac{1}{n}\sum_{j=1}^n\int_{\G^j\bk G}\Psi(g)\Psi(u_{j/n}g)d\mu_{\G^j}(g)-\mu(\Psi)^2+O\left((\cS_{2p,d}(\Psi)^2+\cS_{p,d}(\Psi))y^{1/2}\right),
%\end{align*} 
%where for the second equality we applied \eqref{equ:sobolev}, and for the last equality we used the identity that 
%$$\cS_{2p,d}^{\G^j}(L_{u_{j/n}}\Psi)=\cS_{2p,d}^{\G^j}(\Psi)=\cS_{2p,d}(\Psi).$$

Plugging \eqref{equ:suff} into \eqref{equ:prediscom} and using the identities $\mu_{\G}(\Psi)=\mu_{\G^j}(\Psi)=\mu_{\G^j}(L_{u_{j/n}^{-1}}\Psi)$ (the second equality follows from the left $G$-invariance of the hyperbolic area $\mu_{\G^j}$) we get that
$$D_{n,y}(\Psi)=\frac{1}{n}\sum_{j=0}^{n-1}\left\langle \Psi_0, L_{u_{j/n}^{-1}}\Psi_0\right\rangle_{L^2(\G^j\bk \bH)}+O(\cS(\Psi)y^{1/2}).$$

%It thus remains to show that for each $0\leq j\leq n-1$
%\begin{equation}\label{equ:suff}
%\left\langle \Psi_0, L_{u_{j/n}^{-1}}\Psi_0\right\rangle_{L^2(\G^j\bk \bH)}=\left\langle \Psi_0, \widetilde{T}_{u_{j/n}}(\Psi_0)\right\rangle_{L^2(\G\bk \bH)}.
%\end{equation}
Let $\cF_{\G}\subset \bH$ be a fundamental domain for $\G\bk \bH$. The disjoint union $\bigsqcup_{\gamma\in \G^j\bk \G}\gamma \cF_{\G}$ forms a fundamental domain for $\G^{j}\bk \bH$. Thus we can conclude the proof of \eqref{equ:prehecra} by noting that
\begin{align*}
&\int_{\bigsqcup_{\gamma\in \G^j\bk \G}\gamma \cF_{\G}}\Psi_0(z)\Psi_0(u_{j/n}z)d\mu_{\G^j}(z)=\sum_{\gamma\in \G^j\bk \G}\int_{\gamma\cF_{\G}}\Psi_0(z)\Psi_0(u_{j/n}z)d\mu_{\G^j}(z)\\
&=\int_{\cF_{\G}}\Psi_0(z)\left(\frac{1}{[\G: \G^j]}\sum_{\gamma\in \G^j\bk \G}\Psi_0(u_{j/n}\gamma z)\right)d\mu_{\G}(z)=\int_{\cF_{\G}}\Psi_0(z)\widetilde{T}_{u_{j/n}}(\Psi_0)(z)d\mu_{\G}(z),
\end{align*}
where for the second equation we did a change of variable $z\mapsto \gamma z$, used the left $\G$-invariance of $\Psi$ and the relation $[\G:\G^j]\mu_{\G^j}=\mu_{\G}$, and for the last equality we used the expression \eqref{equ:Hecke}.
Similarly, applying the estimates \eqref{equ:effequ} and \eqref{equ:trivialestimate} and making change of variables we see that $D_{n,y}^{\rm pr}(\Psi)$ equals
\begin{align*}
&\frac{1}{\varphi(n)^2}\sum_{j_1,j_2\in (\Z/n\Z)^{\times}}\int_0^1\Psi(x+\tfrac{j_1}{n}+iy)\Psi(x+\tfrac{j_2}{n}+iy)dx-\mu_{\G}(\Psi)^2+O\left(\cS(\Psi)y^{1/2}\right)\\
&=\frac{1}{\varphi(n)^2}\sum_{j=0}^{n-1}c(j)\int_0^1\Psi(x+iy)\Psi(x+\tfrac{j}{n}+iy)dx-\mu_{\G}(\Psi)^2+O\left(\cS(\Psi)y^{1/2}\right),
\end{align*}
where 
$$c(j):=\#\left\{([j_1], [j_2])\in (\Z/n\Z)^{\times}\times (\Z/n\Z)^{\times}: [j_2]-[j_1]=[j]\right\}.$$
Now similar as before we can apply the estimate \eqref{equ:suff}, the identities $\mu_{\G}(\Psi)=\mu_{\G^j}(\Psi)=\mu_{\G^j}(L_{u_{j/n}^{-1}}\Psi)$ and $\sum_{j=0}^{n-1}c(j)=\varphi(n)^2$ to get
\begin{align*}
D_{n,y}^{\rm pr}(\Psi)&=\frac{1}{\varphi(n)^2}\sum_{j=0}^{n-1}c(j)\left\langle \Psi_0, L_{u_{j/n}^{-1}}\Psi_0\right\rangle_{L^2(\G^j\bk \bH)}+O(\cS(\Psi)y^{1/2})\\
&=\frac{1}{\varphi(n)^2}\sum_{j=0}^{n-1}c(j)\left\langle \Psi_0, \widetilde{T}_{u_{j/n}}(\Psi_0)\right\rangle_{L^2(\G\bk\bH)}+O(\cS(\Psi)y^{1/2}).
\end{align*}
Finally we can finish the proof by noting that for each $0\leq j\leq n-1$, $c(j)\leq \varphi(n)$ (since for each $[j_1]\in (\Z/n\Z)^{\times}$, there is at most one $[j_2]\in (\Z/n\Z)^{\times}$ such that $[j_2]-[j_1]=[j]$). 
%We note that as the pair $([j_1],[j_2])$ runs through $\left(\Z/p\Z\right)^{\times}\times \left(\Z/p\Z\right)^{\times}$, for any $[j]\in \Z/p\Z$, $[j_1-j_2]=[j]$ for $p-1$ times if $[j]=[0]$ and for $p-2$ times otherwise. Hence we have $D_{p,y}^{\rm pr}(\Psi)$ equals
%\begin{align*}
%\frac{1}{p-1}\int_0^1\Psi(x+iy)^2dx+\frac{p-2}{(p-1)^2}\sum_{j=1}^{p-1}\int_0^1\Psi(x+iy)\Psi(x+\frac{j}{p}+iy)dx-\mu_{\G}(\Psi)^2+O(\cS(\Psi)y^{1/2}).
%\end{align*}
%The estimate \eqref{equ:prehecpri} then follows by using the exact same arguments as for \eqref{equ:prehecra} and noting that $\widetilde{T}_{u_{j/p}}=\widetilde{T}_{p}$ for all $1\leq j\leq p-1$ (since for these $j$, the matrices $u_{j/p}$ all have degree $p$). 
%we get $D_{p,y}^{\rm pr}(\Psi)$ equals
%$$\frac{1}{p-1}\langle \Psi, \Psi\rangle +\frac{p-2}{p-1}\langle \Psi, T_{[p]}\Psi\rangle -\mu_{\G}(\Psi)^2+O(\cS(\Psi)y^{1/2}),$$
%finishing the proof.
\end{proof}
\subsection{Second moment estimates}
Combining \propref{prop:secondmoment} and the operator norm bound in \propref{prop:hecbound2} we have the following second moment estimates:
%We now give the proof of \thmref{thm:secmom} using \propref{prop:secondmoment} and the bound for Hecke operators in \propref{prop:hecbound2}.
\begin{Thm}\label{thm:secmom}
For any $n\in \N$, $y>0$ and $\Psi\in C_c^{\infty}(\cM)$ we have
\begin{equation}\label{equ:disest}
\max\left\{D_{n,y}(\Psi), D^{\rm pr}_{n,y}(\Psi)\right\}\ll_{\e}n^{-1+2\theta+\e}\|\Psi\|_2^2+\cS(\Psi)y^{1/2},
\end{equation}
%and
%\begin{equation}\label{equ:disestpri}
%\int_0^1\left|\delta^{\rm pr}_{n,x,y}(\Psi)-\mu_{\cM}(\Psi)\right|^2dx\ll_{\e}n^{-1+2\theta+\e}\|\Psi\|_2^2+\cS(\Psi)y^{1/2}.
%\end{equation}
where $\theta=7/64$ is the best bound towards the Ramanujan conjecture as before and the Sobolev norm $\cS(\Psi)$ is as defined in \eqref{equ:Sobolev}.
%where $\theta=7/64$ is before and 
%$$\cS(\Psi):=\cS^{\G}_{4,4}(\Psi)^2+\cS^{\G}_{2,4}(\Psi)\cS^{\G}_{1,0}(\Psi).$$ 
%Here for any $p>0$ and $d\in\N$, $\cS_{p,d}^{\G}$ is some ''$L^p$, order $d$'' Sobolev norm  defined on $C_c^{\infty}(\G\bk G)$, see \secref{sec:soblev}.
%In particular, for any sequence $\{y_n\}_{n\in \N}$ with $0<y_n\leq C_n^{-1}n^{-1/(1-2\theta)}$ we have
%$$\int_0^1\left|\delta_{x,y_n}(\Psi)-\mu_{\G}(\Psi)\right|^2dx\ll_{\e,\Psi}n^{-1/2+\e}.$$
\end{Thm} 
\begin{rmk}
It is also possible to approach the second moment computation using the spectral bounds on the Fourier coefficients of $\Psi$ from \secref{sec:fouriercoe} rather than Hecke operators. The spectral approach however yields a weaker estimate when $y>0$ is small. For comparison, following the spectral approach, one obtains
$$
\int_0^1 |\delta_{n,x,y}(\Psi)-\mu_\Gamma(\Psi)|^2 dx\ \ll_\e \left(n^{-1}y^{-2(\theta+\e)}+y^{1/2}\right)\mathcal{S}_{2,2}(\Psi).
$$
\end{rmk}

\begin{proof}[Proof of \thmref{thm:secmom}]
First we prove \eqref{equ:disest}. For each $0\leq j\leq n-1$, it is clear that $u_{j/n}$ is of degree $n_j:=n/\gcd(n,j)$, and thus $\widetilde{T}_{u_{j/n}}=\widetilde{T}_{n_j}$. 
%Then by \propref{prop:secondmoment} we have
%\begin{align*}
%D_{n,y}(\Psi)&=\frac{1}{n}\sum_{j=0}^{n-1}\left\langle \Psi, T_{[n_j]}(\Psi)\right\rangle-\left|\mu_{\G}(\Psi)\right|^2+O\left(\cS(\Psi)y^{1/2}\right)\\
%%&=\frac{1}{n}\sum_{j=1}^n\left\langle \Psi, T^0_{[n_j]}(\Psi)\right\rangle-\left|\mu(\Psi)\right|^2+O\left(\left(\cS_{2p,d}(\Psi)^2+\cS_{p,d}(\Psi)\right)y^{1/2}\right)\\
%&=\frac{1}{n}\sum_{j=0}^{n-1}\left\langle\Psi_0, T^0_{[n_j]}\Psi_0\right\rangle +O\left(\cS(\Psi)y^{1/2}\right),
%\end{align*} 
%where for the second equality we used the identity that $\mu_{\G}(T_{[n_j]}(\Psi))=\mu_{\G}(\Psi)$ which follows from the left $G$-invariance of the Haar measure $\mu_{\G}$. 
Applying \eqref{equ:prehecra}, \eqref{equ:prehecpri}, the estimate $\varphi(n)\gg_{\e} n^{-1+\e/2}$ and the operator norm bound in \propref{prop:hecbound2} to the terms $\left\langle\Psi_0, \widetilde{T}_{n_j}\Psi_0\right\rangle$, we get 
\begin{align*}
\max\left\{D_{n,y}(\Psi), D^{\rm pr}_{n,y}(\Psi)\right\}&\ll_{\e}\ n^{-1+\e/2}\sum_{j=0}^{n-1}n_j^{-1+2\theta+\e/4}\|\Psi_0\|_2^2+\cS(\Psi)y^{1/2}.
\end{align*}
For any $d\mid n$, $\#\{0\leq j\leq n-1:n_j=d\}=\varphi(d)$, thus
\begin{align*}
\sum_{j=1}^nn_j^{-1+2\theta+\e/4}=\sum_{d| n}\varphi(d)d^{-1+2\theta+\e/4}< \sum_{d\mid n}d^{2\theta+\e/4}=\sigma_{2\theta+\e/4}(n)\ll_{\e}n^{2\theta+\e/2},
\end{align*}
where for the first inequality we used the trivial bound $\varphi(d)<d$. Finally, we observe that $\|\Psi_0\|_2\leq \|\Psi\|_2$. 
%In view of the inequality $\|\Psi_0\|_2\leq \|\Psi\|_2$, it suffices to show that 
%$$\sum_{j=0}^{n-1}n_j^{-1+2\theta+\e/2}\ll_{\e}n^{2\theta+\e}.$$
%For this noting that for any $d| n$, $\#\{0\leq j\leq n-1:n_j=d\}=\varphi(d)$ we get%as $j$ runs over $\{1,2,\ldots, n\}$ there are exactly $\varphi(n/d)$ have
%\begin{align*}
%\sum_{j=1}^nn_j^{-1+2\theta+\e/2}=\sum_{d| n}\varphi(d)d^{-1+2\theta+\e/2}< \sum_{d\mid n}d^{2\theta+\e/2}=\sigma_{2\theta+\e/2}(n)\ll_{\e}n^{2\theta+\e},
%\end{align*}
%finishing the proof,  %and for the last estimate we used the bound that $\sigma_{1/2+\e/2}(n)\ll_{\e}n^{1/2+\e/2}$. 
%for the sum of positive divisor function $\sigma_x(n):=\sum_{d| n}d^x$. We note that this estimate follows easily from the multiplicative expression of $\sigma_x(n)$; see e.g. \cite[]{HardyWright2008}.
%Now for \eqref{equ:disestpri}, we have by \eqref{equ:prehecpri} and \propref{prop:hecbound2} that
%$$D^{\rm pr}_{p,y}(\Psi)\ll_{\e}\left(\frac{p-2}{p-1}p^{-1+2\theta+\e}+\frac{1}{p-1}\right)\|\Psi_0\|_2^2+\cS(\Psi)y^{1/2}\ll p^{-1+2\theta+\e}\|\Psi\|^2_2+\cS(\Psi)y^{1/2},$$
%finishing the proof.
\end{proof}

We now give a quick
\begin{proof}[Proof of \thmref{thm:equipar}]
%We first prove the existence of such a sequence $\cN\subset \N$ for the sequence $\left\{\delta_{n,x,y_n}\right\}_{n\in\N}$. 
{{Let $\alpha>0$ be the fixed number as in this theorem. Let $\beta:=\min\{\frac{\alpha}{2}, 1-2\theta\}$.}}
 {{Fix $0<c< \beta$ and}} let $\cN\subset \N$ be an unbounded subsequence such that $\sum_{n\in\cN}n^{-c}<\infty$. We want to show that for any $\{y_n\}_{n\in\N}$ satisfying {{$y_n\ll n^{-\alpha}$}} there exists a full measure subset $I\subset \R/\Z$ such that for any $x\in I$, $\delta_{n,x,y_n}(\Psi)\to \mu_{\cM}(\Psi)$ and $\delta^{\rm pr}_{n,x,y_n}(\Psi)\to \mu_{\cM}(\Psi)$ for any $\Psi\in C_c^{\infty}(\cM)$ as $n\in\cN$ goes to infinity. Since the function space $C_c^{\infty}(\cM)$ has a dense countable subset, it suffices to prove the above assertion for a fixed $\Psi$. Now we fix $\Psi\in C_c^{\infty}(\cM)$ and take $\e>0$ sufficiently small such that {{$\beta-2\e>c$}}. For any $n\in\N$ define $I_{n}=I_n^1\cup I_n^2\subset \R/\Z$ such that
$$I_n^1:=\left\{x\in \R/\Z :\left|\delta_{n,x,y_n}(\Psi)-\mu_{\cM}(\Psi)\right|>n^{-\e/2}\right\},$$
and 
$$I_n^2:=\left\{x\in \R/\Z :\left|\delta^{\rm pr}_{n,x,y_n}(\Psi)-\mu_{\cM}(\Psi)\right|>n^{-\e/2}\right\}.$$
Thus by the second moment estimate \eqref{equ:disest}, the assumption that {{$y_n\ll n^{-\alpha}$}} and Chebyshev's inequality we get
%$$\int_0^1\left|\delta_{n,x,y}(\Psi)-\mu_{\cM}(\Psi)\right|^2dx\ll_{\e, \Psi}n^{-1/2+\e}$$
\begin{align*}
\left| I_n\right|\leq \left|I_n^1\right|+\left| I_n^2\right|\leq 2n^{\e}\max\left\{D_{n,y}(\Psi), D^{\rm pr}_{n,y}(\Psi)\right\}\ll_{\e,\Psi}{{n^{-\beta+2\e}<n^{-c}}},
\end{align*}
implying that $\sum_{n\in\cN}|I_n|<\infty$. Hence taking $I\subset \R/\Z$ to be the complement of this limsup set $\limsup_{\substack{n\in\cN\\ n\to\infty}}I_n\subset \R/\Z$ and by the Borel-Cantelli lemma we have $I$ is of full measure. Moreover, for any $x\in I$, $x\in I_n^c$ for all $n\in\cN$ sufficiently large, that is,
$$\max\left\{\left|\delta_{n,x,y_n}(\Psi)-\mu_{\cM}(\Psi)\right|, \left|\delta^{\rm pr}_{n,x,y_n}(\Psi)-\mu_{\cM}(\Psi)\right|\right\}\leq n^{-\e/2},\quad \forall\ n\in\cN\ \text{sufficiently large}.$$ 
In particular for such $x$, $\delta_{n,x,y_n}(\Psi)\to \mu_{\cM}(\Psi)$ and $\delta^{\rm pr}_{n,x,y_n}(\Psi)\to \mu_{\cM}(\Psi)$ as $n\in\cN$ goes to infinity.
%For the sequence $\left\{\delta_{n,x,y_n}^{\rm pr}\right\}_{n\in\n}$, let $\cN$ be an unbounded subsequence of the set of prime numbers such that $\sum_{n\in\cN}n^{-\alpha}<\infty$ for some $0<\alpha<1-2\theta$. The remaining part of the proof is then identical to that for the sequence $\left\{\delta_{n,x,y_n}\right\}_{n\in\N}$ with the estimate \eqref{equ:disest} replaced by \eqref{equ:disestpri}.
\end{proof}
\begin{rmk}\label{rmk:hejhal}
The second moment $D_{n,y}(\Psi)$ is closely related to the sample points \eqref{equ:samplehe} considered in \cite{Hejhal1996}: 
%For simplicity let us assume $\Psi$ is of mean zero, that is $\mu_{\G}(\Psi)=0$. 
Using the extra invariance $\delta_{n,x+1/n,y}(\Psi)=\delta_{n,x,y}(\Psi)$ and applying a change of variable, one can easily check that
$$D_{n,y}(\Psi)=\int_0^1\left|\frac1n\sum_{j=0}^{n-1}\Psi\left(\tfrac{x+j}{n}+iy\right)-\mu_{\G}(\Psi)\right|^2dx.$$
Thus let $\cN\subset \N$ be the fixed sequence as in the above proof, by \thmref{thm:secmom} and the same Borel-Cantelli type argument we have that for almost every $x\in \R/\Z$ the sequence of sample points 
$\{\G(\tfrac{x+j}{n}+iy_n: 0\leq j\leq n-1\}$
equidistributes on $\cM$ with respect to $\mu_{\cM}$ {{as $n\in\cN$ goes to infinity, as long as $\{y_n\}_{n\in\N}$ decays at least polynomially}}. 
%where $S_{n,y,\Psi}(x)$ is given as in \eqref{equ:hejhal}. In particular, \eqref{equ:prehecra} implies that $\int_0^1\left|n^{-1/2}S_{n,y,\Psi}(x)\right|^2dx$ is governed by the term $\sum_{j=0}^{n-1}\left\langle \Psi, \widetilde{T}_{u_{j/n}}(\Psi)\right\rangle$ provided that $ny^{1/2}$ is small. While the heuristic arguments in \cite[(4.6)]{Hejhal1996} shows that
%$$\int_0^1\left|n^{-1/2}S_{n,y,\Psi}(x)\right|^2dx=\mu_{\cM}(|\Psi|^2)+O(ny\log^2 n)+O(\sqrt{ny})$$
%provided that the points $\{\G(\tfrac{x+j}{n}+iy): 0\leq j\leq n-1\}$ become decorrelated when $ny$ is small. In particular, \eqref{equ:prehecra} implies that this decorrelation process fails for $\G=\SL_2(\Z)$.
\end{rmk}

\section{Left regular action of normalizing elements}\label{sec:leftregular}
In this section, $\G$ denotes a congruence subgroup, and we set by $\G_1=\SL_2(\Z)$. We moreover assume that there exists some $h\in\SL_2(\Q)$ normalizing $\G$, that is, $h^{-1}\G h=\G$. It induces the left regular $h$-action on $\G\bk \bH$ given by $\G z\in \G\bk \bH\mapsto \G hz\in \G\bk\bH$. Since $h$ normalizes $\G$, this map is well defined: Suppose $\G z=\G z'$, that is there exists some $\gamma\in \G$ such that $z'=\gamma z$. Then $\G hz'=\G h\gamma z=\G h\gamma h^{-1}hz=\G hz$. The goal of this section is to describe this action on cylindrical cuspidal neighborhoods of $\Gamma\backslash\bH$.

\subsection{Cusp neighborhoods of congruence surfaces}
Since $\G$ is a congruence subgroup, the set of cusps of $\G$ can be parameterized by the coset $\G\bk \left(\Q\cup \{\infty\}\right)$ (see e.g. \cite[p. 222]{LangSL2R}), where the action of $\G$ on $\Q\cup \{\infty\}$ is defined via the M\"{o}bius transformation. We denote by $\Omega_{\G}$ a complete list of coset representatives for $\G\bk \left(\Q\cup \{\infty\}\right)$. For each cusp representative $\fc\in\Omega_{\G}$, its stabilizer subgroup
%\footnote{ More precisely, $\G_{\fc}$ is an index two subgroup of the stabilizer subgroup if $-I_2\in \G$.} 
is given by 
$$\G_{\fc}:=\tau_{\fc}N\tau_{\fc}^{-1}\cap \G,$$ 
where $\tau_{\fc}\in \G_1$ is such that $\tau_{\fc}\infty=\fc$. 
{{(More precisely, $\G_{\fc}$ is an index two subgroup of the stabilizer subgroup of $\fc$ if $-I_2\in \G$.)}}
The existence of such $\tau_{\fc}$ is guaranteed by the transitivity of the action of $\G_1$ on $\Q\cup\{\infty\}$. On the other hand, $\tau_{\fc}$ is only unique up to right multiplication by any element of $\pm N$. We note that $\G_{\fc}$ is independent of the choice of $\tau_{\fc}$, and since $\fc\in \Omega_{\G}$ is a cusp, $\G_{\fc}$ is nontrivial. Moreover, $\tau_{\fc}^{-1}\G_{\fc}\tau_{\fc}$ is a subgroup of $N\cap \G_1=\langle u_1\rangle$. Hence $\tau_{\fc}^{-1}\G_{\fc}\tau_{\fc}$ is a cyclic group generated by a unipotent matrix $u_{\omega_{\fc}}$ for some positive integer $\omega_{\fc}$, which is called the \textit{width} of the cusp $\fc$. 

We can now define cusp neighborhoods on the hyperbolic surface $\G\bk \bH$ around a cusp $\fc\in\Omega_\G$. For any $Y>0$,
%$$
%\cY_{\Gamma,\fc}(z) := \max_{\gamma\in\Gamma} \Im(\tau_\fc^{-1}\gamma z).
%$$
%This function is continuous, $\Gamma$-invariant, and is bounded below by a positive constant depending on $\G$ and $\fc$, see \cite[(3.8)]{Iwaniec2002}. For any $Y, Y'>1$ and for any cusp $\fc\in \Omega_{\G}$, let
%$$
%\cC_Y^{\Gamma,\frak{c}}=\left\{\Gamma z: \cY_{\Gamma,\fc}(z)>Y\right\},\qquad \cC_{Y,Y'}^{\G,\fc}:=\{\Gamma z: Y< \cY_{\Gamma,\fc}(z)<Y'\}.
%$$
$\cC_Y^{\Gamma,\fc}\subset\G\backslash\bH$ denote the projection of the horodisc $\{\tau_\fc z\in\bH: \Im(z)>Y\}$ onto $\G\backslash\bH$. Similarly, for any $Y'>Y>0$, let $\cC_{Y,Y'}^{\G,\fc}$ denote the projection of the cylindrical region $\{\tau_\fc z\in\bH:Y<\Im(z)<Y'\}$ onto $\G\backslash\bH$. We record the following two lemmas for the later purpose of computing the measure of certain unions of cusp neighborhoods.
\begin{Lem}\label{lem:onetoone}
If $Y'>Y>1$, the set $\cC_{Y,Y'}^{\G,\fc}$ is in one-to-one correspondence with the set 
\begin{equation}\label{equ:idenset}
\{\tau_c z\in\bH: \Re(z)\in\R/\omega_\fc\Z,\ \Im(z)\in(Y,Y')\}.
\end{equation}
%For each $\Gamma z\in\Gamma\backslash\bH$, $\G z\in \cC_{Y,Y'}^{\Gamma,\fc}$ if and only if there exists  $z'$ in
%$$\{x+iy\in\bH: 0\leq x<\omega_\fc,\ Y<y<Y'\}$$ such that $\Gamma \tau_\fc z'=\Gamma z$.
In particular, if $-I_2\in\G$ then for any $Y'>Y>1$
\begin{align}\label{equ:volume}
\mu_\Gamma\left(\cC_{Y,Y'}^{\Gamma,\fc}\right) = \frac{3\omega_\fc}{\pi[\Gamma_1:\Gamma]}\left(\frac{1}{Y}-\frac{1}{Y'}\right).
\end{align}
\end{Lem}
\begin{proof}
The one-to-one correspondence is given by the projection of the above rectangular set onto $\G\backslash\bH$. Indeed, since $\G_\fc\subset\G$, this map projects the rectangular set in \eqref{equ:idenset} onto $\cC_{Y,Y'}^{\G,\fc}$. To show that it is also injective, suppose $\Gamma \tau_\fc z=\Gamma \tau_\fc z'$ for some $z, z'$ from this rectangular set. Then there exists some $\gamma\in\G$ such that $\tau_\fc^{-1}\gamma\tau_\fc z=z'$. If $\gamma\in\pm \G_\fc$ then $\tau_\fc^{-1}\gamma\tau_\fc\in\pm \langle u_{\omega_\fc}\rangle$, and this implies that $z=z'$. Otherwise, let $\tau_\fc^{-1}\gamma\tau_\fc=\left(\begin{smallmatrix} a& b\\ c&d \end{smallmatrix}\right)\in \G_1$. Since $\gamma\not\in\pm \G_\fc$, $c\neq0$. We easily see this cannot happen since it would imply 
$$
\Im(z') = \frac{\Im(z)}{|cz+d|^2}=\frac{\Im(z)}{(cx+d)^2 +c^2y^2}\leq \frac{1}{y}\leq 1,
$$
contradicting that $\Im(z')>Y>1$.
%Let $\G z\in \cC_{Y,Y'}^{\Gamma,\fc}$, and set $z_1=\tau_\fc^{-1}\gamma z$, where $\gamma$ is chosen such that $\Im(z_1)=\cY_{\Gamma,\fc}$. We may replace $\gamma$ by $\gamma_1\gamma$, setting $z'=\tau_\fc^{-1}\gamma_1\gamma z$, with $\gamma_1\in\G_\fc$ such that $0\leq \Re(z')<\omega_\fc$. Then moreover $\Im(z')=\Im(z_1)\in(Y,Y')$. Conversely, given $z'$, it suffices to show that $\Im(\tau_\fc^{-1}\gamma \tau_\fc z')>Y$ for some $\gamma\in\Gamma$. Choosing $\gamma\in\Gamma_\fc$, we have that $\gamma$ is a power of $u_{\omega_\fc}$, and $\Im(\tau_\fc^{-1}\gamma \tau_\fc z')=\Im(z')\in(Y,Y').$ 
For the area computation, we use the definition \eqref{equ:norhyarea} of $\mu_{\G_1}$ together with $\mu_{\G_1}=[\G_1:\G]\mu_\G$ (since $-I_2\in \G$).
\end{proof}

\begin{Lem}\label{lem:disjoint}
Given two distinct cusps $\fc_1$, $\fc_2\in\Omega_\Gamma$, and any $Y_1$, $Y_2\geq1$, $\cC_{Y_1}^{\G,\fc_1}\cap\cC_{Y_2}^{\G,\fc_2}=\emptyset.$
\end{Lem}
\begin{proof}
Since $Y_1, Y_2\geq 1$, the sets $\{\tau_{\fc_1} z\in\bH: \Im(z)>Y_1\}$ and $\{\tau_{\fc_2} z\in\bH:\Im(z)>Y_2\}$ are subsets of the interior of the Ford circles based at $\fc_1$ and $\fc_2$ respectively. Two Ford circles are either disjoint or identical. Suppose $\G z\in \cC_{Y_1}^{\G,\fc_1}\cap\cC_{Y_2}^{\G,\fc_2}$. Then there exists an isometry $\gamma\in\G$ that maps the Ford circle at $\fc_1$ to the Ford circle at $\fc_2$. Consequently, we must have $\gamma \fc_1=\fc_2$, which is a contradiction.
\end{proof}

\begin{rmk}\label{rmk:disopen}
We will later consider sets $I_{y,Y,\fc}:=\left\{x\in (0,1):\G(x+iy)\in \cC_Y^{\G,\fc}\right\}$ for some $y>0, Y>1$ and $\fc\in\Omega_{\G}$. This set is the intersection of the line segment $\{x+iy\in\bH:0<x<1\}$ with the preimage of $\cC_Y^{\G,\fc}$ in $\bH$ (under the natural projection from $\bH$ to $\G\bk \bH$).  %Namely, the parts of the horizontal horocycle $\{x+iy\in\bH:0< x<1\}$ visiting the cusp neighborhood $\cC_Y^{\G,\fc}$. 
By definition the preimage of $\cC_Y^{\G,\fc}$ is the disjoint (since $Y>1$) union of the infinitely many horodiscs $\left\{\tau_{\fc'}z\in\bH:\Im(z)>Y\right\}=H^{\circ}_{p/q,1/(2q^2Y)}$ for all cusps $\fc'=p/q\in \G\fc$. %. Moreover, for any $s\in\Q\cup \{\infty\}$ there exists some $\gamma\in \G$ and $\fc\in\Omega_{\G}$ such that $\gamma\tau_{\fc}\infty=s$, implying that $\gamma \tau_{\fc}\left\{z\in\bH:\Im(z)>Y\right\}=H^{\circ}_{p/q,1/(2q^2Y)}$ if $s=p/q$ primitive and $\gamma \tau_{\fc}\left\{z\in\bH:\Im(z)>Y\right\}=\left\{z\in\bH:\Im(z)>Y\right\}$ if $s=\infty$.
Moreover, note that a necessary condition for such a horodisc intersecting the line segment $\{x+iy\in\bH:0<x<1\}$ is that $p/q\in \G\fc\cap (-\tfrac{1}{2Y}, 1+\tfrac{1}{2Y})$ and $1/(q^2Y)>y$, i.e. $q^2<1/(yY)$. Thus there are only finitely many such horodiscs intersecting $\{x+iy\in\bH:0<x<1\}$. Moreover, each such intersection is an open interval and the set $I_{y,Y,\fc}\subset (0,1)$ is thus the disjoint union of these finitely many open intervals. %which in turn can be expressed as a disjoint union of finitely many open intervals. 
Similarly, for any $Y'>Y>1$ the set $\left\{x\in (0,1):\G(x+iy)\in \cC_{Y,Y'}^{\G,\fc}\right\}=I_{y,Y,\fc}\setminus \overline{I}_{y,Y',\fc}$ is also a disjoint union of finitely many open intervals.
\end{rmk}

\subsection{Left regular action of normalizing elements}
Let $h\in\SL_2(\Q)$ be a group element normalizing $\G$. The action of $h$ on $\Q\cup \{\infty\}$ (by M\"{o}bius transformation) induces a well-defined action on $\G\bk \left(\Q\cup \{\infty\}\right)$, the set of cusps of $\G$. %Throughout the remaining of the section, we will fix $\Omega_{\G}\subset \left(\Q\cup \{\infty\}\right)$ such that $\Omega_{\G}$ is stable under the action of $h$, that is, $h\fc\in\Omega_{\G}$ for any $\fc\in\Omega_{\G}$.
\begin{Lem}\label{lem:hec1}
For each $\fc\in \Omega_{\G}$, we have
\begin{equation}\label{equ:conj1}
h\G_{\fc} h^{-1}=\G_{h \fc}
\end{equation}
and
\begin{equation}\label{equ:conj2}
\tau_{h \fc}^{-1}h\tau_{\fc}=\begin{pmatrix}
\sqrt{\omega_{h\fc}/\omega_{\fc}} & *\\
0 & \sqrt{\omega_{\fc}/\omega_{h\fc}}\end{pmatrix}\in \SL_2(\Q).
\end{equation}
\end{Lem}
\begin{proof}
Since $h$ normalizes $\G$ we have $h\G_{\fc} h^{-1}=h\tau_{\fc}N\tau_{\fc}^{-1}h^{-1}\cap \G$. Thus to prove \eqref{equ:conj1} it suffices to show $h\tau_{\fc}N\tau_{\fc}^{-1}h^{-1}=\tau_{h \fc}N\tau_{h \fc}^{-1}$.
We show that $\tau_{h \fc}^{-1}h\tau_{\fc}$ is an upper triangular matrix. Indeed, $\tau_{h \fc}^{-1}h\tau_{\fc}\infty=\tau_{h\fc}^{-1}\left(h\fc\right)=\infty$. %, or equivalently, to show $\bm{e}_2\tau_{\fc}u_{j/n}\tau_{u_{j/n}\cdot \fc}^{-1}=\lambda \bm{e}_2$ for some $\lambda\neq 0$. 
%For this recall that $h\cdot \fc\in \Z^2_{\rm pr}$ is the unique primitive integral vector such that $ \fc h=\lambda h\cdot \fc$ for some $\lambda>0$. Hence we have
%$$\bm{e}_2\tau_{\fc}h\tau_{h\cdot \fc}^{-1}=\fc h\tau_{h\cdot \fc}^{-1}=\lambda (h\cdot \fc) \tau_{h\cdot \fc}^{-1}=\lambda \bm{e}_2,$$
%implying that
This proves \eqref{equ:conj1}. We moreover conclude that
\begin{equation}\label{equ:conj3}
\tau_{h \fc}^{-1}h\tau_{\fc}=\begin{pmatrix}
\lambda & *\\
0 & \lambda^{-1}\end{pmatrix}
\end{equation}
for some $\lambda\neq 0$, and it remains to show that $\lambda^2=\omega_{h\fc}/\omega_{\fc}$.
%To prove \eqref{equ:conj2} it remains to show $\lambda=\sqrt{\omega_{h\cdot \fc}/\omega_{\fc}}$. 
For this we conjugate the subgroup $\tau_{h\fc}^{-1}\G_{h\fc}\tau_{h\cdot \fc}$
%=\left\langle u_{\omega_{u_{1/n}\cdot \fc}}\right\rangle$ 
by the matrix $\tau_{h \fc}^{-1}h\tau_{\fc}$. We obtain with \eqref{equ:conj1} that
$$\tau_{\fc}^{-1}h^{-1}\tau_{h\fc}\left(\tau_{h\fc}^{-1}\G_{h\fc}\tau_{h\fc}\right)\tau_{h\fc}^{-1}h\tau_{\fc}=\tau_{\fc}^{-1}\G_{\fc}\tau_{\fc}=\left\langle u_{\omega_{\fc}}\right\rangle.$$
On the other hand, using \eqref{equ:conj3} and $\tau_{h\fc}^{-1}\G_{h \fc}\tau_{h\fc}=\left\langle u_{\omega_{h\fc}}
%\left(\begin{smallmatrix}
%1 & \omega_{h \fc}\\
%0 & 1\end{smallmatrix}\right)
\right\rangle$,
we have
$$\tau_{\fc}^{-1}h^{-1}\tau_{h \fc}\left(\tau_{h\fc}^{-1}\G_{h \fc}\tau_{h\fc}\right)\tau_{h\fc}^{-1}h\tau_{\fc}=
\left(\begin{smallmatrix} \lambda^{-1} &*\\0&\lambda\end{smallmatrix}\right)   
\left\langle \left(\begin{smallmatrix}
1 & \omega_{h\fc}\\
0 & 1\end{smallmatrix}\right)\right\rangle
\left(\begin{smallmatrix} \lambda &*\\0&\lambda^{-1}\end{smallmatrix}\right)
%\tau_{\fc}^{-1}h^{-1}\tau_{h \fc}\left\langle \left(\begin{smallmatrix}
%1 & \omega_{h\fc}\\
%0 & 1\end{smallmatrix}\right)\right\rangle\tau_{h\fc}^{-1}h\tau_{\fc}
=\left\langle \left(\begin{smallmatrix}
1 & \omega_{h\fc}/\lambda^2\\
0 & 1\end{smallmatrix}\right)\right\rangle.$$
Comparing both equations we conclude that $\lambda^2=\omega_{h\fc}/\omega_{\fc}$. Finally replacing $\tau_{h\fc}$ with $-\tau_{h\fc}$ if necessary, we can ensure $\lambda$ is positive. %implying that $\lambda=\sqrt{\omega_{h\cdot \fc}/\omega_{\fc}}$ (since $\lambda>0$). 
%Hence it remains to show that $\lambda>0$. For this recall that $u_{j/n}\cdot\fc$ is the unique primitive integer vector satisfying that $\fc u_{j/n}=r u_{j/n}\cdot\fc$ with $r>0$. By direct computation we get that the bottom right entry of the matrix $\tau_{\fc}u_{j/n}\tau_{u_{j/n}\cdot \fc}^{-1}$, i.e. $\lambda$ by \eqref{equ:conj3}, equals $r$ which is positive. This finishes the proof.
%this would follow from \eqref{equ:conj2} since an upper triangular matrix normalizes the unipotent subgroup $N$. Now to prove \eqref{equ:conj2} 
\end{proof}

\begin{Prop}\label{prop:leftcuspnbhds}
Let $Y'>Y>0$ and $\fc\in \Omega_{\G}$. If $\G z\in \cC^{\G,\fc}_{\omega_{\fc}Y, \omega_{\fc}Y'}$, then $\G hz\in \cC^{\G,h\fc}_{\omega_{h\fc}Y, \omega_{h\fc}Y'}$. Similarly, if $\G z\in \cC^{\G,\fc}_{\omega_{\fc}Y}$, then $\G hz\in \cC^{\G,h\fc}_{\omega_{h\fc}Y}$.
%\begin{equation}\label{equ:leftcuspnbhds}
%L_{h^{-1}}\chi_{\cC_{\omega_{\fc}Y,\omega_{\fc}Y'}^{\G,\fc}}=\chi_{\cC_{\omega_{h\fc}Y,\omega_{h\fc}Y'}}^{\G,h\fc},
%\end{equation}
%where $\chi_{\cC}$ denotes the indicator function of $\cC\subset \G\bk\bH$. Similarly for any $Y>0$ and for any $\fc\in \Omega_{\G}$ we have
%\begin{equation}\label{equ:leftcuspnbhdsdiff}
%L_{h^{-1}}\chi_{\cC_{\omega_{\fc}Y}^{\G,\fc}}=\chi_{\cC_{\omega_{h\fc}Y}^{\G,h\fc}}.
%\end{equation}
\end{Prop}
\begin{proof}
%First we note that it suffices to show that $$L_{h}\chi_{\cC_Y^{\G,\fc}}\leq \chi_{\cC_{\lambda_{h,\fc}Y}^{\G,h\fc}}$ since then applying the left regular $h$-action on $\chi_{\cC_{\lambda_{h,\fc}Y}^{\G,h\fc}}$ (noting that $h^{-1}$ also normalizes $\G$) we can get
The second statement follows from the first one by taking $Y'\to\infty$. %First we note that this left $h$-action is well-defined: For any $\gamma\in \G$, we have $h\gamma z=h\gamma h^{-1} hz$. Since $h$ normalizes $\G$, this implies that $\G hz=\G h\gamma z$ for any $\gamma\in \G$. 
Since $\G z\in \cC^{\G,\fc}_{\omega_{\fc}Y, \omega_{\fc}Y'}$, by definition there exists $z'=x'+iy'\in \bH$ with $0\leq x'<\omega_\fc$ and $\omega_{\fc}Y<y'< \omega_{\fc}Y'$ and $\G z=\G\tau_{\fc}z'$. Consider $h\tau_\fc z'=\tau_{h\fc}z''$ with $z''=\tau_{h\fc}^{-1}h\tau_{\fc} z'$.
%We note that by our assumption $Y>\max\{1,\omega_{\fc}/\omega_{h\fc}\}$, both $Y$ and $\lambda_{h,\fc}Y$ are greater than one. Hence we can use the geometric description \eqref{equ:setdes} for the two sets in \eqref{equ:leftcuspnbhds} (noting that these two sets are differences of two cusp neighborhoods which can be both described as in \eqref{equ:setdes}). Namely, we can identify $\cC_{Y,Y'}^{\G,\fc}$ with the set
%$$B_1:=\left\{\tau_{\fc} z\in\bH :0\leq \Re(z)< \omega_{\fc}, Y<\Im(z)<Y'\right\},$$
%and identify $\cC_{\lambda_{h,\fc}{Y,\lambda_{h,\fc}Y'}}^{\G,h\fc}$ with the set 
%$$B_2:=\left\{\tau_{h\fc} z\in\bH :t\leq \Re(z)< \omega_{h\fc}+t, \lambda_{h,\fc}Y<\Im(z)<\lambda_{h,\fc}Y'\right\}$$
%for some $t\in\R$ to be determined.
%It thus suffices to show that the map sending $\tau_{\fc}z\in\bH$ to $h\tau_{\fc} z\in \bH$ induces a bijection from $B_1$ to $B_2$ for some $t\in\R$. 
By \eqref{equ:conj2}, we have $\Im(z'') = (\omega_{h\fc}/\omega_\fc)\Im(z')\in (\omega_{h\fc}Y, \omega_{h\fc}Y')$, implying that $\Gamma hz=\G h\tau_\fc z' \in \cC_{\omega_{h\fc}Y,\omega_{h\fc}Y'}^{\Gamma,h\fc}$. % Since $0\leq x''\leq \omega_{h\fc}$ and $\omega_{h\fc}Y<y''<\omega_{h\fc}Y''$, we conclude that $\Gamma hz=\G h\tau_\fc z' \in \cC_{\omega_{h\fc}Y,\omega_{h\fc}Y'}^{\Gamma,h\fc}$.
\end{proof}

\section{Negative results: horocycles expanding arbitrarily fast}\label{sec:neg2}

%\begin{rmk}\label{rmk:subprime}
%
%\end{rmk}
In this section using the results from the previous section, we prove \thmref{thm:negative2} and \thmref{thm:compactsupport} which provide new limiting measures for the sequences $\left\{\delta_{n,x,y_n}\right\}_{n\in\N}$ and $\left\{\delta^{\rm pr}_{n,x,y_n}\right\}_{n\in\N}$, allowing $\{y_n\}_{n\in\N}$ to decay arbitrarily fast. %The main ingredient of the proof is the study of the left regular $u_{1/n}$-action on the set of cusp neighborhoods on the congruence cover $\G_n\bk \bH$, which in turn relies on the analysis from the previous section. 
For any $n\in\N$ we consider the congruence subgroup $\G_n< \SL_2(\Z)$ given by
\begin{equation}\label{equ:latdes}
\G_n:=\left\{\begin{pmatrix}
a & b\\
c & d\end{pmatrix}\in \SL_2(\Z):n^2\mid c,\ a \equiv d\equiv \pm 1 \Mod{n}\right\}.
\end{equation}
It is clear that $\G_1=\SL_2(\Z)$ and that $\G_n$ contains the congruence subgroup 
$$\G_1(n^2):=\left\{\gamma\in \SL_2(\Z):\gamma\equiv \begin{pmatrix}
1 & * \\
0 & 1\end{pmatrix} \Mod{n^2}\right\}.$$  

%As in the previous section, we also denote by $\G_1=\SL_2(\Z)$ throughout this section.
%As mentioned before we note that $\G_1=\SL_2(\Z)$.
%The key step in constructing the counterexamples in \secref{sec:cou1} is \lemref{lem:counter1} which relates a Diophantine inequality to the cusp excursion of the sample points \eqref{equ:samplepoints}. Here we will prove \thmref{thm:negative2} using a similar observation with the role of the Diophantine inequality replaced by cusp excursions on certain congruence covers, see \lemref{lem:cuspanalysis2}. %\begin{equation}\label{equ:latdes}
%$$\G_n:=\bigcap_{1\leq j\leq n}u^{-1}_{j/n}\G u_{j/n}.$$
%\end{equation}
%We note that it is clear from the definition that $u_{j/n}$ normalizes $\G_n$ for any integer $j$. 
%For any $n\in\N$ we define the subgroup $\G_n< \SL_2(\Z)$ by
%\begin{equation}\label{equ:latdes}
%\G_n:=\left\{\begin{pmatrix}
%a & b\\
%c & d\end{pmatrix}\in \SL_2(\Z):n^2 | c, a \equiv d\equiv \pm 1 \Mod{n}\right\}.
%\end{equation}
%It is easy to check that $\G_n$ is a congruence subgroup of $\SL_2(\Z)$ containing the congruence subgroup 
%$$\G_1(n^2):=\left\{\gamma\in \SL_2(\Z):\gamma\equiv \begin{pmatrix}
%1 & \star \\
%0 & 1\end{pmatrix} \Mod{n^2}\right\}.$$ 

\subsection{Basic properties of the congruence subgroups $\G_n$}
%In this subsection we prove some basic properties regarding the congruence subgroup $\G_n$ which will will be needed for our proof. 
%In this subsection we prove some preliminary results about the congruence subgroup $\G_n$. %Although eventually we will only work on the case when $n=p$ is a prime number, here we prove these results for general $n$. %Recall that the lattice $\G_n$ is defined as
%Let 
%$$\mathbb{P}(\R)=\left\{[v]:v\in \R^2\setminus \{0\}, [v_1]=[v_2] \textrm{if there exists $\lambda\neq 0$ such that $v_1=\lambda v_2$}\right\}$$%\{(1,x):x\in \R\}\cup \{(0,1)\}$$ 
%be the real projective line where $G$ acts naturally via the right multiplication. Recall that an element $[v]\in \mathbb{P}(\R)$ is a \textit{cusp} of $\G_n$ if the stabilizer group
%$$\left\{\gamma\in \G_n\ | \ [v]\gamma =[v]\right\}$$
%is nontrivial.
First we show that $\G_n$ is normalized by $u_{j/n}$ for any $j\in \Z$. As mentioned in the introduction this simple fact is the starting point of our proofs to \thmref{thm:negative2} and \thmref{thm:compactsupport}. 

\begin{Lem}\label{lem:explicitdes}
For any $n\in\N$ and for any $j\in\Z$, the unipotent matrix $u_{j/n}$ normalizes $\G_n$.
%we have
%\begin{equation}\label{equ:latdes}
%\G_n=\left\{\begin{pmatrix}
%a & b\\
%c & d\end{pmatrix}\in \G:n^2 | c, a \equiv d \Mod{n}\right\}.
%\end{equation}
\end{Lem}
\begin{proof}
%For any $n\in\N$ we denote by 
%$$\G_n':=\left\{\begin{pmatrix}
%a & b\\
%c & d\end{pmatrix}\in \G:n^2 | c, a \equiv d \Mod{n}\right\},$$
%and we want to show $\G_n=\G_n'$. First we note that when $n=1$ this is true since $\G_1=\G$ and the two congruence conditions in the definition of $\G_n'$ are void, implying that also $\G_1'=\G$. Now we assume $n\geq 2$. By definition we have $\gamma\in \G_n$ if and only if $u_{j/n}\gamma u_{j/n}^{-1}\in\G$ for any $1\leq j\leq n$. 
By direct computation, for any $\gamma=\left(\begin{smallmatrix} 
a & b\\
c& d\end{smallmatrix}\right)\in \G_1$ and for any $j\in\Z$ we have
%\begin{equation}\label{equ:exdes2}
$$u^{-1}_{j/n}\gamma u_{j/n}=\begin{pmatrix}
a-\frac{jc}{n} & b+\frac{(a-d)j}{n}-\frac{j^2c}{n^2}\\
c & d+\frac{jc}{n}\end{pmatrix}.$$
%\end{equation}
Hence if $\gamma\in \G_n$, that is, $n^2\mid c$ and $a\equiv d\equiv \pm 1 \Mod{n}$, all the entries are integers with the bottom left entry divisible by $n^2$, and  
$$a-\frac{jc}{n}\equiv a\equiv d\equiv d+\frac{jc}{n}\equiv \pm 1\ \Mod{n}.$$
This implies that $u^{-1}_{j/n}\G_n u_{j/n}\subset \G_n$. %For the other inclusion, we conjugate $\G_n$ by $u_{-j/n}$ to get $u_{-j/n}\G_n u_{j/n}\subset \G_n$, or equivalently, $\G_n\subset u^{-1}_{j/n}\G_n u_{j/n}$. This finishes the proof.
%it is then clear from \eqref{equ:exdes2} that all the entries of $u_{j/n}\gamma u_{j/n}^{-1}$ are integers, implying that $u_{j/n}\gamma u_{j/n}^{-1}\in \G$ for any $1\leq j\leq n$, or equivalently, $\gamma\in \G_n$. This proves that $\G_n'\subset \G_n$. For the other inclusion, we want to show take any $\gamma=\left(\begin{smallmatrix}
%a & b\\
%c& d\end{smallmatrix}\right)\in \G_n$ and we want to show $n^2 | c$ and $a\equiv d \Mod{n}$. Since $u_{j/n}\gamma u_{j/n}^{-1}\in \G$ for any $1\leq j\leq n$, we have $a+\frac{jc}{n}\in \Z$ for $j=1,2$. Subtracting these two equations we get $n | c$. Next, we also have $b+\frac{(d-a)j}{n}-\frac{j^2c}{n^2}\in\Z$ for $j=1, n-1$. Subtracting these two equations and using the facts that $1^2\equiv (n-1)^2 \Mod{n}$ and $n | c$ we get 
\end{proof}

Next we prove the following index formula for $\G_n$.
\begin{Lem}\label{lem:indexformula}
For any integer $n\geq 3$, we have 
\begin{equation}\label{equ:index}
[\G_1: \G_n]=\frac{n^3}{2}\prod_{\substack{p | n\\ \textrm{prime}}}\left(1-p^{-2}\right).
\end{equation}
\end{Lem}
\begin{proof}
Let $J_n< \left(\Z/n^2\Z\right)^{\times}$ be the subgroup 
\begin{equation}\label{equ:hgh}
J_n:= \left\{[a]\in\left(\Z/n^2\Z\right)^{\times}:a\equiv \pm 1 \Mod{n}\right\}.
\end{equation}
It is easy to check that $\#(J_n)=2n$. Consider the map $h: \G_n \to J_n$ sending $\gamma=\left(\begin{smallmatrix}
a & b\\
c & d\end{smallmatrix}\right)\in \G_n$ to $[a]\in \left(\Z/n^2\Z\right)^{\times}$. Using the definition of $\G_n$, one can check that $h$ is a group homomorphism with the kernel $\ker(h)=\G_1(n^2)$.  %the image of $h$ is
%\begin{equation}\label{equ:image}
%h\left(\G_n\right)=\left\{[d]\in \left(\Z/n^2\Z\right)^{\times}:d^2\equiv 1 \Mod{p}\right\}.
%\end{equation}
%To prove the claim,
% first for any $\gamma=\left(\begin{smallmatrix}
%a & b\\
%c & d\end{smallmatrix}\right)\in \widetilde{\G}(p)$, since $a\equiv d \Mod{p}$ and $ad=1+bc\equiv 1 \Mod{p}$ we have $d^2\equiv ad\equiv 1 \Mod{p}$. This proves that the image $h\left(\widetilde{\G}(p)\right)$ is contained in the right hand side of \eqref{equ:image}. For the other direction, we note that since $p$ is a prime number, $d^2\equiv 1 \Mod{p}$ is equivalent to $d\equiv \pm 1\Mod{p}$. Thus the right hand side of \eqref{equ:image} equals
%\begin{equation}\label{equ:alternative}
%\left\{[kp\pm 1]\in  \left(\Z/p^2\Z\right)^{\times}:0\leq k\leq p-1\right\}.
%\end{equation}
%Moreover, we note that since $p$ is an odd prime number, this set has $2p$ elements. Now 
For each $0\leq k\leq n-1$, set $\gamma_k^{\pm}=\pm\left(\begin{smallmatrix}
1+kn & 1\\
-k^2n^2 & 1-kn\end{smallmatrix}\right)\in \G_n$. Then $h$ surjects the set $\left\{\gamma_k^{\pm}\in \G_n:0\leq k\leq n-1\right\}$ onto $J_n$. Finally we use the index formula for $\G_1(n^2)$ (see e.g. \cite[Section 1.2]{DiamondShurman2005}) to get
\begin{displaymath}
[\G_1 : \G_n]=\frac{[\G_1 : \G_1(n^2)]}{[\G_n : \G_1(n^2)]}=\frac{[\G_1 : \G_1(n^2)]}{\# J_n}=\frac{n^3}{2}\prod_{\substack{p | n\\ \textrm{prime}}}\left(1-p^{-2}\right).\qedhere
\end{displaymath}
\end{proof}
%\begin{rmk}\label{rmk:index2}
%It is easy to check that $\G_2=\G_0(4)$ and we have $[\G_1: \G_2]=6$.
%\end{rmk}
%\begin{rmk}\label{rmk:cosetrep}
%For later use we note here that from the above proof it is clear that for $n\geq 3$, $h$ induces an identification between the set $\left\{\gamma_k^{\pm}\in \G_n:0\leq k\leq n-1\right\}$ and the image $h$ which is isomorphic to $\ker(h)\bk \G_n=\G_1(n^2)\bk \G_n$. Hence the set $\left\{\gamma_k^{\pm}\in \G_n:0\leq k\leq p-1\right\}$ forms a complete list of  coset representatives for $\G_1(n^2)\bk \G_n$. Similarly, when $n=2$, the set $\left\{I_2, \left(\begin{smallmatrix} -1 & 1\\ -4 & 3\end{smallmatrix}\right)\right\}$ forms a complete list of  coset representatives for $\G_1(4)\bk \widetilde{\G}(2)$.
%\end{rmk}

Next, we study the properties of $\G_n$ relative to its cusps. As in \secref{sec:leftregular} we denote by $\Omega_{\G_n}$ the set of cusps of $\G_n$. The following lemma computes the width of each cusp of $\G_n$.
\begin{Lem}\label{lem:hec3}
Let $n\in\N$ and let $\fc=m/l\in \Omega_{\G_n}$ with $\gcd(m,l)=1$ $($if $\fc=\infty$, $m/l$ is understood as $1/0$$)$. Then we have
$$\omega_{\fc}=\frac{n^2}{\gcd (n,l)^2}.$$
\end{Lem}
\begin{proof}
Let $\tau_{\fc}\in \G_1$ be as before such that $\tau_{\fc}\infty=\fc$. Thus the left column of $\tau_{\fc}$ is $\left(\begin{smallmatrix}
m \\
l\end{smallmatrix}\right)$. By direct computation we have 
$$\tau_{\fc}N\tau_{\fc}^{-1}=\left\{\begin{pmatrix}
1-ml t& m^2t\\
-l^2t & 1+mlt\end{pmatrix}\in G:t\in \R\right\}.$$
Thus by \eqref{equ:latdes} an element in $(\G_n)_{\fc}=\tau_{\fc}N\tau_{\fc}^{-1}\cap \G_n$ is of the form $\gamma=\left(\begin{smallmatrix}
1-ml t& m^2t\\
-l^2t & 1+mlt\end{smallmatrix}\right)\in \G_1$ satisfying that $n^2\mid l^2 t$ and $1- mlt\equiv 1+mlt\equiv \pm 1\Mod{n}$. %We claim that the first condition implies the second one. 
%For this we prove that $n^2\mid l^2 t$ implies that $n\mid mlt$ (which would then imply that $1\pm mlt\equiv 1\Mod{n}$). 
Looking at the top right and bottom left entries of $\gamma$, 
we have that $m^2 t, l^2 t\in \Z$. Since $\gcd(m,l)=1$, we have $t\in\Z$. Then the condition $n^2\mid l^2t$ is equivalent to $\frac{n^2}{\gcd(n,l)^2}\mid t$, and the condition $n\mid mlt$ is equivalent to that $\frac{n}{\gcd(n,ml)}\mid t$. Moreover, since $\frac{n}{\gcd(n,ml)}\mid \frac{n^2}{\gcd(n,l)^2}$, the condition $\frac{n}{\gcd(n,ml)}\mid t$ is implied by the condition $\frac{n^2}{\gcd(n,l)^2}\mid t$. We conclude that $n^2\mid l^2 t$ implies $1-mlt\equiv 1+mlt\equiv \pm 1\Mod{n}$. Thus
$$(\G_n)_{\fc}= \left\{\begin{pmatrix}
1-ml t& m^2t\\
-l^2t & 1+ mlt\end{pmatrix}\in\G_1: n^2\mid l^2t\right\}.$$
Conjugating $(\G_n)_{\fc}$ back via $\tau_{\fc}$ and using the equivalence of the two conditions $n^2 \mid l^2 t$ and $\frac{n^2}{\gcd(n,l)^2} \mid t$ we get
$$\tau_{\fc}^{-1}(\G_n)_{\fc}\tau_{\fc}=\left\{u_t=\begin{pmatrix}
1 & t\\
0 & 1\end{pmatrix}\in \G_1:\frac{n^2}{\gcd(n,l)^2}\mid t\right\},$$
%First we show that $t$ is an integer. For this we note since $\gcd(m,l)=1$, there exist $k_1, k_2\in\Z$ such that $k_1m+k_2l=1$. Then we have
%$$t=(k_1m+k_2l)^2t=k_1^2m^2t+2k_1k_2mlt+k_2^2l^2t\in \Z,$$
%which follows from the condition that $m^2t, mlt$ and $l^2t$ are integers. Next 
%Now for $u_t\in \tau_{\fc}\G_n_{\fc}\tau_{\fc}^{-1}$, since $t\in\Z$  Hence we have
%$$\tau_{\fc}\G_n_{\fc}\tau_{\fc}^{-1}=\left\{\begin{pmatrix}
%1 & t\\
%0 & 1\end{pmatrix}\in \G\ \bigg|\ \frac{n^2}{\gcd(n,m)^2}\;|\; t\right\},$$
implying that $\omega_{\fc}=n^2/\gcd(n,l)^2$.
\end{proof}
%\subsubsection{Explicit description of cusps at prime levels}
%While it seems quite difficult for us to describe cusps of $\G_n$ for general $n\in\N$, 
Next we compute the number of cusps of $\G_n$. %The next proposition counts the number of cusps of $\G_n$ for $n\geq 3$.
\begin{Prop}\label{prop:numbercusps}
For any integer $n\geq 3$ we have
$$\#\Omega_{\G_n}=\frac{n^2}{2}\prod_{\substack{p | n\\ \textrm{prime}}}\left(1-p^{-2}\right).$$
\end{Prop}
\begin{rmk}\label{rmk:2cusp}
It is easy to check that $\G_2=\G_0(4)$. Thus $[\G_1: \G_2]=6$ and $\G_2$ has three cusps which can be represented by $\infty$, $1/2$ and $1$ respectively. 
\end{rmk}
%\begin{Lem}\label{lem:cuspdescription}
%Let $(m,l), (m',l')\in \Z^2_{\rm pr}$ be two primitive integral vectors. Then $(m,l)\sim_{\G_n} (m',l')$ if and only if there exists some integers $a\in\Z$ satisfying that $a\equiv 1\Mod{n}$ and $b\in \Z$ such that
%$$(m',l')\equiv \pm (am, \overline{a}l+bm) \Mod{n^2},$$
%%$$m'\equiv am \Mod{n^2}\quad \textrm{and}\quad l'\equiv \overline{a}l \Mod{d},$$
%where $\overline{a}\in\Z$ is the multiplicative inverse of $a$ modulo $n^2$.
%\end{Lem}
To prove \propref{prop:numbercusps} we first prove a preliminary formula for $\#\Omega_{\G_n}$.
\begin{Lem}\label{lem:precuspformula}
For any integer $n\geq 3$ we have
%\begin{equation}\label{equ:precuspfor}
$$\#\Omega_{\G_n}=\sum_{d | n^2}\frac{\varphi(n^2/d)\varphi(d)\gcd(n^2/d,d)}{2n}.$$
%\end{equation}
\end{Lem}
\begin{proof}
Since $-I_2\in \G_n$ and $\G_1(n^2) < \G_n$, we have $\Omega_{\G_n}=\G_n\bk \Omega_{\G_1(n^2)}$. On the other hand, by the analysis in \cite[p. 102]{DiamondShurman2005}, the set $\Omega_{\G_1(n^2)}$ is in bijection with the union of cosets
$\bigsqcup_{d| n^2}\langle \pm I_2\rangle\bk Z_d$, where for each $d\mid n^2$,
$$Z_d:=\left\{([m],[l])^t:[m]\in \left(\Z/d\Z\right)^{\times}, [l]\in \Z/n^2\Z,  \gcd(n^2, l)=d\right\}$$
with $([m],[l])^t$ is the transpose of the row vector $([m],[l])$ and the bijection is induced by the map sending $m/l\in \Q\cup \{\infty\}$ with $\gcd(m,l)=1$ to $([m],[l]))^t\in Z_{d}$ with $d=\gcd(n^2,l)$. Note that $\# Z_d=\varphi(n^2/d)\varphi(d)$.

For each $d\mid n^2$, using the definition of $\G_n$, it is easy to check that the linear action of $\G_n$ on $\Z^2$ (by matrix multiplication) induces a well-defined action of $\G_n$ on $Z_d$ and that the corresponding action of the subgroup $\G_1(n^2)$ is trivial. From the proof of \lemref{lem:indexformula}, we have $\G_n/\G_1(n^2)\cong J_n$, where
\begin{equation}\label{equ:deshn}
J_n=\left\{\pm [1+kn]\in (\Z/n^2\Z)^{\times}:0\leq k\leq n-1\right\},
\end{equation}
which is of size $2n$. Hence the action of $\G_n$ on $Z_d$ induces the action of $J_n$ on $Z_d$ given by 
$$[a]\cdot ([m], [l])^t=([am], [\overline{a}l])^t,$$ 
with $([m],[l])^t\in Z_d$ and $\overline{a}$ the multiplicative inverse of $a$ modulo $n^2$. We note that $[am]\in \left(\Z/d\Z\right)^{\times}$ is well-defined since $d\mid n^2$.

We conclude that $\Omega_{\G_n}=\G_n\bk \Omega_{\G_1(n^2)}$ is in bijection with the union of cosets
$$\bigsqcup_{d| n^2}\G_n\bk Z_d=\bigsqcup_{d| n^2}J_n\bk Z_d,$$
implying that 
$$\#\Omega_{\G_n}=\sum_{d| n^2}\# J_n\bk Z_d.$$
Hence we want to compute the size of the coset $J_n\bk Z_d$ for each $d\mid n^2$. For this we claim that for any for any $([m] ,[l])^t\in Z_d$, the orbit $J_n\cdot ([m],[l])^t$ is of size $2n/\gcd(n^2/d,d)$, implying that 
$$\#J_n\bk Z_d=\frac{\#Z_d}{2n/\gcd(n^2/d,d)}=\frac{\varphi(n^2/d)\varphi(d)\gcd(n^2/d,d)}{2n}.$$ 
We note that \lemref{lem:precuspformula} then follows immediately from this claim. To prove this claim, it suffices to compute the size of the stabilizer 
$$(J_n)_{([m],[l])}:=\left\{[a]\in J_n:[a]\cdot ([m],[l])^t= ([m], [l])^t\in Z_d\right\}.$$
Since by definition $[a]\cdot ([m],[l])^t=([am], [\overline{a}l])^t$, $[a]\in (J_n)_{([m],[l])}$ if and only if $am\equiv m \Mod{d}$ and $\overline{a}l\equiv l\Mod{n^2}$. Since $d=\gcd(n^2,l)$ and $[m]\in \left(\Z/d\Z\right)^{\times}$, these two conditions are equivalent to $a\equiv 1 \Mod{d}$ and $\overline{a}\equiv 1\Mod{n^2/d}$, which are equivalent to $a\equiv 1\Mod{\lcm(n^2/d, d)}$. Hence using the description \eqref{equ:deshn} of $J_n$ and the facts that $n\mid \lcm(n^2/d, d)$ and $\lcm(n^2/d, d)\gcd(n^2/d, d)=n^2$ we have
$$(J_n)_{([m],[l])}=\left\{[1+\lcm(n^2/d, d) j]\in J_n:0\leq j\leq \gcd(n^2/d, d)-1\right\}$$
is of size $\gcd(n^2/d, d)$. This implies that 
$$\#\left(J_n\cdot ([m],[l])^t\right)=\frac{\#J_n}{\#(J_n)_{([m],[l])}}=\frac{2n}{\gcd(n^2/d,d)},$$
proving the claim, and hence also this lemma. 
%\begin{align*}
%\#\left(\Omega_{\G_n}\right)&=\sum_{d | n^2}\frac{\varphi(n^2/d)\varphi(d)\gcd(n^2/d,d)}{2n}.
%\end{align*}
\end{proof}
We can now give the proof of \propref{prop:numbercusps} by simplifying the formula in \lemref{lem:precuspformula}.
\begin{proof}[Proof of \propref{prop:numbercusps}]
Write $n=\prod_{i=1}^k p_i^{\alpha_i}$ in the prime decomposition form and apply \lemref{lem:precuspformula} to get
\begin{align*}
\#\left(\Omega_{\G_n}\right)&=\frac{1}{2n}\sum_{\bm{\beta}\in \Z^k: 0\leq \beta_i\leq 2\alpha_i}\varphi\left(\prod_{i=1}^kp_i^{\beta_i}\right)\varphi\left(\prod_{i=1}^kp_i^{2\alpha_i-\beta_i}\right)\prod_{i=1}^kp_i^{\min\{\beta_i, 2\alpha_i-\beta_i\}},
\end{align*}
where the summation is over all vectors $\bm{\beta}=(\beta_1,\ldots, \beta_k)\in \Z^k$ satisfying $0\leq \beta_i\leq 2\alpha_i$ for all $1\leq i\leq k$, and we used that $\gcd(n^2/d, d)=\prod_{i=1}^kp_i^{{\min\{\beta_i, 2\alpha_i-\beta_i\}}}$ for $d=\prod_{i=1}^k p_i^{\beta_i}$. Using the fact that $\varphi$ is multiplicative and interchanging the summation and product signs we get
\begin{align*}
\#\left(\Omega_{\G_n}\right)&=\frac{1}{2n}\prod_{i=1}^k\left(\sum_{0\leq \beta_i\leq 2\alpha_i}\varphi\left(p_i^{\beta_i}\right)\varphi\left(p_i^{2\alpha_i-\beta_i}\right)p_i^{\min\{\beta_i, 2\alpha_i-\beta_i\}}\right)\\
&=\frac{1}{2n}\prod_{i=1}^k\left(\sum_{1\leq \beta_i\leq 2\alpha_i-1}p_i^{2\alpha_i}(1-p_i^{-1})^2p_i^{\min\{\beta_i, 2\alpha_i-\beta_i\}}+2p_i^{2\alpha_i}(1-p_i^{-1})\right)\\
&=\frac{1}{2n}\prod_{i=1}^kp_i^{2\alpha_i}(1-p_i^{-1})\left((1-p_i^{-1})\sum_{1\leq \beta_i\leq 2\alpha_i-1}p_i^{\min\{\beta_i, 2\alpha_i-\beta_i\}}+2\right),
\end{align*}
where for the second equality we used that for $1\leq \beta_i\leq 2\alpha_i-1$, $\varphi\left(p_i^{\beta_i}\right)\varphi\left(p_i^{2\alpha_i-\beta_i}\right)=p^{2\alpha_i}(1-p_i^{-1})^2$, and for $\beta_i=0$ or $\beta_i=2\alpha_i$, $\varphi\left(p_i^{\beta_i}\right)\varphi\left(p_i^{2\alpha_i-\beta_i}\right)=p^{2\alpha_i}(1-p_i^{-1})$ and $\min\{\beta_i, 2\alpha_i-\beta_i\}=0$. We note that the term $\sum_{1\leq \beta_i\leq 2\alpha_i-1}p_i^{\min\{\beta_i, 2\alpha_i-\beta_i\}}$ equals
\begin{align*}
&\sum_{1\leq \beta_i\leq \alpha_i}p_i^{\beta_i}+\sum_{\alpha_i< \beta_i\leq 2\alpha_i-1}p_i^{2\alpha_i-\beta_i}=\sum_{1\leq \beta_i\leq \alpha_i}p_i^{\beta_i}+\sum_{1\leq \beta_i< \alpha_i}p_i^{\beta_i}\\
&=2\sum_{1\leq \beta_i\leq \alpha_i}p_i^{\beta_i}-p_i^{\alpha_i}=\frac{2p_i(p_i^{\alpha_i}-1)}{p_i-1}-p_i^{\alpha_i}.
\end{align*}
Hence we have
\begin{align*}
\#\left(\Omega_{\G_n}\right)&=\frac{1}{2n}\prod_{i=1}^kp_i^{2\alpha_i}(1-p_i^{-1})\left((1-p_i^{-1})\left(\frac{2p_i(p_i^{\alpha_i}-1)}{p_i-1}-p_i^{\alpha_i}\right)+2\right)\\
&=\frac{1}{2n}\prod_{i=1}^kp_i^{2\alpha_i}(1-p_i^{-1})p_i^{\alpha_i}(1+p_i^{-1})=\frac{n^2}{2}\prod_{i=1}^k(1-p_i^{-2}),
\end{align*}
finishing the proof.
%Let $\cC_{\G_1(n^2)}$ be the set of cusps of $\G_1(n^2)$. By the analysis in \cite[]{DiamondShurman2005}, the set $\cC_{\G_1(n^2)}$ is in bijection with the coset
%$\left(\bigsqcup_{d| n^2}Z_d\right)/\langle \pm I_2\rangle,$%=\bigsqcup_{d| n^2}\left(Z_d/\langle \pm I_2\rangle\right),$
%where the action of $-I_2$ on $\bigsqcup_{d| n^2}Z_d$ is given by the right multiplication and the bijection is given by identifying $(m,l)\in \Omega_{\G_n}$ with $\pm ([m], [l])\in \left(\Z/n^2\Z\right)\times \left(\Z/d\Z\right)^{\times}$ with $d=\gcd(m,n^2)$. Since by \rmkref{rmk:cosetrep} the set $\left\{\gamma_k^{\pm}\in \G_n:0\leq k\leq n-1\right\}$ forms a complete list of coset representatives for $\G_1(n^2)\bk \G_n$ and $-I_2$ is contained in this set, we have the set of cusp of $\G_n$ is in bijection with the coset
%$$\left(\bigsqcup_{d| n^2}Z_d\right)/\left\{\gamma_k^{\pm}\in \G_n:0\leq k\leq n-1\right\}=\left(\bigsqcup_{d| n^2}Z_d\right)/A_n=\bigsqcup_{d| n^2}\left(Z_d/A_n\right),$$
%where the action of $\gamma_k^{\pm}$ is given by right multiplication and the action of $A_n$ is defined such that $a\cdot ([m], [l]):=([am], [\overline{a}l])$ for $a\in A_n$ and $([m],[l]\in \bigsqcup_{d| n^2}Z_d$.
\end{proof}
\subsection{Proof of \thmref{thm:negative2}}
For simplicity of notation, we abbreviate the cusp neighborhoods $\cC_Y^{\G_n,\fc}$ and $\cC_{Y,Y'}^{\G_n,\fc}$ by $\cC_Y^{n,\fc}$ and $\cC_{Y,Y'}^{n,\fc}$ respectively and the set of cusps $\Omega_{\G_n}$ by $\Omega_n$. We first prove the following key lemma which says that if $\G_n z$ visits a cusp neighborhood on $\G_n\bk \bH$, then all companion points $\G_1 u_{j/n} z, 0\leq j\leq n-1$ make excursions to some cusp neighborhood on $\cM=\G_1\bk\bH$, the modular surface. We recall that $\cC_Y$ is the projection onto $\cM$ of the region
$\{z\in\bH:\Im(z)>Y\}$.

\begin{Lem}\label{lem:cuspanalysis2}
 Let $Y>0$ and $n\in \N$. If $\G_n z\in \cC_{\omega_{\fc}Y}^{n,\fc}$ for some $\fc\in \Omega_{n}$ then $\G_1 u_{j/n}z\in \cC_Y$ for all $0\leq j\leq n-1$. 
%and let $\fc\in\cC_{\widetilde{\G}(p)}$ be a cusp of type $\rom{1}$ or $\rom{2}$. If $\widetilde{\G}(p)g\in \cC_Y^{\widetilde{\G}(p), \fc}$ for some $g\in G$, then $\G u_{j/p}g\in\cC_Y$ for any $1\leq j\leq p$.
%\begin{enumerate}
%\item if $\fc\in \cC_{\rom{1}}$ is of type $\rom{1}$, then $\G u_{j/p}g\in \cC_Y$ for all $j\in \{1,2,\ldots, p\}$;
%\item if $\fc\in \cC_{\rom{2}}$ is of type $\rom{2}$, then $\G u_{j/p}g\in \cC_Y$ for all but one  $j\in \{1,2,\ldots, p\}$.
%\end{enumerate}  
\end{Lem}
\begin{proof}
%Let $\fc\in\Omega_{\G_n}$ be the cusp such that $\G_n z\in \cC_{\omega_{\fc}Y}^{n,\fc}$. Then by \rmkref{rmk:disjointuniopen} and the estimate $\omega_{\fc}\geq 1$ we have $\G_1 z\in \cC_{\omega_{\fc}Y}\subset \cC_Y$. It thus remains to treat the case when $j\neq 0$. 
Fix $0\leq j\leq n-1$. By \lemref{lem:explicitdes}, $u_{j/n}$ normalizes $\G_n$. Assuming that $\G_n z\in \cC_{\omega_{\fc}Y}^{n,\fc}$ and applying \propref{prop:leftcuspnbhds} to $h=u_{j/n}$, we get $\G_n u_{j/n}z\in \cC_{\omega_{h\fc}Y}^{n,h\fc}$. By definition, there exists $z'\in\bH$ with $\Im(z')>\omega_{h\fc}Y\geq Y$ such that $\Gamma_n \tau_{h\fc} z' = \Gamma_n u_{j/n} z$. Since $\tau_{h\fc}\in\Gamma_1$, this implies $\Gamma_1 u_{j/n}z = \G_1 z'\in \cC_Y$.
%Again by \rmkref{rmk:disjointuniopen} and using the estimate $\omega_{h\fc}\geq 1$ we get $\G_1 u_{j/n}z\in \cC_{Y}$, finishing the proof.
%Let $Y'=\omega_{\fc}Y$ such that $\lambda_{h,\fc}Y'=\omega_{h\fc}Y$. Hence both $Y'$ and $\lambda_{h,\fc}Y'$ are greater or equal to $Y>1$, and we can apply \eqref{equ:leftcuspnbhdsdiff} to left regular $h$-action on the cusp neighborhood $\cC_{Y'}^{n,\fc}=\cC_{\omega_{\fc}Y}^{n,\fc}$ to get $\G_n u_{j/n}z=\G_n hz\in \cC_{\lambda_{h,\fc}Y'}^{n,h\fc}=\cC_{\omega_{h\fc}Y}^{n,h\fc}$. Similarly, by \rmkref{rmk:disjointuniopen} and the estimate $\omega_{h\fc}Y\geq Y>1$ we have $\G_1 u_{j/n} z\in \cC_{\omega_{h\fc}Y}\subset \cC_Y$, finishing the proof.
 \end{proof}

We can now give the 
\begin{proof}[Proof of \thmref{thm:negative2}]
For any $n\in \N$ let $Y_n=\max\{\log n, 1\}$,
%Let $\{Y_p\}_{p\in\mathbb{P}}$ be an non-decreasing unbounded sequence of positive numbers such that $Y_p\geq 1$ and the series 
%\begin{equation}\label{equ:diverser}
%\sum_{p\in\mathbb{P}}\frac{1}{(p+1)Y_p}=\infty.
%\end{equation}
%For example, for any $k\in \N$ let $p_k\in \mathbb{P}$ be the $k$-th odd prime number. Then using the estimate $p_k\ll k\log k$ we can take $Y_3=1$ and $Y_{p_k}=\max\left\{\log\log(k),1\right\}$. 
and let $\Psi_n$ %=\sum_{\fc\in \cC_{\widetilde{\G}(p)}^{\rom{1}}\cup\cC_{\widetilde{\G}(p)}^{\rom{2}}}\chi_{\widetilde{C}_{Y_p}^{\fc}}$ 
be the indicator function of the union 
$$\bigcup_{\fc\in \Omega_n}\cC_{\omega_{\fc}Y_n,2\omega_{\fc}Y_n}^{n,\fc}\subset \G_n\bk \bH.$$
%$\widetilde{E}_{Y_n}^n\subset \G_n\bk G$ defined as in \rmkref{}.
%of the cusp neighborhoods $\widetilde{C}_{Y_p}^{\widetilde{\G}(p), \fc}$ for all the cusps of $\widetilde{\G}(p)$ of type $\rom{1}$ and $\rom{2}$. 
Since for any cusp $\fc\in \Omega_{n}$, $\omega_{\fc}Y_n\geq Y_n\geq 1$, {{by \lemref{lem:onetoone}, each $\cC_{\omega_{\fc}Y_n,2\omega_{\fc}Y_n}^{n,\fc}$ is a Borel set with boundary of measure zero; and by \lemref{lem:disjoint} the above union is disjoint. Thus $\Psi_n$ is the indicator function of a Borel set with boundary of measure zero. Moreover,}} applying the volume formula  \eqref{equ:volume}, the index formula in \lemref{lem:indexformula} and the cusp number formula in \propref{prop:numbercusps} (see also \rmkref{rmk:2cusp} for the case when $n=2$) we have for any $n\in \N$,
\begin{equation}\label{equ:volucom}
\mu_{\G_n}\left(\Psi_n\right)=\sum_{\fc\in \Omega_n}\mu_{\G_n}\left(\cC_{\omega_{\fc}Y_n,2\omega_{\fc}Y_n}^{n,\fc}\right)=\sum_{\fc\in \Omega_n}\frac{3\omega_{\fc}}{\pi [\G_1 : \G_n]}\times \frac{1}{2\omega_{\fc}Y_n}=\frac{3}{2\pi Y_n}\frac{\# \Omega_n}{[\G_1 : \G_n]}\asymp\frac{1}{nY_n}.
\end{equation} 
%Moreover, by \propref{prop:primecusp}, \rmkref{rmk:uaction} and the volume formula in \lemref{lem:cuspnbhdvolume} we have
%\begin{equation}
%\mu_{\widetilde{\G}(p)}(\Psi_p)=\sum_{\fc\in \cC_{\widetilde{\G}(p)}^{\rom{1}}\cup\cC_{\widetilde{\G}(p)}^{\rom{2}}}\frac{6\omega_{\fc}}{\pi p(p^2-1)Y_p}=\frac{3}{\pi(p+1)Y_p}.
%%\frac{2}{p(p^2-1)}\left(\sum_{\fc\in \cC_{\rom{1}}\cup\cC_{\rom{2}}}\omega_{\fc}\right)\mu_{\G}(\widetilde{C}_{Y_p})=\frac{\mu_{\G}(\widetilde{C}_{Y_p})}{p-1}
%\end{equation}
For any $n\in\N$ and $0<y<1$ we define 
$$I_n(y):=\left\{x\in \R/\Z:\Psi_n(x+iy)=1\right\}.$$ 
%We first note that for any fixed integer $n\in \N$, $\cI_n(y)\subset [0,1)$ is a finite disjoint union of open intervals. More precisely, it is the intersection of the horizontal segment 
%$$\{z\in\mathbb{H}:0\leq \Re(z) <1\}$$ 
%(naturally identified with the unit interval $[0,1)$) with the horodiscs $\mathcal{H}^{\circ}_{m/n,1/(2Y_pn^2)}\subset \mathbb{H}$, where $m/n$ runs through the finite set of primitive rational numbers in $[0,1)$ satisfying that $p| n$ and $1/Y_pn^2> y$. 
By definition, $x\in I_n(y)$ if and only if $\G_n(x+iy)\in \cC_{\omega_{\fc}Y_n,2\omega_{\fc}Y_n}^{n,\fc}\subset \cC_{\omega_{\fc}Y_n}^{n,\fc}$ for some $\fc\in \Omega_n$. Thus \lemref{lem:cuspanalysis2} implies that
$$
I_n(y) \subset \{x\in\R/\Z: \cR_n(x,y)\subset \cC_{Y_n}\}.
$$
% we get $\cR_n(x, y)\subset \cC_{Y_n}$ for any $x\in I_n(y)$.
This, together with our choice that $Y_n=\max\{\log n, 1\}$ and the distance formula \eqref{equ:distance}, implies that for any $n\geq 3$ and for any $x\in I_n(y)$  
$$\inf_{\G_1 z\in \cR_n(x, y)}d_{\cM}(\G_1 z_0,\G_1 z)\geq \log (Y_n)+O(1)=\log\log n+ O(1).$$
It thus suffices to show that there exists a sequence $\{y_n\}_{n\in \N}$ satisfying that $0<y_n<c_n$ for all $n\in\N$ and that the limsup set $\limsup_{n\to\infty}I_n(y_n)$ is of full Lebesgue measure in $\R/\Z$. 

For this, we will construct a sequence $\{y_n\}_{n\in\N}$ decaying sufficiently fast and then apply the quantitative Borel-Cantelli lemma \corref{cor:quanbc} to the sequence $\{I_n(y_n)\}_{n\in\N}\subset \R/\Z$.  To ensure the quasi-independence condition \eqref{equ:quasicond2} in \corref{cor:quanbc}, we need, for every pair $1\leq m<n\in\N$, the two quantities $\left|I_m(y_m)\cap I_n(y_n)\right|$ and $\left|I_m(y_m)\right|\left|I_n(y_n)\right|$ to be sufficiently close to each other. The key observations for this are the following two relations that
\begin{equation}\label{equ:keyob1}
\left|I_n(y_n)\right|=\int_0^1\Psi_n(x+iy_n)dx
\end{equation}
and
\begin{equation}\label{equ:keyob2}
\left|I_{m}(y_m)\cap I_n(y_n)\right|=\int_0^1\Psi_{m}(x+iy_m)\Psi_n(x+iy_n) dx=\int_{I_{m}(y_m)}\Psi_n(x+iy_n) dx.
\end{equation}
Assuming the limit equation \eqref{equ:equipie} holds for the pairs $((0,1), \Psi_n)$ and $(I_m(y_m), \Psi_n)$ (we will verify this later), then by relation \eqref{equ:keyob2} the quantity $\left|I_{m}(y_m)\cap I_n(y_n)\right|$ is close to the quantity $|I_m(y_m)|\mu_{\G_n}(\Psi_n)$ which in turn is close to $|I_m(y_m)||I_n(y_n)|$ by relation \eqref{equ:keyob1}, provided that $y_n>0$ is sufficiently small.

We now implement the above ideas rigorously. We first claim that there exists a sequence $\{y_n\}_{n\in\N}$ satisfying, for all $n\in\N$, $0<y_n<c_n$ and
\begin{equation}\label{equ:conscond}
\left|\frac{1}{\left|I\right|}\int_{I}\Psi_n(x+iy_n)dx-\mu_{\G_n}(\Psi_n)\right|\leq \frac{\mu_{\G_n}(\Psi_n)
}{2n^2},
\end{equation}
for any subset $I\subset \R/\Z$ taken from the finite set
\begin{equation}\label{equ:finiteset}
\left\{(0,1)\right\}\bigcup\left\{I_{m}(y_{m}): 1\leq m<n\right\}.
\end{equation}
%$$\left| \left|I_{p}(y_p)\cap \cI_{p'}(y_{p'}\right|-\left|I_{p}(y_p}\right|\left|\cI_{p'}(y_{p'}\right|\right| \leq \left|I_{p}(y_p}\right|\left|\cI_{p'}(y_{p'}\right|$$
%Before proving the claim, we first note that since $\omega_{\fc}Y_n\geq 1$ for any $n\in\N$ and for any $\fc\in \Omega_n$, in view of the one-to-one correspondence in \lemref{lem:onetoone}, one can easily construct a sequence of compactly supported and continuous functions $\{\Psi_{n,j}^{\pm}\}_{j\in\N}$ on $\G_n\bk \bH$ such that $\Psi_{n,j}^-\leq \Psi_n\leq \Psi_{n,j}^+$ for every $j\in\N$ and $\lim\limits_{j\to\infty}\mu_{\G_n}\left(\Psi_{n,j}^{\pm}\right)=\mu_{\G_n}(\Psi_n)$. 
{{For this, first note that}} by {{\rmkref{rmk:apparg}}} for any $I\subset \R/\Z\cong [0,1)$ a disjoint union of finitely many open intervals, we have
\begin{equation}\label{equ:equdistest}
\lim\limits_{y\to 0^+}\frac{1}{|I|}\int_{I}\Psi_n(x+iy)dx= \mu_{\G_n}(\Psi_n).
\end{equation}
%We may approximate the indicator function $\Psi_n$ by two positive smooth functions $\Psi^\pm_n$ in its neighborhood such that $\Psi_n^-\leq \Psi_n\leq \Psi_n^+$. Consequently $\mu_{\Gamma_n}(\Psi_{n}^-)\leq \mu_{\Gamma_n}(\Psi_n)\leq \mu_{\Gamma_n}(\Psi_n^+).$
%%%%%%We first note that it is clear from the definition of the cusp neighborhoods that for any $n\in\N$, the boundary of $\supp \Psi_n$ is of measure zero. 
We now construct such a sequence successively. For the base case $n=1$ since \eqref{equ:equdistest} %(see also \rmkref{rmk:indica}) 
holds for the pair $((0,1),\Psi_1)$ on $\cM=\G_1\bk\bH$, %we have
%$$\int_0^1\Psi_1^\pm(x+iy)dx\ \xrightarrow{y\to 0^+}\ \mu_{\G_1}(\Psi_1^\pm).$$
there exists $0<y_1<c_1$ sufficiently small such that 
$$\left|\int_0^1\Psi_1(x+iy_1)dx- \mu_{\G_1}(\Psi_1)\right|< \frac12\mu_{\G_1}(\Psi_1).$$
For a general integer $n\geq 2$, suppose that we already have chosen $0<y_m<c_m$ satisfying \eqref{equ:conscond} for all the positive integers $m< n$. By \rmkref{rmk:disopen} the set $I_m(y_m)\subset \R/\Z$ is a disjoint union of finitely many open intervals for any $m<n$. Thus \eqref{equ:equdistest} is satisfied for all the pairs 
$$\left((0,1), \Psi_n\right), (I_m(y_m), \Psi_n), 1\leq m< n$$
on $\G_n\bk\bH$. Since there are only finitely many such pairs, we can take $0<y_n<c_n$ sufficiently small such that \eqref{equ:conscond} is satisfied for all $I\in \left\{(0,1)\right\}\bigcup\left\{I_{m}(y_{m}): 1\leq m<n\right\}$, which is the set in \eqref{equ:finiteset}. This finishes the proof of the claim. 
%Let $I\subset \R/\Z$ be an arbitrary disjoint union of open intervals in $\R/\Z\cong [0,1)$. Thus again in view of \rmkref{rmk:indica} we can apply \eqref{equ:equipie} to the pair $(I, \Psi^\pm_n)$ on $\G_n\bk \bH$ to get
%$$\frac{1}{\left|I\right|}\int_{I}\Psi^\pm_n(x+iy)dx\ \xrightarrow{y\to 0^+}\  \mu_{\G_n}(\Psi^\pm_n).$$
%In particular, in view of \rmkref{rmk:disopen} we can take $0<y_n<c_n$ sufficiently small such that \eqref{equ:conscond} is satisfied for all $I$ coming from the finite set in \eqref{equ:finiteset}. This finishes the proof of the claim. 

Now let $\{y_n\}_{n\in\N}$ be as in the claim. For any $n\in\N$ apply \eqref{equ:conscond} to the pair $((0,1),\Psi_n)$ we get
%$$\left|I_n(y_n)\right|=\int_0^1\Psi_n(x+iy_n)dx.$$
%Moreover, by \eqref{equ:conscond},
\begin{equation}\label{equ:ine2}
\left| \left|I_n(y_n)\right|-\mu_{\G_n}(\Psi_n)\right| \leq\frac{\mu_{\G_n}(\Psi_n)}{2n^2}.%\leq \frac12\mu_{\widetilde{\G}(p)}(\Psi_p).
\end{equation}
By the triangle inequality, this implies 
\begin{equation}\label{equ:boundmu}
\mu_{\Gamma_n}(\Psi_n)\leq\ 2|I_n(y_n)|.
\end{equation}
More generally, for each $1\leq m< n$
%$$\left|I_{m}(y_m)\cap I_n(y_n)\right|=\int_0^1\Psi_{m}(x+iy_m)\Psi_n(x+iy_n) dx=\int_{I_{m}(y_m)}\Psi_n(x+iy_n) dx,$$
apply \eqref{equ:conscond} to the pair $(I_{m}(y_{m}),\Psi_n)$ we get
\begin{equation}\label{equ:ine1}
\left|\left| I_{m}(y_{m})\cap I_n(y_n)\right|-\left| I_{m}(y_{m})\right|\mu_{\G_n}(\Psi_n)\right|\leq \frac{\left| I_{m}(y_{m})\right|\mu_{\G_n}(\Psi_n)}{2n^2}.
\end{equation}
Using the inequalities \eqref{equ:ine2}, \eqref{equ:boundmu}, \eqref{equ:ine1} together with the triangle inequality we get
\begin{equation}\label{equ:qbc2}
\left|\left| I_{m}(y_{m})\cap I_n(y_n)\right|-\left| I_{m}(y_{m})\right|\left|I_n(y_n)\right|\right|\leq \frac{\left|I_{m}(y_{m})\right|\mu_{\G_n}(\Psi_n)}{n^2}\leq \frac{2\left| I_{m}(y_{m})\right|\left|I_n(y_n)\right|}{n^2}.
\end{equation}
Hence the sequence $\left\{I_n(y_n)\right\}_{n\in\N}\subset \R/\Z$ satisfies the quasi-independence condition \eqref{equ:quasicond2} (with the subset $\mathbb{S}=\N$ and the exponent $\eta=2$). %In particular, by \rmkref{rmk:furcond} it satisfies the quasi-independence condition \eqref{equ:qbccond}.
%We now want to apply \lemref{lem:bcd} (the quantitative Borel-Cantelli lemma) to the sequence $\{\cI_n(y_n)\}_{n\in\N}$. For this we let $(X,\cB,\mu)=([0,1),\cB,|\cdot|)$ with $\cB$ the usual Borel $\sigma$-algebra on $[0,1)$ generated by open intervals. For any $i\in\N$ we let $A_i=\cI_p(y_p)$ with $p\in\mathbb{P}$ the $i$-th odd prime number. For any $i,j\in\N$, let $R_{i,j}=\mu(A_i\cap A_j)-\mu(A_i)\mu(A_j)$ be as in \lemref{lem:bcd}, and we note that \eqref{equ:qbc2} implies that for any $1\leq i<j$, $R_{i,j}\leq 2\mu(A_i)\mu(A_j)/j^2$. Hence in view of \rmkref{rmk:furcond} the sequence $\{A_i\}_{i\in\N}$ satisfies the condition \eqref{equ:qbccond}. 
%We want to check that $\{A_i\}_{i\in\N}$ satisfies the quasi-independence condition \eqref{equ:qbccond}. First using the estimate $R_{i,i}=\mu(A_i)-\mu(A_i)^2\leq \mu(A_i)$ and the symmetry $R_{i,j}=R_{j,i}$ for any $i,j\in\N$, we have for any $k_2> k_1\geq 1$
%\begin{align*}
%\sum_{i,j=k_1}^{k_2}R_{i,j}&=\sum_{i=k_1}^{k_2}R_{i,i}+\sum_{k_1\leq i< j\leq k_2}R_{i,j}+\sum_{k_1\leq j<i\leq k_2}R_{i,j}\leq \sum_{i=k_1}^{k_2}\mu(A_i)+2\sum_{k_1\leq i< j\leq k_2}R_{i,j}.
%%&\leq \sum_{i=k_1}^{k_2}\mu(A_i)+2\sum_{k_1\leq i< j\leq k_2}\mu(A_i)\mu(A_j)/j^2
%\end{align*}
%Hence it suffices to show
%$$\sum_{k_1\leq i< j\leq k_2}R_{i,j}\ll \sum_{i=k_1}^{k_2}\mu(A_i).$$
%the sequence $\{\cI_p(y_p)\}_{p\in\mathbb{P}}$ satisfies the \textit{quais-independence} condition that for any pair of prime numbers $p'<p$,
%$$\left| \cI_{p'}(y_{p'})\cap \cI_p(y_p)\right|\leq \frac52 \left| \cI_{p'}(y_{p'})\right|\left|\cI_p(y_p)\right|.$$
Moreover, using the inequality \eqref{equ:ine2}, the volume computation \eqref{equ:volucom} and the estimate that $Y_n\asymp \log n$ we have that
$$\sum_{n\in \N}\left| I_n(y_n)\right|\geq \sum_{n\in\N}\frac12\mu_{\G_n}(\Psi_n)\asymp \sum_{n\in\N}\frac{1}{n\log n}=\infty.$$
%diverges which follows from our choice of the sequence $\{Y_p\}_{p\in \mathbb{P}}$. 
Thus by \corref{cor:quanbc}, $\limsup_{n\to\infty}I_n(y_n)\subset \R/\Z$ is of full Lebesgue measure, finishing the proof.
\end{proof}

\begin{rmk}\label{rmk:fastestrate}
It is not clear to us whether the rate $\log\log n$ is the fastest excursion rate for generic translates. We note that in principle it can be proved (or disproved) if one can compute the volume of the set
$$\cE_Y^n:=\left\{\G_n z\in \G_n\bk \bH:\G_1 u_{j/n}z\in \cC_Y\ \textrm{for all $0 \leq j\leq n-1$}\right\}.$$
For instance, if one can show $\mu_{\G_n}(\cE_Y^n)\asymp 1/(nY)$ for all $n\in\N$ and for all $Y\geq 1$, then \thmref{thm:negative2} together with a standard application of the Borel-Cantelli lemma would imply that the inequality in \eqref{equ:loglaw} is indeed an equality for almost every $x\in \R/\Z$. We also note that our analysis (\lemref{lem:disjoint} and \lemref{lem:cuspanalysis2}) shows that for any $n\in\N$ and for any $Y\geq 1$
$$\bigsqcup_{\fc\in \Omega_n}\cC_{\omega_{\fc}Y}^{n,\fc}\subset \cE_Y^n\subset \bigsqcup_{\fc\in \Omega_n}\cC_{Y}^{n,\fc},$$
implying that $1/(nY)\ll \mu_{\G_n}\left(\cE_Y^n\right)\ll 1/Y $. On the other hand using some elementary arguments (which relies on the width computation \lemref{lem:hec3}) one can show that any $\langle u_{1/n}\rangle$-orbit contains at least one cusp of width one. %for any cusp $\fc\in \Omega_{\G_n}$ there exists some $0\leq j\leq n-1$ such that $\omega_{u_{j/n}\cdot \fc}=1$. 
This fact together with the fact that $1\leq \omega_{\fc}\leq n^2$ implies that $\cE_Y^n= \bigsqcup_{\fc\in \Omega_n}\cC_{Y}^{n,\fc}$ when $Y\geq n^2$ . However, both estimates are not sufficient for the purpose of obtaining an upper bound. %implying that $\mu_{\G_n}\left(\cE_Y^n\right)\asymp 1/(nY)$ for $Y\geq n$.
\end{rmk}

\begin{rmk}\label{rmk:andreas}
Here we give a very brief sketch of the argument communicated to us by Str\"ombergsson: For each $n\in\N$ and $y>0$, it is not difficult to see that $\G_n(x+iy)\in \mathcal{C}_{\omega_{\fc}Y_n}^{n,\fc}$ for some $\fc=\tfrac{p}{q}\in\Omega_n$ with $\gcd(p,q)=1$ if and only if
\begin{equation}\label{equ:equinum}
\left|x-\frac{p}{q}\right|^2<\frac{y}{\omega_{\fc}Y_nq^2}-y^2=\frac{y\gcd(n,q)^2}{n^2Y_n q^2}-y^2.
\end{equation}
Here $Y_n=\max\{\log n, 1\}$ is as in the above proof. Define 
$$\tilde{I}_n(y):=\left\{x\in\R/\Z:\ \textrm{$\exists$ primitive $\frac{p}{q}$ s.t. $n\mid q$, $q<\frac{1}{2\sqrt{yY_n}}$, $\left|x-\frac{p}{q}\right|<\frac{\sqrt{y}}{2\sqrt{Y_n }q}$}\right\}.$$
One can easily check that elements in $\tilde{I}_n(y)$ satisfy the inequality \eqref{equ:equinum}. Hence by \lemref{lem:cuspanalysis2} we have 
\begin{equation}\label{equ:relationtow}
\tilde{I}_n(y)\subset  \{x\in\R/\Z: \cR_n(x,y)\subset \cC_{Y_n}\}.
\end{equation} 
Moreover, using some standard techniques from analytic number theory one can show that for any subinterval $I\subset \R/\Z$ (or more generally, any finite disjoint union of subintervals), 
$$\lim\limits_{y\to 0^+}|I|^{-1}\left|\tilde{I}_n(y)\cap I\right|=\frac{c_n}{Y_n}$$
with $c_n=\frac{3}{\pi^2}\frac{\varphi(n)}{n^2}\prod_{p\nmid n}(1-p^{-2})^{-1}\gg \frac{\varphi(n)}{n^2}$. This limit equation is the analog of \eqref{equ:equdistest}. Another input is the divergence of the series $\sum_{n\in\N}\frac{c_n}{Y_n}\gg \sum_{n\in\N}\frac{\varphi(n)}{n^2\log n}$, which follows from the estimate $\varphi(n)\gg n/\log\log n$.
With these two inputs one can then mimic the arguments in the above proof to construct a sequence $\{y_n\}_{n\in\N}$ decaying sufficiently fast and then apply \corref{cor:quanbc} to get a full measure limsup set $\limsup_{n\to\infty}\tilde{I}_n(y_n)\subset \R/\Z$. Finally, we note that the relation \eqref{equ:relationtow} can be checked directly using the definition of the set $\tilde{I}_n(y)$. Hence this argument can be carried over without going into the congruence covers $\G_n\bk\bH$.
\end{rmk}

\subsection{Proof of \thmref{thm:compactsupport}}
We prove \thmref{thm:compactsupport} in this subsection. The strategy is similar to that of \thmref{thm:negative2} with the sequence of cuspidal sets approaching the cusps replaced by a sequence of compact cylinders approaching certain closed horocycles. Let $n\in\N$ be an integer and let $\G_n z\in \G_n\bk\bH$ be a point close to a cusp $\fc\in\Omega_n$. For any $0\leq j\leq n-1$, the analysis in \secref{sec:leftregular} gives exact information about the height of the companion point $\G_n u_{j/n}z$ with respect to the cusp $u_{j/n}\fc$. While this is sufficient for \thmref{thm:negative2} (cusp excursions), to realize the limiting measure $\nu_{m,Y}$ in \thmref{thm:compactsupport} one needs more refined information about the spacing of these companion points along the closed horocycles they lie on. For this, we further analyze the left regular $u_{1/n}$-action on points near certain type cusps which we now define. 
%We first define a family of cusps that we will further analyze. Let $n\in\N$ be an integer. %While it seems quite difficult to understand $\langle u_{1/n}\rangle$-orbit for the $u_{1/n}$-action on $\Omega_{\G_n}$ completely, it is relatively easy to understand such orbits for certain cusps. 

We say $\fc\in\Omega_n$ is of \textit{simple type} if $\fc$ can be represented by a primitive rational number $m/q$ satisfying that $\gcd(n^2,q)\mid n$, and we denote by $\Omega_n^{\rm sim}\subset \Omega_n$ the set of simple type cusps. {{(This notion of simple type cusps is closely related to the condition $n\in\N_q$ in \thmref{thm:fullrange}. In fact, let $p/q$ be a primitive rational number then the condition $n\in \N_q$ is equivalent to that the cusp $\fc\in\Omega_n$ represented by $p/q$ is of simple type.)}} If $m'/q'$ is another representative for $\fc$, that is, $m'/q'$ is primitive and $m'/q'=\gamma (m/q)$ for some $\gamma\in\G_n$, then using the definition of $\G_n$, it is easy to check that $\gcd(n^2, q)=\gcd(n^2, q')$. Hence the simple type cusps are well-defined.

As mentioned in \secref{sec:fullrange} the condition $\gcd(n^2,q)\mid q$ implies the further decomposition $q=kl$ with $l=\gcd(n,q)\mid n$ and $k=q/l$ satisfying $\gcd(k,n)=1$. We can thus reparameterize a simple type $\fc$ by $m/(kl)$ with $\gcd(m,kl)=\gcd(k,n)=1$ and $l\mid n$. The main new ingredient of our proof to \thmref{thm:compactsupport} is the following decomposition of the sample points which generalizes \eqref{equ:symreg}.
%which relies on more refined analysis of the left regular $u_{1/n}$-action on points near cusps of simple type. 
%We will give an explicit description of the $\langle u_{1/n}\rangle$-orbits for cusps of simple type, see \lemref{lem:symmetrysimpletype}. For now we first compute the the size of $\Omega_{\G_n}^s$.
%\begin{Lem}\label{lem:simpletype}
%For any integer $n\geq 3$ we have
%$$\#\Omega^{\rm sim}_{\G_n}=\frac{n^2}{2}\prod_{\substack{p | n\\ \textrm{prime}}}\left(1-p^{-1}\right)=\frac{n\varphi(n)}{2}.$$
%\end{Lem}
%\begin{proof}
%It is clear from the proof of \lemref{lem:precuspformula} that
%\begin{align*}
%\#\Omega^{\rm sim}_{\G_n}=\sum_{d\mid n}\frac{\varphi(n^2/d)\varphi(d)\gcd(n^2/d,d)}{2n}.
%\end{align*}
%Now using the identities $\varphi(n^2/d)=\varphi(n)\tfrac{n}{d}$ and $\gcd(n^2/d,d)=d$ for any $d\mid n$ and the identity $\sum_{d\mid n}\varphi(d)=n$ we get
%\begin{displaymath}
%\#\Omega^{\rm sim}_{\G_n}=\sum_{d\mid n}\frac{\varphi(n)\times \tfrac{n}{d}\times \varphi(d)\times d}{2n}=\frac{\varphi(n)}{2}\sum_{d\mid n}\varphi(d)=\frac{n\varphi(n)}{2}.\qedhere
%\end{displaymath}
%%Then the formula in this lemma follows from similar computations as in the above proof of \propref{prop:numbercusps}.
%\end{proof}

\begin{Prop}\label{prop:symmetrysimpletype}
Fix $n\in\N$, $z=x+iy\in \bH$ and $\fc\in \Omega_n^{\rm sim}$. Then
%Let $z'=x'+iy'\in \bH$ such that $\G_n z=\G_n \tau_{\fc}z'$ for some $\fc\in \Omega_n^{\rm sim}$ of simple type. Then we have
$$\cR_n(x,y)=\bigcup_{d\mid n}\cR_{n/d}^{\rm pr}(x'_{\fc, d}, d^2y'/\omega_{\fc}),$$
where $z'=x'+iy'\in\bH$ is such that $\G_n z=\G_n \tau_{\fc}z'$, and $x'_{d,\fc}\in \R/\Z$ depends only on $x'$, $\fc$ and $d$.
\end{Prop}
%\begin{rmk}\label{rmk:gener}
%\propref{prop:symmetryset} is a special case of \propref{prop:symmetrysimpletype}: Let $x=p/q$ and $n\in\N_q$ be as in \propref{prop:symmetryset}. Then the cusp $\fc\in\Omega_n$ represented by $p/q$ is a simple type cusp, and for any $y>0$ we have by direct computation the relation $\G_n(x+iy)=\G_n\tau_{\fc}(x'+iy')\in \G_n\bk\bH$ with $y'=\tfrac{1}{q^2y}$. Further using the relation $\omega_{\fc}=n^2/\gcd(n,q)^2$ one can see that \propref{prop:symmetrysimpletype} recovers \propref{prop:symmetryset} when $x=p/q$.
%\end{rmk}
We first prove a simple lemma computing the width of elements in the orbits $\langle u_{1/n}\rangle \fc$ when $\fc\in \Omega_n^{\rm sim}$ is of simple type. 
\begin{Lem}\label{lem:widthsimple}
Fix $n\in\N$ and $\fc\in \Omega_n^{\rm sim}$ a simple type cusp. Then for any $0\leq j\leq n-1$ we have 
$$\omega_{u_{j/n}\fc}=\gcd(m\tfrac{n}{l}+jk,n)^2,$$
where $m/(kl)$ is a representative for $\fc$ with $\gcd(m,kl)=\gcd(k,n)=1$ and $l\mid n$. 
\end{Lem}
\begin{proof}
For any $0\leq j\leq n-1$, 
$$u_{j/n}\fc=\frac{m}{kl}+\frac{j}{n}=\frac{m\tfrac{n}{l}+jk}{kn}=:\frac{p_j}{q_j}$$
with $\gcd(p_j,q_j)=1$. Let $d_{j}:=\gcd(m\tfrac{n}{l}+jk,kn)$ such that $q_j=kn/d_{j}$. Since $\gcd(mn,k)=1$, we have $\gcd(m\tfrac{n}{l}+jk,k)=\gcd(m\tfrac{n}{l},k)=1$. Hence $d_{j}=\gcd(m\tfrac{n}{l}+jk,n)\mid n$. Now by \lemref{lem:hec3} and the assumption that $\gcd(k,n)=1$ we have
\begin{displaymath}
\omega_{u_{j/n}\fc}=\frac{n^2}{\gcd(n,kn/d_{j})^2}=d_{j}^2=\gcd(m\tfrac{n}{l}+jk,n)^2.\qedhere
\end{displaymath}
\end{proof}
We can now combine ideas from \secref{sec:fullrange} and \secref{sec:leftregular} to give the
\begin{proof}[Proof of \propref{prop:symmetrysimpletype}]
Assume $\fc=m/(kl)$ with $\gcd(m,kl)=\gcd(k,n)=1$ and $l\mid n$. Up to changing the representatives for $\fc$, we may assume $mkl\neq 0$. Let $\tau_{\fc}=\left(\begin{smallmatrix}
m & a\\
kl & b\end{smallmatrix}\right)\in \G_1$, and for each $1\leq j\leq n-1$ let $\tau_{u_{j/n}\fc}=\left(\begin{smallmatrix}
p_j & v_j\\
q_j & w_j\end{smallmatrix}\right)\in \G_1$, where $p_j,q_j$ are as in the proof of \lemref{lem:widthsimple}, $a, b, v_j, w_j$ are some integers such that $\tau_{\fc}, \tau_{u_{j/n}\fc}\in\G_1$, that is, 
\begin{equation}\label{equ:tworela}
mb-kla=1\quad \textrm{and}\quad (m\tfrac{n}{l}+jk)w_j-knv_j=d_{j}
\end{equation}
with $d_j=\gcd(m\tfrac{n}{l}+jk,n)$ as in the proof of \lemref{lem:widthsimple}.
%We note that $a, b$ depend only on $\fc$, while $v,w$ also depend on $j$. 
By direct computation and using \lemref{lem:hec1} and \lemref{lem:widthsimple} (and the relation $\omega_{\fc}=d_0^2=n^2/l^2$) we have
$$\tau_{u_{j/n}\fc}^{-1}u_{j/n}\tau_{\fc}=\begin{pmatrix}
d_{j}l/n & w_ja+b(\tfrac{jw_j}{n}-v_j)\\
0 & n/(d_{j}l)\end{pmatrix}.$$
Using the relations in \eqref{equ:tworela} the top right entry becomes
\begin{align*}
w_ja+b\left(\frac{jw_j}{n}-v_j\right)&=w_ja+\frac{1+kla}{m}\left(\frac{jw_j}{n}-v_j\right)=\frac{a(w_jmn+jw_jkl-klv_jn)+jw_j-v_jn}{mn}\\
&=\frac{ad_{j}l}{mn}+\frac{1}{mn}\left(\frac{d_{j}l-w_jmn}{kl}\right)=\frac{bd_{j}}{nk}-\frac{w_j}{kl}.
\end{align*}
%where for the third equality we used the second equation in \eqref{equ:tworela} and for the last equality we used the first equation in \eqref{equ:tworela}. 
(Here we used the assumption that $mkl\neq 0$.) Hence we have for any $0\leq j\leq n-1$
\begin{align}\label{equ:relationshard}
\G_n u_{j/n}z&=\G_n u_{j/n}\tau_{\fc}(x'+iy')=\G_n\tau_{u_{j/n}\fc}\tau_{u_{j/n}\fc}^{-1}u_{j/n}\tau_{\fc}(x'+iy')\\
&=\G_n\tau_{u_{j/n}\fc}\left(\tfrac{d_{j}^2l^2}{n^2}x'+\tfrac{d^2_{j}lb}{n^2k}-\tfrac{d_{j}w_j}{kn}+i\tfrac{d_{j}^2l^2}{n^2}y'\right).\nonumber
\end{align}
Here for the first equality we used the assumption that $\G_n z=\G_n \tau_{\fc} z'$ and the fact that $u_{j/n}$ normalizes $\G_n$.
Now as in the proof of \propref{prop:symmetryset} for any $d\mid n$, we define 
$$D_d:=\{0\leq j\leq n-1: d_{j}=d\}$$
so that 
%Since $d_{j,\fc}=\gcd(m\tfrac{n}{l}+jk,n)$ and $\gcd(k,n)=1$, as $j$ runs through $\{0,1,\ldots, n-1\}$, $[m\tfrac{n}{l}+jk]\in \Z/n\Z$ runs though all the congruence classes (exactly once) in $\Z/n\Z$. This implies that 
%\begin{displaymath}
%\left\{[m\tfrac{n}{l}+jk]\in (\Z/n\Z)^{\times}: j\in D_d\right\}=\left\{[l]\in (\Z/n\Z)^{\times}: [l]\in \Z/n\Z,\ \gcd(l, n)=d\right\}.
%\end{displaymath}
%In particular, we have 
\begin{equation}\label{equ:predecompss}
\cR_n(x,y)=\bigcup_{d\mid n}\left\{\G_1 u_{j/n}z\in \cM: j\in D_d\right\},
\end{equation}
and
\begin{equation}\label{equ:congrucopr}
\left\{[(m\tfrac{n}{l}+jk)/d]\in (\Z/(n/d)\Z)^{\times}: j\in D_d\right\}=(\Z/(n/d)\Z)^{\times}.
\end{equation}
Use the second relation in \eqref{equ:tworela} to get for $j\in D_d$,
$$w_j\left((m\tfrac{n}{l}+jk)/d\right)\equiv 1\Mod{k\tfrac{n}{d}}.$$
Solving the above congruence equation as in the proof of \lemref{lem:symmetryrational} we get 
$$w_j\equiv dl\overline{mn}\tfrac{n}{d}e+\left(\left(m\tfrac{n}{l}+jk\right)/d\right)^*kf \Mod{k\tfrac{n}{d}},$$
where for any integer $t$, $\overline{t}$ denotes the multiplicative inverse modulo $k$, $t^*$ denotes the multiplicative inverse modulo $n/d$, and $e=e_d, f=f_d\in \Z$ are two fixed integers such that $e\tfrac{n}{d}+fk=1$. Plugging this relation into \eqref{equ:relationshard} and using the relation $\omega_{\fc}=n^2/l^2$ we get for any $d\mid n$ and for any $j\in D_d$,
$$\G_n u_{j/n}z= \G_n\tau_{u_{j/n}\fc}\left(x'_{\fc,d}-\tfrac{\left(\left(m\tfrac{n}{l}+jk\right)/d\right)^* f}{n/d}+i\tfrac{d^2y'}{\omega_{\fc}}\right),$$
where $x'_{\fc,d}:=\tfrac{d^2l^2}{n^2}x'+\tfrac{d^2lb}{n^2k}-\tfrac{dl\overline{mn}e}{k} \Mod{\Z}\in\R/\Z$. Since $\tau_{u_{j/n}\fc}\in \G_1$ we have
$$\left\{\G_1 u_{j/n}z\in \cM: j\in D_d\right\}=\left\{\G_1\left(x'_{\fc,d}-\tfrac{\left(\left(m\tfrac{n}{l}+jk\right)/d\right)^* f}{n/d}+i\tfrac{d^2y'}{\omega_{\fc}}\right)\in \cM: j\in D_d\right\}.$$
Thus in view of \eqref{equ:predecompss} and the above relation it suffices to show 
\begin{displaymath}
\left\{-[\left((m\tfrac{n}{l}+jk)/d\right)^*f]\in (\Z/(n/d)\Z)^{\times}: j\in D_d\right\}=(\Z/(n/d)\Z)^{\times}.
\end{displaymath}
But this follows from \eqref{equ:congrucopr} and the fact that $\gcd(f, \tfrac{n}{d})=1$ (since $\gcd(f, \tfrac{n}{d})=\gcd(fk,\tfrac{n}{d})=\gcd(1-e\tfrac{n}{d}, \tfrac{n}{d})=1$), and we have thus finished the proof.
\end{proof}
We will also need the following lemma estimating the number of cusps in $\Omega_n^{\rm sim}$ satisfying certain restrictions on the width.
\begin{Lem}\label{lem:simpletype}
Let $m\in \N$ be a fixed integer and let $n=m\ell\geq 3$ for some prime number $\ell$ not dividing $m$. Then we have
$$\#\left\{\fc\in \Omega^{\rm sim}_{n}: \omega_{\fc}\geq m^2\right\}\geq \frac{\varphi(m)(\ell-1)^2}{2}.$$
\end{Lem}
\begin{proof}
Recall from the proof of \lemref{lem:precuspformula} that $\Omega_n$ is in bijection with the disjoint union $\bigsqcup_{d\mid n^2} J_n\backslash Z_d$. On the other hand, by definition of the simple type cusps, $\Omega_n^{\rm sim}$ corresponds to the subset $\sqcup_{d\mid n} J_n\backslash Z_d$. Moreover, let $\fc=m/l\in\Omega_n^{\rm sim}$ with $\gcd(m,l)=1$ be a simple type cusp corresponding to an element in $J_n\bk Z_d$ for some $d\mid n$, that is, $d=\gcd(n^2,l)$. Since $d\mid n$, this implies that $d=\gcd(n^2,l)=\gcd(n,l)$. Hence by \lemref{lem:hec3}, $\omega_\fc = n^2/d^2$. Therefore for each $d\mid n$
$$
\#\{\fc\in\Omega_n^{\rm sim} : \omega_\fc=n^2/d^2 \} = |J_n\backslash Z_d|\ =\ \frac{\varphi(n^2/d)\varphi(d)\gcd(n^2/d,d)}{2n} = \frac{\varphi(n)\varphi(d)}{2},
$$
where for the last equality we used the identities $\gcd(n^2/d,d)=d$ (since $d\mid n$) and
$$\varphi\left(\tfrac{n^2}{d}\right)=\frac{n^2}{d}\prod_{\substack{p\mid (n^2/d)\\ prime}}(1-p^{-1})=\frac{n}{d}\times n\prod_{\substack{p\mid n\\ prime}}(1-p^{-1})=\frac{n\varphi(n)}{d},$$ 
where for the second equality we used the fact that $n^2/d$ and $n$ share the same set of prime divisors.
%First we note here that for any cusp $\fc=m/l\in\Omega_{\G_n}^{\rm sim}$ of simple type with $\gcd(m,l)=1$, similar as in the above proof we have $\omega_{\fc}=n^2/d^2$ with $d:=\gcd(n^2, l)$. On the other hand, by the proof of \lemref{lem:precuspformula} we have for such $d=\gcd(n^2, l)$
%$$\#\{\fc\in \Omega_{\G_n}^{\rm sim} : \omega_{\fc}=n^2/d^2\}=\frac{\varphi(n^2/d)\varphi(d)\gcd(n^2/d,d)}{2n}=\frac{\varphi(n)\times \tfrac{n}{d}\times \varphi(d)\times d}{2n}=\frac{\varphi(n)\varphi(d)}{2},$$
%where for the second equality we used the the identities $\gcd(n^2/d,d)=d$ (since $d\mid n$ by the assumption that $\fc$ is of simple type) and $\varphi(n^2/d)=\varphi(n)\tfrac{n}{d}$ (since $n^2/d$ and $n$ have the same set of prime divisors). 
Hence for $n=m\ell$ we have
\begin{align*}
\#\left\{\fc\in \Omega^{\rm sim}_{n}: \omega_{\fc}\geq m^2\right\}&=\frac{\varphi(n)}{2}\sum_{\substack{d\mid n\\ n^2/d^2\geq m^2}}\varphi(d)\geq \frac{\varphi(n)\varphi(\ell)}{2}=\frac{\varphi(m)(\ell-1)^2}{2}.\qedhere
\end{align*}
%where for the last equality we used the assumptions that $\ell$ is a prime number and $\gcd(\ell, m)=1$.
\end{proof}

\begin{Lem}\label{lem:choiceyn}
Let $m\in\N$ and $Y>0$ satisfy that $m^2Y>1$. Let 
$$\bP_m=\{n=m\ell\in\N: \textrm{$\ell$ is a prime number and $\ell\nmid m$}\}$$
be as in \eqref{equ:primesub}. Then there exist sequences of positive numbers $\{Y_n\}_{n\in\bP_m}$ and $\{Y_n'\}_{n\in\bP_m}$ satisfying that \begin{enumerate}
\item $Y_n'>Y>Y_n>m^{-2}$ for any $n\in\bP_m$ and $\lim\limits_{\substack{n\in\bP_m\\ n\to\infty}}Y_n=\lim\limits_{\substack{n\in\bP_m\\ n\to\infty}}Y'_n=Y$;
\item $\sum_{n\in\bP_m}\frac{1}{n}\left(\frac{1}{Y_n}-\frac{1}{Y_n'}\right)=\infty$.
\end{enumerate}
\end{Lem}
\begin{proof}
For each $n=m\ell\in\bP_m$, take $Y'_n:=(1-(2t_n)^{-1})^{-1}Y$ and $Y_n:=(1+(2t_n)^{-1})^{-1}Y$ with
$$t_n=\max\{(m^2Y-1)^{-1}, \log\log \ell\}.$$
We note that the first condition is guaranteed by the facts that $t_n\geq (m^2Y-1)^{-1}$ and that $\lim\limits_{\substack{n\in\bP_m\\ n\to\infty}}t_n=\infty$. For the second condition, we note that by the definitions of $Y_n$ and $Y_n'$, $\frac{1}{Y_n}-\frac{1}{Y_n'}=\frac{1}{Yt_n}$. Moreover, using the fact that there are only finitely many prime numbers dividing $m$ we get
$$\sum_{n\in\bP_m}\frac{1}{n}\left(\frac{1}{Y_n}-\frac{1}{Y_n'}\right)\asymp_{m,Y}\sum_{\substack{\ell\in \bP_1\\ \ell\nmid m}}\frac{1}{\ell\log\log \ell}=\sum_{\ell\in \bP_1}\frac{1}{\ell\log\log \ell}+O_m(1)=\infty,$$
where the divergence of the rightmost series follows from the estimate $\ell_j\asymp j\log j$ which is an easy consequence of the prime number theorem. Here $\ell_j\in \bP_1$ denotes the $j$-th prime number.
%$Y_j'>Y>Y_j>m^{-2}$ and
%\begin{equation}\label{equ:relationY}
%Y_j^{-1}-(Y_j')^{-1}=\ell_j^{-1}Y^{-1}.
%\end{equation} 
\end{proof}
We now give the 
\begin{proof}[Proof of \thmref{thm:compactsupport}]
Fix throughout the proof $m\in\N$ and $Y>0$ with $m^2Y>1$ and let $\bP_m$ be as above. Let $\{Y_n\}_{n\in\bP_{m}}$ and $\{Y_n'\}_{n\in\bP_m}$ be two sequences satisfy the conditions in \lemref{lem:choiceyn}.
%and let $\bP_m\subset \N$ be as in the theorem. For each $j\in\N$ let $p_j\in\bP_1$ be the $j$-th prime number not dividing $m$. Thus $n_j:=mp_j$ is the $j$-th number in $\bP_m$. Let $\{\ell_j\}_{j\in\N}$ be a non-decreasing sequence of positive numbers satisfying that $\ell_j> \tfrac{1}{2}(m^2Y-1)^{-1}$ for all $j\in\N$, $\lim\limits_{j\to\infty}\ell_j=\infty$ and the series $\sum_{j\in\N}\frac{1}{p_j\ell_j}$ diverges. For instance, using the estimate $p_j\ll j\log j$ (which is an easy consequence of the prime number theorem) we can take 
%$$\ell_j=\max\{(m^2Y-1)^{-1}, \log\log (j+2)\}.$$ 
%For any $j\in\N$ define $Y'_j:=(1-(2\ell_j)^{-1})^{-1}Y$ and $Y_j:=(1+(2\ell_j)^{-1})^{-1}Y$ such that
%$Y_j'>Y>Y_j>m^{-2}$ and
%\begin{equation}\label{equ:relationY}
%Y_j^{-1}-(Y_j')^{-1}=\ell_j^{-1}Y^{-1}.
%\end{equation} 
%Here the inequality $Y_j>m^{-2}$ is guaranteed by our choice of $\{\ell_j\}_{j\in\N}$ such that $\ell_j>\tfrac{1}{2}(m^2Y-1)^{-1}$. 
For any $n\in\bP_m$, let $\Psi_n\in L^2(\G_n\bk \bH)$ such that $\Psi_n$ is the indicator function of the union 
$$\bigcup_{\substack{\fc\in \Omega_{n}^{\rm sim}\\ \omega_{\fc}\geq m^2}}\cC_{\omega_{\fc}Y_n, \omega_{\fc}Y'_n}^{n,\fc}\subset \G_{n}\bk \bH.$$
Since $Y_n>m^{-2}$ for any $n\in\bP_m$, $\omega_{\fc}Y_n>1$ for any $\fc\in \Omega_{n}^{\rm sim}$ with $\omega_{\fc}\geq m^2$. Hence {{similar as in the proof of \thmref{thm:negative2}, by \lemref{lem:onetoone} and \lemref{lem:disjoint} the above union is disjoint and $\Psi_n$ is the indicator function of a Borel set with boundary of measure zero. By the disjointness and}} the volume formula \eqref{equ:volume} we have for any $n\in\bP_m$   
%We note that since we only assume $Y>0$, $\supp \Psi_j$ is not necessarily a disjoint union. However, by our assumption that $m^2Y>1$ and the fact that $\lim\limits_{j\to\infty}Y_j=Y$ (since $\lim\limits_{j\to\infty}\ell_j=\infty$) we have $m^2Y_j>1$ for all $j\gg 1$ sufficiently large. Hence for all $j\gg 1$ sufficiently large, in view of \lemref{lem:disjoint} we have the following subset of $\supp \Psi_j$ which is a disjoint union:
%$$\bigsqcup_{\substack{\fc\in \Omega_{\G_{n_j}}^{\rm sim}\\ \omega_{\fc\geq m^2}}}\cC_{\omega_{\fc}Y_j, \omega_{\fc}Y'_j}^{n_j,\fc}\subset \supp \Psi_j.$$
%This, together with the volume formula \eqref{equ:volume} implies that for all $j\gg 1$ sufficiently large
$$\mu_{\G_{n}}\left(\Psi_n\right)= \frac{3}{\pi}\frac{\#\left\{\fc\in \Omega^{\rm sim}_{n}: \omega_{\fc}\geq m^2\right\}}{[\G_1: \G_{n}]}\left(\frac{1}{Y_n}-\frac{1}{Y_n'}\right).$$
Note that for $n=m\ell\in \bP_m$, by \lemref{lem:indexformula}, $[\G_1: \G_{n}]\asymp_m \ell^3$. Hence by \lemref{lem:simpletype} and the above relation we get for any $n=m\ell\in\bP_m$
\begin{equation}\label{equ:voluforprime}
\mu_{\G_{n}}\left(\Psi_n\right)\gg_{m,Y} \frac{1}{\ell}\left(\frac{1}{Y_n}-\frac{1}{Y_n'}\right)\asymp_m\frac{1}{n}\left(\frac{1}{Y_n}-\frac{1}{Y_n'}\right).
\end{equation}
Similar as in the proof of \thmref{thm:negative2} for any $n\in\bP_m$ and $0<y<1$ we define 
$$I_n(y):=\left\{x\in \R/\Z:\Psi_n(x+iy)=1\right\}.$$
%such that 
%$$\left|\cI_j(y)\right|=\int_0^1\Psi_j(x+iy)dx$$
%By \lemref{lem:noncuspidal} for any $x\in I_j(y)$, all but possibly one points in $\cR_{p_j}(x,y)$ are contained the region $\cC_{Y_j, Y_j'}\subset \cM$. Since $\lim\limits_{j\to\infty}Y_j'=\lim\limits_{j\to\infty}Y_j=Y$, the region $\cC_{Y_j, Y_j'}$ converges to the closed horocycle $\cH_Y$ as $j$ goes to infinity. 
We first show that there exists a sequence $\{y_n\}_{n\in\bP_m}$ satisfying that $0<y_n<c_n$ for all $n\in\bP_m$ and that the limsup set $\limsup_{\substack{n\in\bP_m\\ n\to\infty}}I_n(y_{n})\subset \R/\Z$ is of full measure. As in the proof of \thmref{thm:negative2}, we can use \rmkref{rmk:apparg}, together with \rmkref{rmk:disopen} and \lemref{lem:onetoone}, to construct a sequence $\{y_{n}\}_{n\in\bP_m}$ successively satisfying for any $n\in\bP_m$, $0<y_n<c_n$ and that
\begin{equation}\label{equ:conscondprime}
\left|\frac{1}{\left|I\right|}\int_{I}\Psi_n(x+iy_{n})dx-\mu_{\G_{n}}(\Psi_n)\right|\leq \frac{\mu_{\G_{n}}(\Psi_n)
}{2n^2}
\end{equation}
for all subsets $I\subset \R/\Z$ taken from the finite set $\left\{(0,1)\right\}\bigcup\left\{I_{l}(y_l): l\in \bP_m,\  l<n\right\}$. Again as before one can show that condition \eqref{equ:conscondprime} implies that the sequence $\{I_n(y_{n})\}_{n\in\bP_m}\subset \R/\Z$ satisfies the quasi-independence condition \eqref{equ:quasicond2} (with the subset $\mathbb{S}=\bP_m$ and exponent $\eta=2$). Moreover, using the estimate \eqref{equ:voluforprime} and our assumptions on $\{Y_n\}_{n\in\bP_m}$ and $\{Y_n'\}_{n\in\bP_m}$) we have
$$\sum_{n\in\bP_m}\left|I_n(y_{n})\right|\asymp \sum_{n\in\bP_m}\mu_{\G_{n}}\left(\Psi_n\right)\gg_{m,Y} \sum_{n\in\bP_m}\frac{1}{n}\left(\frac{1}{Y_n}-\frac{1}{Y_n'}\right)=\infty.$$
Hence by \corref{cor:quanbc}, $\limsup_{\substack{n\in\bP_m\\ n\to\infty}}I_n(y_{n})\subset \R/\Z$ is of full Lebesgue measure. 

Now take $x\in \limsup_{\substack{n\in\bP_m\\ n\to\infty}}I_n(y_{n})$, then there exists an unbounded subsequence $\cN_x\subset \bP_m$ such that $x\in I_{n}(y_{n})$ for all $n\in \cN_x$. It thus suffices to show that for any $\Psi\in C_c^{\infty}(\cM)$, 
$$\lim\limits_{\substack{n\in \cN_x\\ n\to\infty}}\delta_{n,x,y_{n}}(\Psi)=\nu_{m,Y}(\Psi)$$
with $\nu_{m,Y}$ defined as in \eqref{equ:limmea}.
%Let $\nu_x$ be a limiting measure of $\{\delta_{p_j,x,y_{p_j}}\}_{j\in S_x}$ and passing to a subsequence if necessary, we may assume that $\delta_{p_j,x,y_{p_j}}$ weakly converges to $\nu_x$ along $S_x$. It thus remains to show that $\nu_x$ is the probability uniformly distributed measure on $\cH_Y$, that is, $\nu_x(\Psi)=\int_0^1\Psi(x+iY)dx$ for any $\Psi\in C_c^{\infty}(\cM)$. 
For any $n\in \cN_x\subset \bP_m$, since $x\in I_n(y_n)$ by definition we have $\G_{n}(x+iy_n)\in \cC_{\omega_{\fc}Y_n, \omega_{\fc}Y'_n}^{n,\fc}$ for some $\fc\in \Omega_{n}^{\rm sim}$ of simple type, that is, there exist  some $\fc\in \Omega_{n}^{\rm sim}$ and $z_n'=x'_n+i\omega_{\fc}y'_n\in \bH$ satisfying that $\G_{n}(x+iy_{n})=\G_{n}\tau_{\fc} z'_n$ with $Y_n< y'_n<Y_n'$. Then by \propref{prop:symmetrysimpletype}, we have
$$\cR_{n}(x,y_{n})=\bigcup_{d\mid n}\cR_{n/d}^{\rm pr}(x'_{n,\fc, d}, d^2y'_n)$$
for some $x'_{n,\fc,d}\in \R/\Z$. This implies that for any $n\in \cN_x$
%the set $\cR_{p_j}(x, y_{p_j})$ differs from the set $\cR_{p_j}(x_j, t_j)$ only by one element, implying that
$$\delta_{n,x, y_{n}}(\Psi)=\frac{1}{n}\sum_{d\mid n}\varphi\left(\tfrac{n}{d}\right)\delta^{\rm pr}_{n/d, x'_{n,\fc,d}, d^2y'_n}(\Psi).$$
Since $y'_n, Y\in (Y_n , Y_n')$, $\max\{y'_n/Y, Y/y'_n\}\leq Y_n'/Y_n$. Thus by the intermediate value theorem we can estimate for $n\in \cN_x$
\begin{align*}
\delta_{n,x, y_{n}}(\Psi)&=\frac{1}{n}\sum_{d\mid n}\varphi\left(\tfrac{n}{d}\right)\left(\delta^{\rm pr}_{n/d, x'_{n,\fc,d}, d^2Y}(\Psi)+O\left(\cS^{\G_1}_{\infty,1}(\Psi)\log\left(Y_n'/Y_n\right)\right)\right)\\
&=\frac{1}{n}\sum_{d\mid n}\varphi\left(\tfrac{n}{d}\right)\delta^{\rm pr}_{n/d, x'_{n,\fc,d}, d^2Y}(\Psi)+O_{\Psi}\left(\log\left(Y_n'/Y_n\right)\right),
%\delta_{p_j, x_j, Y}(\Psi)+O\left(\cS^{\G_1}_{\infty,1}(\Psi)\log\left(Y_j'/Y_j\right)\right). 
%&=\frac{1}{p_j}\sum_{k=0}^{p_j-1}\Psi(x_j+k/p_j+iY)+O\left(\cS_{\infty,1}(\Psi)\log\left(Y_j'/Y_j\right)\right)
\end{align*}
where for the second estimate we used the identity $\sum_{d\mid n}\varphi(n/d)=n$.
Thus for $n=m\ell\in \cN_x$ sufficiently large such that $\Psi$ vanishes on the cusp neighborhood $\cC_{\ell^2Y}$ we have
\begin{align*}
\delta_{n,x, y_{n}}(\Psi)&=\frac{1}{m\ell}\sum_{d\mid m}\varphi\left(\tfrac{m\ell}{d}\right)\delta^{\rm pr}_{m\ell/d, x'_{n,\fc,d}, d^2Y}(\Psi)+O_{\Psi}\left(\log\left(Y_n'/Y_n\right)\right)\\
%&=\frac{\ell-1}{m\ell}\sum_{d\mid m}\varphi\left(\tfrac{m}{d}\right)\delta^{\rm pr}_{m\ell/d, x'_{n,\fc,d}, d^2Y}(\Psi)+O_{\Psi}\left(\log\left(Y_n'/Y_n\right)\right),
&=\frac{\ell-1}{m\ell}\sum_{d\mid m}\varphi\left(\tfrac{m}{d}\right)\left(\mu_{d^2Y}(\Psi)+O_{\Psi, m,Y,\e}\left(\ell^{-1+\e}\right)\right)+O_{\Psi}\left(\log\left(Y_n'/Y_n\right)\right)\\
&=\frac{\ell-1}{\ell}\nu_{m,Y}(\Psi)+O_{\Psi, m,Y,\e}\left(\ell^{-1+\e}+\log\left(Y_n'/Y_n\right)\right),
\end{align*}
where for the second equality we used the facts that %identity that 
%$$\varphi(mp_j/d)=\varphi(m/d)\varphi(p_j)=(p_j-1)\varphi(m/d)$$ 
$\ell$ is a prime number and $\gcd(m,\ell)=1$ and applied the effective estimate \eqref{equ:prmeffbound} to each of the term $\delta^{\rm pr}_{m\ell/d, x'_{n,\fc,d}, d^2Y}(\Psi)$.
%\begin{align*}
%\delta_{n,x, y_{n}}(\Psi)&=\frac{\ell-1}{m\ell}\sum_{d\mid m}\varphi\left(\tfrac{m}{d}\right)\left(\mu_{d^2Y}(\Psi)+O_{\Psi, m,Y,\e}\left(\ell^{-1+\e}\right)\right)+O_{\Psi}\left(\log\left(Y_n'/Y_n\right)\right)\\
%&=\frac{\ell-1}{\ell}\nu_{m,Y}(\Psi)+O_{\Psi, m,Y,\e}\left(\ell^{-1+\e}+\log\left(Y_n'/Y_n\right)\right).
%\end{align*}
%Now apply the definition of $\delta_{p_j, x_j, Y}$, the left $u_1$-invariance of $\Psi$ and the mean value theorem again to get
%$$\delta_{p_j, x_j, Y}(\Psi)=\frac{1}{p_j}\sum_{k=0}^{p_j-1}\Psi(x_j+k/p_j+iY)=\int_0^1\Psi(x+iY)dx+O_Y\left(p_j^{-1}\cS^{\G_1}_{\infty,1}(\Psi)\right).$$
%To conclude we have
%\begin{align*}
%\delta_{p_j,x, y_{p_j}}(\Psi)=\int_0^1\Psi(x+iY)dx+O_Y\left(\cS^{\G_1}_{\infty,1}\left(\Psi)(p_j^{-1}+\log\left(Y_j'/Y_j\right)\right)\right),
%\end{align*}
%where the term $p_j^{-1}\|\Psi\|_{\infty}$ is absorbed into $p_j^{-1}\cS_{\infty,1}(\Psi)$. 
We now conclude by taking $n\to\infty$ along the subsequence $\cN_x$ and noting that $\lim\limits_{\substack{n\in \cN_x\\ n\to\infty}}\log\left(Y_n'/Y_n\right)=0$ (since $\lim\limits_{\substack{n\in\bP_m\\ n\to\infty}}Y_n'/Y_n=1$ which follows from the assumption $\lim\limits_{\substack{n\in\bP_m\\ n\to\infty}}Y_n=\lim\limits_{\substack{n\in\bP_m\\ n\to\infty}}Y'_n=Y$).
\end{proof}

\begin{rmk}\label{rmk:combine}
It is clear that we can take a sequence $\{y_n\}_{n\in\N}$ decaying sufficiently fast such that the conditions \eqref{equ:conscond} and \eqref{equ:conscondprime} (for any finitely many pairs $(m,Y)$ with $m^2Y>1$) are all satisfied and hence (noting that the intersection of finitely many full measure sets is still of full measure) for such a sequence the conclusions of \thmref{thm:negative2} and \thmref{thm:compactsupport} (for any finitely many pairs $(m,Y)$ with $m^2Y>1$) hold simultaneously.
\end{rmk}

\bibliographystyle{alpha}
\bibliography{DKbibliog}

\begin{thebibliography}{EBHL18}

\bibitem[Ber20]{Bersudsky2020}
M.~Bersudsky.
\newblock On the image in the torus of sparse points on dilating analytic
  curves.
\newblock {\em arXiv preprint arXiv:2003.04112}, 2020.

\bibitem[Bor72]{Borel1972}
A.~Borel.
\newblock Some metric properties of arithmetic quotients of symmetric spaces
  and an extension theorem.
\newblock {\em J. Differential Geometry}, 6:543--560, 1972.

\bibitem[COU01]{ClozelOhUllmo2001}
L.~Clozel, H.~Oh, and E.~Ullmo.
\newblock Hecke operators and equidistribution of {H}ecke points.
\newblock {\em Invent. Math.}, 144(2):327--351, 2001.

\bibitem[CU04]{ClozelUllmo2004}
L.~Clozel and E.~Ullmo.
\newblock Equidistribution des pointes de hecke.
\newblock In {\em Contribution to automorphic forms, geometry and number
  theory}, pages 193--254. Johns Hopkins Univ. Press, 2004.

\bibitem[DA14]{Demirci2014}
E.~Demirci~Akarsu.
\newblock Short incomplete {G}auss sums and rational points on metaplectic
  horocycles.
\newblock {\em Int. J. Number Theory}, 10(6):1553--1576, 2014.

\bibitem[DS05]{DiamondShurman2005}
F.~Diamond and J.~Shurman.
\newblock {\em A first course in modular forms}, volume 228 of {\em Graduate
  Texts in Mathematics}.
\newblock Springer-Verlag, New York, 2005.

\bibitem[EBHL18]{BazHuangLee2019}
D.~El-Baz, B.~Huang, and M.~Lee.
\newblock Effective joint equidistribution of primitive rational points on
  expanding horospheres.
\newblock {\em arXiv preprint arXiv:1811.04019}, 2018.

\bibitem[ELS20]{EinsiedlerLuethiShah2020}
M.~Einsiedler, M.~Luethi, and N.~Shah.
\newblock Primitive rational points on expanding horocycles in products of the
  modular surface with the torus.
\newblock {\em Ergodic Theory and Dynamical Systems}, pages 1--45, 2020.

\bibitem[EMSS16]{EinsiedlerMozesShahShapira2016}
M.~Einsiedler, S.~Mozes, N.~Shah, and U.~Shapira.
\newblock Equidistribution of primitive rational points on expanding
  horospheres.
\newblock {\em Compos. Math.}, 152(4):667--692, 2016.

\bibitem[GM03]{GoldsteinMayer2003}
D.~Goldstein and A.~Mayer.
\newblock On the equidistribution of {H}ecke points.
\newblock {\em Forum Math.}, 15(2):165--189, 2003.

\bibitem[GR70]{GarlandRaghunathan1970}
H.~Garland and M.~S. Raghunathan.
\newblock Fundamental domains for lattices in ({R}-)rank {$1$} semisimple {L}ie
  groups.
\newblock {\em Ann. of Math. (2)}, 92:279--326, 1970.

\bibitem[Hej96]{Hejhal1996}
D.~A. Hejhal.
\newblock On value distribution properties of automorphic functions along
  closed horocycles.
\newblock In {\em X{VI}th {R}olf {N}evanlinna {C}olloquium ({J}oensuu, 1995)},
  pages 39--52. de Gruyter, Berlin, 1996.

\bibitem[Hej00]{Hejhal2000a}
D.~A. Hejhal.
\newblock On the uniform equidistribution of long closed horocycles.
\newblock {\em Asian J. Math.}, 4(4):839--853, 2000.
\newblock Loo-Keng Hua: a great mathematician of the twentieth century.

\bibitem[HW08]{HardyWright2008}
G.~H. Hardy and E.~M. Wright.
\newblock {\em An introduction to the theory of numbers}.
\newblock Oxford University Press, Oxford, sixth edition, 2008.
\newblock Revised by D. R. Heath-Brown and J. H. Silverman, With a foreword by
  Andrew Wiles.

\bibitem[Iwa02]{Iwaniec2002}
H.~Iwaniec.
\newblock {\em Spectral methods of automorphic forms}, volume~53 of {\em
  Graduate Studies in Mathematics}.
\newblock American Mathematical Society, Providence, RI; Revista Matem\'{a}tica
  Iberoamericana, Madrid, second edition, 2002.

\bibitem[Jan21]{Jana2021}
S.~Jana.
\newblock Joint equidistribution on the product of the circle and the unit
  cotangent bundle of the modular surface.
\newblock {\em J. Number Theory}, 226:271--283, 2021.

\bibitem[KK18]{KelmerKontorovich2018}
D.~Kelmer and A.~Kontorovich.
\newblock Effective equidistribution of shears and applications.
\newblock {\em Math. Ann.}, 370(1-2):381--421, 2018.

\bibitem[KK20]{KelmerKontorovich2020}
D.~Kelmer and A.~Kontorovich.
\newblock Exponents for the equidistribution of shears and applications.
\newblock {\em J. Number Theory}, 208:1--46, 2020.

\bibitem[KS03]{KimSarnak03}
H.~H. Kim and P.~Sarnak.
\newblock Refined estimates towards the {R}amanujan and {S}elberg conjectures,
  appendix to {H}. {H}. {K}im, {F}unctoriality for the exterior square of
  {$GL_4$} and the symmetric fourth of {$GL_2$}, with appendix 1 by {D}inakar
  {R}amakrishnan and appendix 2 by {K}im and {P}eter {S}arnak.
\newblock {\em J. Amer. Math. Soc.}, 16(1), 2003.

\bibitem[KY19]{KelmerYu2019a}
D.~Kelmer and S.~Yu.
\newblock Shrinking target problems for flows on homogeneous spaces.
\newblock {\em Trans. Amer. Math. Soc.}, 372(9):6283--6314, 2019.

\bibitem[Lan75]{LangSL2R}
S.~Lang.
\newblock {\em {${\rm SL}\sb{2}({\bf R})$}}.
\newblock Addison-Wesley Publishing Co., Reading, Mass.-London-Amsterdam, 1975.

\bibitem[Li15]{Han2015}
H.~Li.
\newblock Effective limit distribution of the {F}robenius numbers.
\newblock {\em Compos. Math.}, 151(5):898--916, 2015.

\bibitem[LM18]{LeeMarklof2018}
M.~Lee and J.~Marklof.
\newblock Effective equidistribution of rational points on expanding
  horospheres.
\newblock {\em Int. Math. Res. Not. IMRN}, (21):6581--6610, 2018.

\bibitem[Lue21]{Luethi2021}
M.~Luethi.
\newblock Primitive rational points on expanding horospheres in {H}ilbert
  modular surfaces.
\newblock {\em J. Number Theory}, 225:327--359, 2021.

\bibitem[Mar04]{Margulis2004}
G.~A. Margulis.
\newblock {\em On some aspects of the theory of {A}nosov systems}.
\newblock Springer Monographs in Mathematics. Springer-Verlag, Berlin, 2004.
\newblock With a survey by Richard Sharp: Periodic orbits of hyperbolic flows,
  Translated from the Russian by Valentina Vladimirovna Szulikowska.

\bibitem[Mar10]{Marklof2010}
J.~Marklof.
\newblock The asymptotic distribution of {F}robenius numbers.
\newblock {\em Invent. Math.}, 181(1):179--207, 2010.

\bibitem[MS03]{MarklofStrombergsson2003}
J.~Marklof and A.~Str\"{o}mbergsson.
\newblock Equidistribution of {K}ronecker sequences along closed horocycles.
\newblock {\em Geom. Funct. Anal.}, 13(6):1239--1280, 2003.

\bibitem[Sar81]{Sarnak1981}
P.~Sarnak.
\newblock Asymptotic behavior of periodic orbits of the horocycle flow and
  {E}isenstein series.
\newblock {\em Comm. Pure Appl. Math.}, 34(6):719--739, 1981.

\bibitem[Spr79]{Sp79}
V.~G. Sprind{\v z}uk.
\newblock {\em Metric theory of {D}iophantine approximations}.
\newblock V. H. Winston \&\ Sons, Washington, D.C.; A Halsted Press Book, John
  Wiley \&\ Sons, New York-Toronto, Ont.-London, 1979.
\newblock Translated from the Russian and edited by Richard A. Silverman, With
  a foreword by Donald J. Newman, Scripta Series in Mathematics.

\bibitem[Str04]{Strombergsson2004}
A.~Str\"{o}mbergsson.
\newblock On the uniform equidistribution of long closed horocycles.
\newblock {\em Duke Math. J.}, 123(3):507--547, 2004.

\bibitem[SU15]{SarnakUbis2015}
P.~Sarnak and A.~Ubis.
\newblock The horocycle flow at prime times.
\newblock {\em J. Math. Pures Appl. (9)}, 103(2):575--618, 2015.

\bibitem[Zag81]{Zagier1981}
D.~Zagier.
\newblock Eisenstein series and the {R}iemann zeta function.
\newblock In {\em Automorphic forms, representation theory and arithmetic
  ({B}ombay, 1979)}, volume~10 of {\em Tata Inst. Fund. Res. Studies in Math.},
  pages 275--301. Tata Inst. Fundamental Res., Bombay, 1981.

\end{thebibliography}
\end{document}